\documentclass[11pt,titlepage]{article}
\usepackage{amscd,amsmath,amsxtra,amssymb,theorem,latexsym,amsfonts}
\usepackage[all]{xy}
\usepackage{epic,eepic}

\usepackage{oldgerm}
\usepackage{bbm}

\def\hom{\text{Hom}}

\def\ext{\text{Ext}}

\def\ker{\text{Ker}}

\def\coker{\text{Coker}}
\def\supp{\text{Supp}}

\def\rank{\text{rank}}
\def\Span{\text{span}}
\def\length{\text{length}}
\def\codim{\text{codim}}

\newtheorem{satz}{Satz}[section]
\newtheorem{lemma}[satz]{Lemma}
\newtheorem{theorem}[satz]{Theorem}
\newtheorem{defin}[satz]{Definition}
\newtheorem{corollary}[satz]{Corollary}
\newtheorem{prop}[satz]{Proposition}

\newtheorem{remark}[satz]{Remark}
\newtheorem{construction}[satz]{Construction}
\newtheorem{thought}[satz]{}

\newcommand{\Mgot}{{\mathfrak{M}}}
\renewcommand{\sl}{{\mathfrak{sl}(2)}}
\renewcommand{\gg}{{\mathfrak{g}}}
\newcommand{\gh}{{\mathfrak{h}}}
\newcommand{\gf}{{\mathfrak{f}}}

\renewcommand{\P}{{\mathbb P}}

\newcommand{\G}{{\mathbb G}}

\newcommand{\Id}{{\mathbbm 1}}

\newcommand{\ka}{{\mathcal A}}
\newcommand{\kb}{{\mathcal B}}
\newcommand{\kc}{{\mathcal C}}

\newcommand{\ke}{{\mathcal E}}
\newcommand{\kf}{{\mathcal F}}
\newcommand{\kg}{{\mathcal G}}
\newcommand{\kh}{{\mathcal H}}
\newcommand{\ki}{{\mathcal I}}

\newcommand{\kk}{{\mathcal K}}
\newcommand{\kl}{{\mathcal L}}
\newcommand{\km}{{\mathcal M}}
\newcommand{\kn}{{\mathcal N}}
\newcommand{\ko}{{\mathcal O}}

\newcommand{\kt}{{\mathcal T}}
\newcommand{\ku}{{\mathcal U}}

\newcommand{\kx}{{\mathcal X}}

\newcommand{\bkf}{{\boldsymbol{\mathcal F}}}
\newcommand{\kbf}{{\boldsymbol{\mathcal F}}}

\newcommand{\bkn}{{\boldsymbol{\mathcal N}}}
\newcommand{\bke}{{\boldsymbol{\mathcal E}}}

\newcommand{\qed}{{\begin{flushright}$\Box$\end{flushright}}}

\newcommand{\Fgot}{{\mathfrak F}}

\CompileMatrices
\SelectTips{cm}{11}

\begin{document}



\thispagestyle{empty}

\begin{center}
{\huge Universit\"at Kaiserslautern} \\[1ex]
{\huge Fachbereich Mathematik} \\[30ex]
{\Huge\bf  Diploma Thesis} \\[10ex]
{\Large by Alexander Getmanenko}\\[10ex]
{\Huge  Geometry and Moduli Space of Certain \\
       Rank-4 Vector Bundles on $\P^4$} \\[10ex]
{\Large under the Supervision of\\
Prof. Dr. G. Trautmann} \\[17ex]
June, 2000
\end{center}

\newpage

\pagenumbering{roman} 
\vspace*{3cm}
\begin{center}
{\Large \bf Acknowledgements}
\end{center}

\vspace{3cm}
{\it 
I thank Prof.Dr.G.Trautmann, my supervisor for this Diploma
thesis, for the interesting problem, for many hours of discussions and for
deep mathematical ideas I could learn form him.

I express my gratitude to Dr.B.Kreu\ss ler who took a lot of care
of me and my education in the field of algebraic geometry, and to all members
and students of the  Group for Algebraic Geometry at the University 
Kaiserslautern for their sympathy and helpfulness, mathematical and
personal. 

My studies in Germany were made possible by the scholarship I was granted by
the Faculty of Mathematics of the University Kaiserslautern as a participant of 
the program "Mathematics International".

}

\newpage

\vspace*{7cm}

\tableofcontents

\newpage

\pagenumbering{arabic}



\section{Introduction} 
\label{first}

The topic of this diploma thesis are instanton bundles of rank 4 with the 
second Chern class (quantum number) 2 on $\P^4$ which are, as will be shown in 
the section \ref{intro}, precisely the bundles appearing as cohomologies of
short monads
$$ 
     2\Omega_{\P^4}^4(4) \stackrel{M}{\longrightarrow} 2\Omega_{\P^4}^1(1)
      \stackrel{N}{\longrightarrow} 2\ko_{\P^4} 
$$
 
We use the following general definition of instanton bundles: 

{\bf Definition \ref{Definst}} {\it 
A mathematical instanton bundle on $\P^s$ with quantum number
$n$ or an $n$-instanton is an algebraic vector bundle $\ke$ on $\P^s$ which 
satisfies\\
(i) $\ke$ has Chern polynomial 
$(1-h^2)^{-n}$  and $s-1 \le r = \rank\,\ke\, \le (s-1)n$ \\
(ii) for any $d$ in the range $-s \le d \le 0$ we have $h^i\ke (d) \ne 0$ for 
at most one $i$.
}

In \cite{script} you can find a collection of general results and 
basic cinstructions that are useful to work with bundles with these properties.

The moduli problem for instanton bundles on $\P^3$ is being intensively studied
for different Chern classes. Originally instanton bundles of rank 2 on $\P^3$ 
appeared as liftings of bundles on $S^4$ via the twistor construction, see
\cite{AHDM}. They have been studied in a number of papers, e.g. \cite{NarTrm},
\cite{HaHi}, \cite{HNar}, \cite{Nuessler}, \cite{MaTrm},  \cite{Barth}.
The most recent result in the direction is that in \cite{KO}
where it is proved the regularity of the moduli space 
of mathematical instanton bundles on $\P^3$ with $c_2=5$. Together with the 
moduli problem, different questions about the geometry can be considered, 
e.g. geometry of zero sets of sections and sets of jumping lines, cf. 
\cite{NarTrm},\cite{BHir1}, \cite{BHir2}. 

Instanton bundles of rank $s-1$ on $\P^s$ for odd $s$ have been studied
in \cite{SpindlerTrm}, \cite{OkSpindler}, \cite{AO1}, \cite{OTrm}, \cite{AO2}. 

Another reason making mathematical instanton bundles interesting is, that
they give nontrivial examples of simple vector bundles on
projective spaces of different dimension.

In this diploma thesis I am discussing instanton bundles of rank 4 with 
quantum number 2 on $\P^4$. This case is similar to that considered in 
papers \cite{NarTrm}, \cite{HNar} on $\P^3$ and to special 't Hooft bundles on 
$\P^5$ in \cite{SpindlerTrm}.

The main results of my diploma thesis are:

{\bf Theorem \ref{ExistModuli}\ }{\it
There exixt a 29-dimensional irreducible quasi-projective 
variety $\km$ that is a coarse moduli space 
of instanton bundles on $\P^4$ with quantum number 2 and of 
rank 4.}\\

{\bf Theorem \ref{Msmooth}\ }{\it The moduli space $\km$ is smooth.}\\

{\bf Theorem \ref{DualTh}\ }{\it
Let an instanton bundle $\ke$ of rank 4 with quantum number 2 on $\P^4=\P V$ 
($V$ is a 5-dimensional vector space)
be given by a monad
$$ 
     2\Omega_{\P^4}^4(4) \stackrel{M}{\longrightarrow} 2\Omega_{\P^4}^1(1)
      \stackrel{N}{\longrightarrow} 2\ko_{\P^4} 
$$
and let the matrix $P$ be defined as the kernel in
$$ 0 \to k^p \stackrel{P}{\longrightarrow}  
    k^2\otimes V \stackrel{\tilde M}{\longrightarrow} 
    k^2\otimes \bigwedge^4 V,
$$
where $\tilde M$ is the contraction with the matrix $M$ and can be proved
to be of rank 7 or 8.
Then the dual bundle $\ke^\vee$ is either 
the cohomology bundle of a monad
$$ 
2\Omega^4(4) \stackrel{M^\vee}{\longrightarrow} 2\Omega^1(1)
      \stackrel{P^\vee}{\longrightarrow} 2\ko 
$$
if $\rank\,\tilde M\,=8$, and is therefore an instanton bundle, too,
or, if $\rank\,\tilde M\,=7$, it is an extension
$$ 
0 \to \ko \to \ke^\vee \to \ke^\vee_1 \to 0
$$
where $\ke_1^\vee$ is a reflexive sheaf appearing as the cohomology of the 
monad
$$ 
2\Omega^4(4) \stackrel{M^\vee}{\longrightarrow} 2\Omega^1(1)
      \stackrel{P^\vee}{\longrightarrow} 3\ko, 
$$
}\\

{\bf Theorem \ref{JumpTh}\ } {\it
In the notation of theorem \ref{DualTh},
let $\ell$ be a line in $\P^4$ identified via the Pl\"ucker embedding with
an element of $\P\bigwedge^2 V \, \supset \, \G (2,5)\, \ni\, \ell$. 
If $\langle \ell \rangle = \langle x \wedge y \rangle$ in $\P\bigwedge^2 V$ and
$M=(m_{i,j})_{i,j=1,2}$, denote by $M\wedge \ell$ a matrix 
$(m_{i,j}\wedge x \wedge y)_{i,j=1,2} \in Mat_{2\times 2}(\bigwedge^5 V) \cong
 Mat_{2\times 2}(k)$ (well defined up to  proportionality). 
Then:\\
(i) if $\rank\, M\wedge\ell = 2$, then $\ell$ is not a jumping line; \\
(ii) if $\rank\, M\wedge\ell = 1$, then $\ell$ is a jumping line with splitting
$$ \ke_\ell \cong 
          \ko_\ell(-1) \oplus \ko_\ell^{\oplus 2} \oplus \ko_\ell(1); $$
(iii) if $ M\wedge\ell = 0$, then $\ell$ is a jumping line with
$$\ke_\ell 
     \cong \ko_\ell(-1)^{\oplus 2} \oplus \kf
          \text{\quad or\quad } 
\ke_\ell\cong\ko_\ell(-2) \oplus \ko_\ell \oplus \kf$$
where $\kf \cong \ko_\ell(1)^{\oplus 2}$ \quad or 
          \quad $\kf\cong\ko_\ell \oplus \ko_\ell(2)$; \\
(iv) in the situation of (iii), the jumping lines $\ell$ with the property 
$\ke_\ell=\ko_\ell(-2)\oplus\ko_\ell\oplus\kf$
form a smooth conic on the Grassmannian $\G (2,5)$; \\
(v) in the situation of (iii), the jumping lines $\ell$ with the property 
$\kf = \ko_\ell\oplus\ko_\ell(2)$ form either a smooth conic, if 
$\rank\,\tilde M\,=8$, or, if $\rank\,\tilde M\,=7$, a surface 
 on the Grassmannian $\G (2,5)$.
}\\

{\bf Theorem \ref{PlaneTh}\ \it In the notation of theorem \ref{DualTh},
denote by $W \, =\, \Span\, N$ the vector subspace of $V$ generated by the 
entries of $N$.
If $\dim\,W=4$ and $H=\P W \cong \P^3$ is a hyperplane in $\P^4$,
then $\ke|_H \cong \kh \oplus 2\ko_H$, where $\kh$ is an instanton 
bundle on $\P^3$  of rank 2 with $c_2=2$.
}\\

In addition to these theorems, the following results are obtained:

$\bullet\ $ computation of dimensions and proof of irreducibility of certain 
subsets of the moduli space that naturally arise in the previous steps 
(section \ref{subsets}), namely:\\
$\km^3\ :=\ $ sheaves in whose monads $\dim\,\Span\,N\,=\,3$  \\   
$\km^7\ :=\ $ sheaves $\ke$ in whose monads $rank\,\tilde M\, =\, 7$\, or,
equivalently (cf. corollary \ref{rkMtilde}), $H^0 \ke^\vee = 1$ \\
$\km^{sd}\ :=\ $ self-dual bundles $\ke$, i.e. bundles with
$\ke \cong \ke^\vee$.\\
In particular, it is shown (theorem \ref{M7thm}) that $\km^7\subset \km^3$.\\

$\bullet\ $ construction of a fibration of the moduli space over an open 
dense subset of the set of singular hyperquadrics in $\P^4$ 
(section \ref{fibration});\\

$\bullet\ $ some remarks on geometry of sections of the kernel bundle of monads
defining instanton bundles on $\P^4$ with $c_2=2$ (section \ref{scroll}).

\newpage

\section{Basics}
\label{intro}


In this preparatory section we are giving some basic definitions and facts 
about mathematical instanton bundles in general we shall use in the sequel and 
also describe the objects we are going to work with in our 
situation.

The main reference here is \cite{script}.

We are working with algebraic schemes over an algebraicly closed field $k$ 
of characteristic 0 and therefore may do most of our proofs for closed points 
only.

Given a projectivization $\P W$ of a $k$-vector space $W$, we denote by 
$\langle w \rangle$ the image of $w\in W$ under the canonical projection 
$W-\{ 0\}\,\to\,\P W$.

Given a $k$-vector space, a sheaf, a module etc. $A$ and an integer $n\ge 0$, 
the symbols $A^{\oplus n}$, $nA$ and $k^n\otimes A$ mean the same; in each 
particular situation we try to use the most convenient notation.

To fix the notation, we take the $s$-dimensional projective space $\P^s$ to
be the projectivization of an $(s+1)$-dimensional vector space $V$ and write
$\P^s = \P V$.

\begin{defin} \label{Definst}  
A mathematical instanton bundle on $\P^s$ with quantum number
$n$ or an $n$-instanton is an algebraic vector bundle $\ke$ on $\P^s$ which 
satisfies\\
(i) $\ke$ has Chern polynomial 
$(1-h^2)^{-n}$  and $s-1 \le r = \rank\,\ke\, \le (s-1)n$ \\
(ii) for any $d$ in the range $-s \le d \le 0$ we have $h^i\ke (d) \ne 0$ for 
at most one $i$.
\end{defin}

\begin{lemma} \label{Hinst} \cite[lemma 1.3]{script}
 Assuming condition (i), conditon (ii) is equivalent to condition \\
$\begin{array}{ccl}
(ii') & (a) & h^i\ke = 0 \ \text{for} \ i\ne 1 \ \text{and} \ 
            h^1\ke =(s-1)n-r  \\
 & (b) & h^i\ke (-1)=0 \ \text{for} \ i\ne 1 \ \text{and} \ h^1\ke (-1)=n \\
 & (c) & h^i\ke (d)=0 \ \text{for} \ -s < d < -1  \ \text{and all} \ i \\
 & (d) & h^i\ke (-s)=0 \ \text{for} \ i\ne s-1 \ \text{and} \ h^{s-1}\ke (-s)=n  
\end{array}$
\end{lemma}

{\bf The objects of this study are rank-4 instanton bundles on $\P^4$ with 
quantum number $2$. }

The conditions can be written in this case as:\\
(i) The Chern polynomial is $c(\ke ) = (1-h^2)^{-2} = 1 + 2 h^2 + 3 h^4$ \\
$\begin{array}{ccl}
(ii') & (a) & h^i\ke = 0 \ \text{for} \ i\ne 1 \ \text{and} \ 
            h^1\ke = 2  \\
 & (b) & h^i\ke (-1)=0 \ \text{for} \ i\ne 1 \ \text{and} \ h^1\ke (-1)=2 \\
 & (c) & h^i\ke (d)=0 \ \text{for}\ -4 < d < -1  \ \text{and all} \ i \\
 & (d) & h^i\ke (-4)=0 \ \text{for}\ i\ne 3 \ \text{and} \ h^3\ke (-4)=2
\end{array}$


There is a convenient way to describe intanton bundles through Beilinson 
monads.

{\bf Beilinson monads }

One of the usual methods to describe sheaves is by constructing some exact 
sequences made up of simplier sheaves such that the required sheaf is the 
kernel, cokernel, etc. of a certain morphism.   
For our problem the building blocks will be the sheaves 
$\Omega^i(j)$ on $\P^s$,
where
$\Omega^i = \bigwedge^i \kt^\vee$
and $\kt$ is a tangent bundle of $\P^s$.
Invariants of the sheaves $\Omega^i(j)$ and their morphisms  are 
classically  known (the so called Bott formulae can be found, for example, in 
\cite[Ch.I.1]{OSS}). We shall represent instanton bundles as 
cohomologies 
of complexes of sheaves called monads and then use this representation to 
introduce  parametrization of the bundles by matrices

\begin{defin}
If $X$ is a scheme, a {\bf monad} over $X$ is a complex
$$
 \dots \kf_{-1} \to \kf_0 \to \kf_1 \to \dots
$$
of coherent sheaves over $X$ which is exact everywhere except at $\kf_0$.
The cohomology $\ke$ at $\kf_0$ is called the {\bf cohomology of} this 
{\bf monad}, and we speak of a {\bf monad for $\ke$ }.
\end{defin}

An intermediate step we shall use to construct a monad for instanton bundles
will be the Beilinson spectral sequence. The theorem of Beilinson proves the 
existence of such a spectral sequence. Because of importance of its 
construction we quote this 
theorem together with the sketch of its proof done in details in 
\cite[Ch.II.3]{OSS}.

\begin{theorem} \label{BeilSpectral}
Let $\ke$ be a vector bundle over $\P^s$. There is a  
{\bf Beilinson spectral sequence of \linebreak type II} $E^{p,q}_{r}$ with 
$E^1$-term 
$$ E^{pq}_1 = H^q(\ke(p))\otimes \Omega^{-p}(-p) $$
which converges to $(E^i)$ with  
$$ E^i = \left\{ \begin{array}{l} 
\ke \ \text{ for \ } i=0, \\
0 \ \text{ otherwise,}
\end{array} \right.
$$
i.e. $E^{p,q}_{\infty}=0$ for $p+q \ne 0$ and 
$\bigoplus_{p=0}^{n} E^{-p,p}_{\infty}$ is the associated graded sheaf of a 
filtration of $E$.
\end{theorem} 

{\bf Sketch of the proof.} Consider the product projections 
$p_1,p_2 : \P^s \times \P^s \to \P^s$. 
We write \linebreak
$\ka \boxtimes \kb = p_1^*\ka \otimes p_2^*\kb$ for any bundles $\ka,\kb$.

There is a free (Koszul) resolution of $\ko_\Delta$ where $\Delta$ is the 
diagonal in $\P^s \times \P^s$:
$$ 
0 \to \bigwedge^s(\ko(-1) \boxtimes \Omega^1(1)) \to 
\bigwedge^{s-1}(\ko(-1) \boxtimes \Omega^1(1)) \to \dots \to
\ko(-1) \boxtimes \Omega^1(1) \to \ko_{\P^s\otimes \P^s} \to \ko_\Delta \to 0
$$
or, written in another form,
$$ 
0 \to \ko(-s) \boxtimes \Omega^s(s) \to 
\ko(-s+1) \boxtimes \Omega^{s-1}(s-1) \to \dots \to
\ko(-1) \boxtimes \Omega^1(1) \to \ko_{\P^s\otimes \P^s} \to \ko_\Delta \to 0.
$$

Let $k(\langle u \rangle,\langle v \rangle ) \, \cong \, k$ be the residue 
field at a point
$(\langle u \rangle,\langle v \rangle )\in \P^s \times \P^s = \P V \times \P V$
where $u,v\in V$. 
We think of $k(\langle u \rangle,\langle v \rangle )$
as of a skyscraper sheaf at the point 
$(\langle u \rangle,\langle v \rangle ) \in \P^s\times\P^s$ and tensorizing 
with it means simply taking a fiber at this point.

The maps of this Koszul complex are defined up to a twist by the following maps
on fibers of the vector bundles 
$(\langle u \rangle,\langle v \rangle ) \in \P^s\times\P^s$
$$\begin{array}{rcl}
 {( \ko(-1) \boxtimes \Omega^j(j) ) 
    \otimes k(\langle u \rangle,\langle v \rangle )} & \longrightarrow &
   {( \ko \boxtimes \Omega^{j-1}(j-1) ) 
    \otimes k(\langle u \rangle,\langle v \rangle )} \\
\bigwedge^jV^\vee \supset \bigwedge^j (V/v)^\vee \ \cong \ 
  (k\cdot u) \otimes \bigwedge^j (V/v)^\vee & \longrightarrow &
  k \otimes \bigwedge^{j-1} (V/v)^\vee \ =\ \bigwedge^{j-1} (V/v)^\vee \\
x \otimes \theta & \longmapsto & 1 \otimes (x \lrcorner \theta ) 
\end{array}$$
for every point $(\langle u \rangle,\langle v \rangle )\in \P^s \times \P^s$
and $x \lrcorner \theta$ standing for the result of contraction
$$
V \otimes \bigwedge^{j}V^\vee \stackrel{\lrcorner}{\longmapsto} 
\bigwedge^{j-1}V^\vee
$$
Here we used the canonical description of the fibers of the vector bundles 
$\ko(-1)$ and $\Omega^j(j)$ over a point $\langle x \rangle \in \P^s$ as
$$\begin{array}{rclcl}
\ko(-1)\otimes k(\langle x \rangle) & = & (k\cdot x)  &\subset  & V \\
\Omega^j(j)\otimes k(\langle x \rangle) & = & \bigwedge^j (V/x)^\vee & 
 \subset & \bigwedge^j V^\vee
\end{array}$$
Here $k\cdot x$ means the subspace of $V$ generated by a vector $x$. 

Write $\kc^{-k}=\ke(-k) \boxtimes \Omega^{k}(k)$, where $0\le k \le s$. The 
Koszul sequence above gives after tensorizing with $p_1^* \ke$ yields a complex
$$ 0 \to \kc^{-s} \to \kc^{-s+1} \to \dots \to \kc^0 \to 0. $$

By the theory of hyperdirect images of sheaves \cite{EGAIII}, there are two 
spectral sequences 
$$ \begin{array}{c}
 'E^{p,q}_{2} = H^p(R^q p_{2*}(\kc^\cdot)) \ , \quad
 {''E}^{p,q}_{2} = R^p p_{2*}(\kc^\cdot))
\end{array}$$
converging to the $(p+q)$-th hyperdirect image sheaf $R^{p+q}p_{2*}(\kc^\cdot)$.
A reader may consult \cite{EGAIII} for definitions of higher direct images and of 
hypercohomologies; they are not that important for this exposition. Important 
is that the two spectral sequences have the same limit.

Since $\kc^\cdot$ is a locally free resolution of $p_1^*\ke |_\Delta$, it 
follows that 
$$ H^q(\kc^\cdot) \ \cong \ 
\left\{ \begin{array}{ll}
\ke & \text{ for $p=q=0$, } \\ 0 & \text{ otherwise } 
\end{array} \right.
$$
and therefore 
$$ {''E}^{p,q}_2 \, = \, R^pp_{2*}(H^q(\kc^\cdot) \, \cong \, 
\left\{ \begin{array}{ll}
\ke & \text{ for $p=q=0$, } \\ 0 & \text{ otherwise } 
\end{array} \right.
$$
The spectral sequence ${''E}^{p,q}_2$ converges to $R^{p+q}p_{2*}(\kc^\cdot)$ 
and thus
$$ R^ip_{2*}(\kc^\cdot) \ \cong \  
\left\{ \begin{array}{ll}
\ke & \text{ for $p=q=0$, } \\ 0 & \text{ otherwise } 
\end{array} \right.
$$
On the other hand, one proves that we can construct the $E_1$-level of 
the spectral sequence $'E_1$ by setting
$$\begin{array}{rcl}
'E^{p,q}_1 = R^qp_{2*}(\kc^p) & = & R^qp_{2*}(\ke(p)\boxtimes \Omega^{-p}(-p) \\
& = & R^qp_{2*}(p_1^* \ke(p)\boxtimes \Omega^{-p}(-p) \\
& = & H^q(\P^s, \ke(p))\otimes \Omega^{-p}(-p), 
\end{array}
$$
where $0\le q\le s$ and $-s \le p \le 0$, with maps 
$'E^{p,q}_1 \to 'E^{p+1,q}_1$ induced by maps in the Koszul resolution above 
and described as
$$\begin{array}{rclcl} 'E^{p,q}_1 & = & 
H^q(\P^s, \ke(p))\otimes \Omega^{-p}(-p) \\
& \to &  H^q(\P^s, \ke(p))\otimes V^\vee \otimes \Omega^{-p-1}(-p-1) \\
& \cong &
 H^q(\P^s, \ke(p))\otimes H^0(\P^s, \ko(1))\otimes \Omega^{-p-1}(-p-1) \\ 
& \to & H^q(\P^s, \ke(p+1))\otimes \Omega^{-p-1}(-p-1) & = & 'E^{p+1,q}_1
\end{array}$$
Here the first map is induced by the contraction homomorphism
$\Omega^{-p}(-p) \otimes V \to \Omega^{-p-1}(-p-1)$, the isomorphism
$H^0(\P^s, \ko(1))\cong V^\vee$ is canonical and the last map is the 
cup-product in cohomologies.
\qed

There is also another variant of the Beilinson theorem:
\begin{theorem} Let $\ke$ be a vector bundle over $\P^s$. There is a 
{\bf Beilinson spectral sequence of type I} $E^{p,q}_{r}$ with $E^1$-term 
$$ E^{pq}_1 = H^q(\ke\otimes \Omega^{-p}(-p))\otimes \ko(-p) $$
which converges to 
$$ E^i = \left\{ \begin{array}{l} 
\ke \ \text{ for \ } i=0, \\
0 \ \text{ otherwise.}
\end{array} \right.
$$
\end{theorem}

The next step is to pass from the Beilinson spectral sequence of type II to a 
monad for instanton bundles given by the invariants from the lemma \ref{Hinst}.

On the $E_1$-level we have only 3 non-zero terms:
$$\begin{array}{lcl}
E^{-s,s-1}_1 & = & H^{s-1}(\ke(-s))\otimes \Omega^s(s) \\
E^{-1,1}_1 & = & H^1(\ke(-1))\otimes \Omega^1(1) \\
E^{0,1}_1 & = & H^1(\ke)\otimes \ko 
\end{array} $$ 
and when written out together with differentials looks like
$$ \xymatrix@C5ex@R3ex{
0 \ar[r] & E^{-s,s-1}_1 \ar[r] \ar@{}[dr]|{\ddots} & 0 \ar[r] & {\dots} \\
& {\dots} \ar[r] & 0 \ar[r] & 0 \ar[r] & 0 \ar[r] & {\dots} \\
& {\dots} \ar[r] & 0 \ar[r] & E^{-1,1}_1 \ar[r]^{\nu} 
                 & E^{0,1}_1 \ar[r] & {0} 
}$$
Computing the cohomologies of these differentials, we see 
$E^{-s,s-1}_2=E^{-s,s-1}_1,\ E^{-1,1}_2 = \ker \,\nu,\  
E^{0,1}_2 = \coker \,\nu$. If $s>2$, the differentials in $E_2$ are zero:
$$ \xymatrix@C5ex@R3ex{
0 \ar[rrd] & {\dots} \\
0 & 0\ar[rrd] & E^{-s,s-1}_k \ar[rrd] & 0 & {\dots} \\
0 & 0 & 0 \ar[rrd] & 0 \ar[rrd] & 0 & {\dots}\\
& {\dots} & 0 & 0 & {\ker\,\nu} \ar[rrd] & 
        {\coker\,\nu} \ar[rrd] & 0 & {\dots} \\
&&&& {\dots} & 0 & 0 & 0  
}$$
and analogous for all $E_k$, $1<k<s$. The last level with nonzero differentials
is $E_s$:
$$ \xymatrix@C5ex@R3ex{
0 \ar[rrrdd] & {\dots}\\
&0 & 0 & 0 & 0 \\
&0\ar[rrrdd] & 0\ar[rrrdd] & E^{-s,s-1}_s \ar[rrrdd]^{\mu} & 
     0 \ar[rrrdd] & 0 & {\dots} \\
&0 & 0 & 0 & 0 & & 0 & 0 & {\dots}\\
&& {\dots} & 0 & 0 & 0 & {\ker\,\nu} \ar[rrrdd] & 
        {\coker\,\nu} \ar[rrrdd] & 0 & {\dots} \\
&&&&&& {\dots} & 0 & 0 & 0 \\
&&&&&&&&& 0 & 0 
}$$
where $\mu: E^{-s,s-1}_1 = E^{-s,s-1}_1 \to \ker\,\nu$ is some 
homomorphism. 

We get the following information about $E_\infty$: 
$$ 
E^{-s,s-1}_\infty \ =\ \ker\, \mu \quad ; \quad
E^{-1,1}_\infty \ =\ \coker\, \mu \quad ; \quad
E^{0,1}_\infty \ =\ \coker\, \nu ,
$$
and all other terms vanish.
In this circumstances the Beilinson theorem \ref{BeilSpectral} implies: \\
$ E^{-s,s-1}_\infty \ =\ \ker\, \mu = 0 $ hence 
$ \mu: E^{-s,s-1}_1 \to ker\, \nu \subset E^{-1,1}_1 $
is injective; \\
$E^{0,1}_\infty \ =\ \coker\, \nu = 0$ hence 
$ \nu: E^{-1,1}_1 \to E^{0,1}_1 $ is surjective, and \\
$\ke \cong E^{-1,1}_\infty\, = \, \coker\, \mu$,
i.e. we have the

\begin{theorem} \label{InstMonad}
Every mathematical instanton bundle with invariants as in the 
definition \ref{Definst} and lemma \ref{Hinst}, is the cohomology of the monad
$$
0 \to H^{s-1}(\ke(-s)) \otimes \Omega^s(s) \stackrel{\mu}{\to}
    H^1( \ke(-1) ) \otimes \Omega^1(1) \stackrel{\nu}{\to}
 H^1(\ke) \otimes \ko \to 0 $$
where $\mu$ is a subbundle and $\nu$ is surjective.
\end{theorem}

An analogous procedure can be carried out to get a monad from the spectral 
sequence of type I.

There is also a stronger version of the Beilinson theorem, leading to the 
same monad in our case.

\begin{theorem} {\bf ( Beilinson ) } \  \cite{AO3} \
Any coherent $\ko$-module $\kf$ on $\P^s$ 
is the cohomology  at 0 of a complex 
$$ 0 \to C^{-s} \to \dots \to C^0 \to C^1 \to \dots \to C^s \to 0 $$
which is exact except at 0 and for which the sheaves are given by
$$ C^p = \bigoplus_i H^i(\kf \otimes \Omega^{i-p}(i-p)) \otimes \ko (p-i)$$
or by
$$ C^{-p} = \bigoplus_i H^i(\kf(-i-p))\otimes \Omega^{p+i}(p+i) $$
We call such complexes {\bf Beilinson monads of types I and II} respectively.
\end{theorem}

Note that the homomorphisms of the monad are not {\it a priori} canonically 
given by $\kf$ in general; any two such monads for a sheaf $\kf$ are known 
only to be homotopically equivalent. The reason is basically that the 
higher differentials of a spectral sequence are induced in a non-unique way.
In the case of mathematical instantons, however, one can express the 
homomorphisms in the monad via some natural maps in cohomologies, see later 
the supplement to the Beilinson theorem \ref{BeilSuppl}.  

We shall be working with the monad of type II most of the time. 
Inserting the invariants for our rank-4 bundles on $\P^4$ with $c_2=2$, we can 
rewrite this monad as
$$ 0 \to k^2\otimes\Omega^4(4) \stackrel{\mu}{\to} k^2\otimes\Omega^1(1) 
                 \stackrel{\nu}{\to} k^2\otimes\ko \to 0. $$

\begin{thought} \label{MtxMN}
The morphisms of sheaves $\mu,\nu$ from the monad
$$
0 \to H^{s-1}(\ke(-s)) \otimes \Omega^s(s) \stackrel{\mu}{\to}
    H^1( \ke(-1) ) \otimes \Omega^1(1) \stackrel{\nu}{\to}
 H^1(\ke) \otimes \ko \to 0 $$
can be identified with matrices 
$M\in Mat_{n\times n}(\bigwedge^{s-1}V)$ and
$N\in Mat_{((s-1)n-r)\times n}(V)$ of suitable sizes using the canonical 
ismomorphism
$\hom_{\P^s}(\Omega^k(k),\Omega^j(j)) = \bigwedge^{k-j}V$ 
for any $0\le j\le k \le n$.
\end{thought}

There are at least two equivalent ways to describe this correspondence between 
sheaf homomorphism and matrices.

For the first construction note that 
$\hom_{\P^s}(\Omega^k(k),\Omega^j(j)) = \hom_{\P^s}(\Omega^k(k+1),\Omega^j(j+1))$.
The sheaves $\Omega^k(k+1),\Omega^j(j+1)$ are generated by their global sections, and 
$$ H^0(\Omega^k(k+1))=\bigwedge^{k+1}V^\vee , \quad 
H^0(\Omega^j(j+1))=\bigwedge^{j+1}V^\vee.$$  
Take an exterior form $\omega \in \bigwedge^{k-j}V$ and define a $k$-linear map by contraction 
with $\omega$:
$$\xymatrix{
H^0(\Omega^k(k+1)) \ar@{}[r]|{=}  & 
    {\bigwedge^{k+1}V^\vee} \ar[r]^{\lrcorner \omega} \ar[d]_{\cong} &
    {\bigwedge^{j+1}V^\vee} \ar@{}[r]|{=} \ar[d]^{\cong} &  
    H^0(\Omega^j(j+1)) \\
& {\bigwedge^{s-k}V} \ar[r]^{\wedge \omega} & {\bigwedge^{s-j}V} 
}$$ 
One can prove that such linear map defines a sheaf homomorphism.

The second definition uses the Euler sequences 
$$ 0 \to \Omega^k(k) \to \bigwedge^{k}V^\vee\otimes\ko \to \Omega^k(k+1) \to 0, $$
$$ 0 \to \Omega^j(j) \to \bigwedge^{j}V^\vee\otimes\ko \to \Omega^j(j+1) \to 0, $$
One can define a morphism 
$$\xymatrix{
{\phi} \ar@{}[r]|{:\qquad} & 
    {\bigwedge^{k}V^\vee \otimes \ko} \ar[r]^{\lrcorner \omega} \ar[d]_{\cong} &
    {\bigwedge^{j}V^\vee \otimes \ko} \ar[d]^{\cong} \\ 
& {\bigwedge^{s-k}V\otimes\ko} \ar[r]^{\wedge \omega} & {\bigwedge^{s-j}V\otimes\ko} 
}$$ 
and show that $\phi(\Omega^k(k)) \subset \Omega^j(j)$.


\begin{thought} \label{BeilSuppl}
{\bf Supplement / Special case of the Beilinson theorem.} 
\ (see \cite{script}) 
Let $\mu$,$\nu$ be 
the morphsism from the type-II Beilinson monad for the instanton bundle $\ke$
as in theorem \ref{InstMonad} and 
$M\in Mat_{n\times n}(\bigwedge^{s-1}V)$, 
$N\in Mat_{((s-1)n-r)\times n}(V)$ $M,N$ be matrices of $\mu,\nu$ constructed 
as in \ref{MtxMN} , where $n=h^{s-1}(\ke(-s))=h^1(\ke(-1))$ and  $(s-1)n-r=h^1(\ke)$. 
Then via the identifications \small
$$
Mat_{n\times n}(\bigwedge^{s-1}V)\, = \,
   \hom_k (\,H^{s-1}(\ke(-s))\,,\, H^1(\ke(-1))\otimes \bigwedge^{s-1}V\,) 
\, = \,
  \hom_k (\,H^{s-1}(\ke(-s))\otimes \bigwedge^{s-1}V^\vee\, ,\, H^1(\ke(-1))\,)
$$ \normalsize
respectively
$$
Mat_{((s-1)n-r)\times n}(V)\, = \,
   \hom_k (\,H^1(\ke(-1))\,,\, H^1(\ke)\otimes V\,) \, = \,
   \hom_k (\,H^1(\ke(-1))\otimes V^\vee\, ,\, H^1(\ke)\,)
$$ 
$M$ and $N$ can be canonically written as:
$$
  N \ : \ H^1\ke(-1) \otimes V^\vee \ =\  H^1\ke(-1) \otimes H^0\ko(1) 
      \stackrel{\text{cup}}{\longrightarrow} H^1\ke
$$
$$
  M \ : \ H^{s-1}\ke(-s) \otimes \bigwedge^{s-1}V^\vee \ = \ 
         H^{s-1}\ke(-s) \otimes H^0(\ke \otimes \Omega^{s-2}(s-1)) 
         \stackrel{\text{cup}}{\longrightarrow}
         H^{s-1}(\ke \otimes \Omega^{s-2}(-1)) \ \cong \
         H^1(\ke(-1)) 
$$
where the last isomorphism the inverse to the composition
$$
 H^1\ke(-1) \stackrel{\cong}{\to} H^2(\ke\otimes\Omega^1(-1))
\stackrel{\cong}{\to} H^3(\ke\otimes\Omega^2(-1))  \stackrel{\cong}{\to} \dots
 \stackrel{\cong}{\to}  H^{s-1}(\ke\otimes\Omega^{s-2}(-1))
$$
of connecting homomorphisms of the  long exact homology sequences associated 
to exact sequences of sheaves 
$$
0 \to \ke\otimes\Omega^p(-1) \to \bigwedge^pV^\vee \otimes \ke(-p-1) \to 
  \ke \otimes \Omega^{p-1}(-1) \to 0.  
$$
\end{thought}

The connecting homomorphisms are isomorphisms because $h^i\ke(p-1)=0$ for
all $i$.


\begin{construction}
There is a direct procedure to construct a monad (almost) of type I for an 
instanton bundle $\ke$ from a monad 
of type II for the same bundle $\ke$. In the additional assumption $H^0\ke=0$ 
we can also do the inverse 
procedure, i.e. construct a monad of type II from that of type I.
\end{construction}

Given a monad of type II 
$$ 0 \to A \otimes \Omega^s(s) \stackrel{\mu}{\longrightarrow} 
   B \otimes \Omega^1(1) \stackrel{\nu}{\longrightarrow}
   C \otimes \ko \to 0,
$$
with some $k$-vector spaces $A,B,C$,
we write the Euler sequence tensorized with $B$
$$
0 \to B \otimes \Omega^1(1) \longrightarrow B \otimes V^\vee \otimes \ko 
  \longrightarrow B \otimes \ko(1) \to 0
$$
and embed the first map into the commutative diagram
$$\xymatrix{
& & B \otimes \Omega^1(1) \ar[d] \ar[r]^{\nu} & 
    C\otimes \ko \ar@{=}[d] \ar[r] & 0\\
0 \ar[r] & H \otimes \ko \ar[r] & B\otimes V^\vee \otimes \ko \ar[r]^{N} &
    C\otimes \ko \ar[r] & 0\\ 
}$$
where $N$ is a matrix for $\nu$ as in \ref{MtxMN} and we define
$H := \ker (B\otimes V^\vee \stackrel{N}{\longrightarrow} C )$. Now we can 
complete the diagram to the display
$$\xymatrix{
& & & 0 \ar[d] \\
& 0 \ar[r] \ar[d] & A\otimes \Omega^s(s) \ar[r]^{\mu} \ar@{=}[ld] &
  B\otimes \Omega^1{1} \ar[dd] \ar[r]^{\nu} & 
  C\otimes \ko \ar@{=}[dd] \ar[r] & 0 \\
& A\otimes \bigwedge^{s+1}V^\vee \otimes \ko(-1) \ar[d]^{\alpha} \\
0 \ar[r] & H \otimes \ko \ar[rr] \ar[d]^{\beta} && 
  B \otimes V^\vee \otimes \ko \ar[r]^{\quad N} \ar[d] & C\otimes\ko \ar[r] 
  & 0 \\
& B\otimes\ko(1) \ar@{=}[rr] \ar[d] && B\otimes\ko(1) \ar[d] \\
& 0 && 0
}$$

Diagram chase yields that the cohomology sheaf of the left column is 
again $\ke$.

It can be shown that 
\begin{lemma} \cite[lemma 1.10]{script} The left column when twisted by -1 
gives the Beilinson monad of $\ke(-1)$.
\end{lemma}

{\bf Proof.} The proof of the lemma consists essentially in showing that the 
discription of 
maps in the monad of type I as of cohomological operators coincides with
the definition of $\alpha$ and $\beta$ by this display. \qed

Conversely, if we start from a monad of the form
$$
  0 \to A\otimes\bigwedge^{s+1}V^\vee\otimes\ko(-1) \longrightarrow
    H \otimes \ko \longrightarrow B\otimes \ko(1) \to 0
$$
with the cohomology sheaf $\ke$ such that $H^0\ke=0$, than we can factorize 
$\beta$ through $B\otimes \ko$ 
$$\xymatrix{
H\otimes \ko \ar[r]^{b\quad} \ar[d]_{\beta} & 
     B\otimes V^\vee \otimes \ko \ar[d]^{\text{can}} \\
B\otimes\ko(1) \ar@{=}[r] & B\otimes\ko(1)
}$$
and define $C:=\coker\, b$. Further, $b$ is injective since
$$
 \ker b = \ker\, H^0(\beta) = H^0(\kn) = H^0(\ke) = 0  
$$
where $\kn$ is the kernel sheaf of the monad, i.e. $\kn = \kk er(\beta)$
and the last equality comes from the exact sequence
$$
  0 \to A\otimes\bigwedge^{s+1}V^\vee\otimes\ko(-1) \longrightarrow \kn
  \longrightarrow \ke \to 0
$$
as $H^1\ko(-1)=0$. All this said, we can reconstruct the display as above with
the monad of type I as its upper row.


{\bf Open conditions on $M$ and $N$}

\begin{thought}
Let $M\in Mat_{n\times n}(\bigwedge^{s-1}V)$, 
$N\in Mat_{((s-1)n-r)\times n}(V)$ be two matrices defining maps $\mu$ and 
$\nu$ in
$$
0 \to k^n \otimes \Omega^s(s) \stackrel{\mu}{\to}
    k^n \otimes \Omega^1(1) \stackrel{\nu}{\to}
 k^{((s-1)n-r)} \otimes \ko \to 0. $$
In order to define a monad of an instanton bundle, $M$ and $N$, or, 
equivalently, $\mu$ and $\nu$ must subject to the following conditions:\\
(i) $\nu \circ \mu = 0$ or, equivalently, 
$N\wedge M=0\in Mat_{((s-1)n-r)\times n}(\bigwedge^sV) $;
(ii) $\mu$ is a subbundle ; \\
(iii) $\nu$ is surjective.
\end{thought}

In the condition (i), the operation $\wedge$ between matrices means their 
product as matrices over the exterior algebra $\bigwedge^\cdot V$. 
The condition (i) here is a closed condition, and the conditions (ii),(iii) are
open conditions with respect to the Zarisski topology on the affine space 
$ Mat_{n\times n}(\bigwedge^{s-1}V) \, \times \, Mat_{((s-1)n-r)\times n}(V)$.

Now it is time to describe these two open conditions in algebraic terms.

\begin{defin} 
(i) Nontrivial linear combinations of columns respectively rows of a matrix 
are also called its {\bf generalized}, as in \cite[Lecture 9]{harris}, or 
{\bf generated}, as in \cite{script}, columns respectively rows. \\
(ii) We say that a column or a row of a matrix with entries in some vector
spece $U$ has the form 
$\lambda \otimes z$, where $z \in U$, $\lambda \in k^p$, if it is equal to
$(\lambda_1\cdot z,\lambda_2\cdot z,\dots,\lambda_p\cdot z)$.  
\end{defin}

\begin{lemma} \cite[2.1]{script} \label{OpenN}
The following conditions are equivalent: \\
(i) The homomorphism $k^n\otimes \Omega^1(1) \stackrel{\nu}{\to} 
k^{(s-1)n-r}\otimes\ko$ is surjective;\\
(ii) There is no nontrivial linear combinations of the \underline{rows} of 
$N$ of the 
form $\lambda \otimes x$, $x \in V$, $\lambda \in k^n-\{0\}$; \\
\end{lemma} \qed

\begin{remark} \label{remOpenN} 
In the case when $N$ is a $2\times 2$-matrix we easily get one
more condition equivalent to the previous two:\\
(ii') There is no nontrivial linear combinations of the \underline{columns} of 
$N$ of the 
form $\lambda \otimes x$, $x \in V$, $\lambda \in k^n-\{0\}$ 
\end{remark}

\begin{lemma} \cite[2.2]{script} \label{OpenM}
The following conditions are equivalent: \\
(i) The homomorphism $k^n\otimes \Omega^s(s) \stackrel{\mu}{\to} 
k^n\otimes\Omega^1(1)$ is a subbundle;\\
(ii) For any $0\ne x\in V$ the matrix $x\wedge M$ which induces the map 
$$\xymatrix{
 \mu(x) \ar@{}[r]|{:} & \Omega^s(s)(x) \ar@{}[r]|{=}  & 
        {\bigwedge^s(V/x)^\vee} \ar[d]_{\cong} \ar[r]^{\lrcorner M} & 
        (V/x)^\vee \ar@{}[r]|{=} \ar@{^{(}->}[d]  & \Omega^1(1)(x) \\
 & & {\bigwedge^{s+1}V^\vee} \ar[r]^{\quad \lrcorner (x\wedge M)} & V^\vee
}$$
is an injective linear operator; \\  
(iii) There is no generalized \underline{column} of $M$ of the form 
$x\wedge m$ where
$m$ is a column with entries in $\bigwedge^{s-1}V$ and $x \in V$; \\ 
(iv) The matrix $P$ representing the kernel of $M$
$$ 0 \to k^p \stackrel{P}{\longrightarrow}  
    k^n\otimes V \stackrel{M}{\longrightarrow} 
    k^n\otimes \bigwedge^s V
$$
has no generalized column of the form $\lambda\otimes x$, $x\in V$, 
$\lambda \in k$.
\end{lemma} \qed


For the case $\dim V=5$ and $M \in Mat_{2\times 2}(\bigwedge^3 V)$ which is 
relevant to instantons with our invariants we can develop another criterion 
determining wheather $M$ defines a subbundle. This criterion will be used, e.g.,
in the proof of \ref{ClassNprim} .

Suppose we are given a generated column $(\xi,\eta)^T$ of the matrix $M$ and
we want to know if this column is forbidden by the condition (iii) of the 
previous lemma. 

\begin{lemma} \label{WedgeLemma}
Let $V$ be a 5-dimensional $k$-vector space. Two exterior 
forms $\eta, \xi \in \bigwedge^3 V$ have a common 
linear factor if and only it the two equalities hold:
$$
   (\xi^{*2})^*\wedge \eta =0 \ ; \ (\eta^{*2})^*\wedge \xi =0 . \quad (*)
$$
where * denotes a dualization map 
$\bigwedge^{\cdot}V \stackrel{\cong}{\longrightarrow} 
 \bigwedge^{5-\cdot}V^\vee$.
\end{lemma}

{\bf Proof.} If $\xi$ or $\eta$ is a zero form, we have nothing to prove.\\
Choose a basis of $V$ such that $\xi=e_{012}$ or 
$\xi=e_0 \wedge(e_{12}+e_{34})$  
(we abbreviate $e_i \wedge e_j \wedge \dots \wedge e_k$ by $e_{ij\dots k}$). 

Case 1 : $\xi=e_{012}$ (we say that $\xi$ is decomposable). Then up to a scalar 
multiple $\xi^*=e^*_{45}$  and $\xi^{*2}=0$. 
If $\eta$ is decomposable, too, then, on one hand, conditions (*) are fulfilled and, on the other hand, 
$k \cdot \xi=\bigwedge^3 V_{\xi}$, $k \cdot \eta=\bigwedge^3 V_{\eta}$,
where  $V_{\xi}, V_{\eta}$ are 3-dimensional subspaces in $V$. Now a vector
$0 \ne v \in V_{\xi} \cap V_{\eta}$ fulfills 
$\xi = v \wedge \xi '\,,\, \eta = v \wedge \eta '$.

Suppose one of $\xi,\eta$, say, $\xi$ is indecomposable and arrive at the \\
Case 2 : $\xi = e_{012} + e_{034}$. Then up to a scalar multiple 
$\xi^* = e^*_{34} + e^*_{12}$, $\xi^{*2} = e^*_{1234}$ and, finally, 
$(\xi^{*2})^* = e_0$. We see that the linear factor of $\xi$ is precisely
$(\xi^{*2})^*$. It is also a linear factor in $\eta$ iff 
$(\xi^{*2})^*\wedge \eta =0$. \\
If $\eta$ is decomposable, then the second equality is trivially satisfied, 

If both $\xi$ and $\eta$ are indecomposable, they both have linear factors 
$(\xi^{*2})^* \ni 0$ and $(\eta^{*2})^* \ni 0$ correspondently, and our equalities,
equivalent in this situation, mean simply that they are collinear.\\
Hence the lemma. \qed

\footnotesize

The similar lemmata can be proved also for $\bigwedge^2 k^4$ and for $\bigwedge^4 k^6$, i.e.
for spaces of exterior forms relevant to the cases of instantons on $\P^3$ and $\P^5$ respectively.

{\bf Lemma.} {\it Two exterior forms $\eta, \xi \in \bigwedge^4 k^6$ have a 
common 
linear factor if and only it the four equalities hold:
$$
   \xi^{*3}=0 \ ; \ \eta^{*3}=0 \ ; \ (\xi^{*2})^*\wedge \eta =0 \ ; \ (\eta^{*2})^*\wedge \xi =0 . 
$$ }

{\bf Sketch of the proof.} The idea is the same as above. The main points are:\\
For every form $0 \ne \xi \in \bigwedge^4 k^6$ we can choose a basis of $k^6$ such that  
$\xi^*=e^*_{12}$ or $\xi^*=e^*_{12}+e^*_{34}$ or $\xi^*=e^*_{12}+e^*_{34}+e^*_{56}$.\\
The condition $\xi^{*3}=0$ is equivalent to $V_\xi \ne 0$. If so, then either $\xi^{*2}=0$ and 
$\dim\,V_\xi=4$ or $\bigwedge^2 V_\xi = k\cdot (\xi^{*2})^*$. \qed

{\bf Lemma.} {\it Two exterior forms $\eta, \xi \in \bigwedge^2 k^4$ have a 
common linear factor if and only it the three equalities hold:
$$
 \xi^2=0 \ ; \ \eta^2=0 \ ; \  \xi \wedge\eta =0.
$$ }

{\bf Proof} is obvious. 

\normalsize

\begin{prop} \label{ColWedge} Let a matrix  
$$ 
M = \left( \begin{array}{cc} m_{11} & m_{12} \\ m_{21} & m_{22} 
    \end{array} \right) 
$$
define a sheaf homomorphism $2\Omega^4_{\P^4}(4) \to 2\Omega^1_{\P^4}(1)$, 
where $m_{ij} \in \bigwedge^3 V$, $V\cong k^5$. Then $M$ defines a subbundle 
if and only if for no $(s,t)\in \P^1$ the equalities
$$
((s\cdot m_{11}+t\cdot m_{12})^{*2})^* \wedge 
      (s\cdot m_{21}+t\cdot m_{22}) = 0; \quad
((s\cdot m_{21}+t\cdot m_{22})^{*2})^* \wedge 
      (s\cdot m_{11}+t\cdot m_{12}) = 0
$$
are satisfied simultaneously.
\end{prop}

{\bf Proof} follows from lemmata \ref{WedgeLemma} and \ref{OpenM}. \qed

{\bf The Group Action}

From now on we consider only the case with our fixed invariants, i.e.,
rank-4 bundles on $\P^4$ with $c_2=2$.
Our next aim is to construct a space parametrizing these instanton bundles.

\begin{thought} There are two natural group operations on the set of 
all monads
$$ 0\to A\otimes \Omega^4(4) \stackrel{M}{\longrightarrow}
    B\otimes \Omega^1(1) \stackrel{N}{\longrightarrow}
    C\otimes \ko \to 0
$$
where $A,B,C$ are 2-dimensional vector spaces, namely, those of 
$GL(A)\times GL(B)\times GL(C)$ and of $GL(V)$, where, as usual, $\P^4=\P V$.
\end{thought}

Let us first define the operation of the
group $GL(A)\times GL(B)\times GL(C)$. 
If $(g_1,g_2,g_3) \in GL(A)\times GL(B)\times GL(C)$, we get a commutative 
diagram
$$\xymatrix{
0 \ar[r] & A\otimes \Omega^4(4) \ar[r]^{M} \ar[d]_{g_1\otimes \Id}&
    B\otimes \Omega^1(1) \ar[r]^{N} \ar[d]_{g_2\otimes \Id} 
    & C\otimes \ko \ar[r] \ar[d]_{g_3\otimes \Id} & 0 \\
0 \ar[r] & A\otimes \Omega^4(4) \ar[r]^{g_2Mg_1^{-1}} &
    B\otimes \Omega^1(1) \ar[r]^{\, g_3Ng_2^{-1}} & C\otimes \ko \ar[r] & 0
}$$
and the cohomologies of the both rows are clearly isomoprphic.
In particular, if we replace matrix $M$ by $\lambda M$ for 
some $\lambda\in k-\{0\}$, the cohomology sheaf will not change, and the same 
is true for N. Slightly reformulating, we see that the whole $G$-orbit of 
$$(\langle M \rangle,\langle N \rangle) \in \P Mat_{2\times 2}(\bigwedge^3 V)
\times \P Mat_{2\times 2}(V)$$ 
defines the same cohomology sheaf, where $G=SL(2)\times SL(2)\times SL(2)$
and the operation is given by
$$ (g_1,g_2,g_3)(\langle M \rangle, \langle N \rangle) := 
  (\langle g_2^{-1} M g_1 \rangle , \langle g_3^{-1} N g_2 \rangle   ) .  $$
 
Later on we shall see, that there exist a 1:1 correspondence between 
$G$-orbits and the cohomology bundles.

We also have a $GL(V)$-operation on our monads that corresponds to changing
coordinates in our $\P^4$ and pull-backs of the cohomology sheaves under
these coordinate transformations.

On the set $ Mat_{2\times 2}(V) $ of matrices $N$ we have therefore the 
restricted group action of \linebreak 
$GL(B)\times GL(C) \times GL(V)$ which should 
be thought of as of 
$$\left( \begin{array}{c} 
\text{column} \\ \text{transformations}
\end{array} \right) \times 
\left( \begin{array}{c} 
\text{row} \\ \text{transformations}
\end{array} \right) \times 
\left( \begin{array}{c} 
\text{change of } \\ \text{basis in $V$}
\end{array} \right) .
$$

\begin{thought} \label{MotivClass}
The questions listed below may be solved on representatives of 
$GL(B)\times GL(C) \times GL(V)$-orbits, i.e. if we know an answer/proof for
one particular matrix $N_0$, we also have it for all matrices $N$ equivalent 
to $N_0$ modulo this group action: \\
(i) to describe the set of matrices $M\in Mat_{2\times 2}(\bigwedge^3V)$ such
that $N\wedge M=0$ (cf. \ref{SyzN} and, for different monads, proposition 
\ref{ClassNprim}); \\
(ii) to study geometry of cohomology sheaves of monads 
$$ 0\to A\otimes \Omega^4(4) \stackrel{M}{\longrightarrow}
    B\otimes \Omega^1(1) \stackrel{N}{\longrightarrow}
    C\otimes \ko \to 0
$$
for all $M$'s with $N\wedge M=0$; \\
(iii) to study properties of the parameter space $X_0$ (cf. \ref{X0smooth}) 
and moduli
space (defined in the section \ref{formal}) at points corresponding to the 
pairs of matrices $(M,N)$ with $N\wedge M$=0,\\
etc.
\end{thought}

After the next lemma we shall see that for our purposes we need to consider 
only two \linebreak $GL(B)\times GL(C) \times GL(V)$-orbits and hence have to 
do a lot of things for two fixed matrices $N$ only.

\begin{lemma} \label{span3} 
Given a matrix $A \in Mat_{2\times 2}(U)$ with $U$ being a $k$-vector space. 
Suppose that 
$$ \Span\, A\ :=\ \Span\,\{ a_{11}, a_{12}, a_{21}, a_{22} \} \ \cong \ k^3 $$
and that $A$ does not generate column of the form $\lambda\otimes v$ for 
$\lambda \in k^2,\, v\in U$. 
Then there exist an element $g\in GL(2)$ and three linearly independent 
vectors $e_1,e_2,e_3$ such that
$$ A\cdot g = \left( \begin{array}{cc}
e_1 & e_2 \\ e_3 & e_1  
\end{array} \right) . $$
\end{lemma}

{\bf Proof.} Note first, that left multiplication of a matrix $A$ with a 
group element $g\in GL(2)$ corresponds to elementary column operations in $A$.

Without loss of generality, 
$A = \left( \begin{array}{cc} a & b \\ 
             c & \lambda a + \mu b + \nu c \end{array} \right)$, 
where $a,b,c \in U;\ \lambda,\mu,\nu \in k$, where $a,b,c$ is a basis of 
$\Span\,N$ (otherwise permute rows and columns)
Then we have the chain of transformations shown symbolically on the diagram
and explained in details beneath:
$$\xymatrix {
{\left( \begin{array}{cc} a & b \\ c
& \lambda a + \mu b + \nu c \end{array} \right)}
\ar@{~>}[r]_{\txt{\tiny \qquad  ${}^{\text{elem.}}_{\text{oper.}}$}}
& 
{\left( \begin{array}{cc} a & b-\nu a \\ c
& \lambda a + \mu b  \end{array} \right)}
\ar@{~>}[r]_{\txt{\tiny new $b, \lambda$}}
&
{\left( \begin{array}{cc} a & b \\ c
& \lambda a + \mu b \end{array} \right)}
\ar@{~>}[r]_{\txt{\tiny elem.op. \\ \tiny $\lambda \ne 0$ for $M \in X_0$ }} 
&
{\left( \begin{array}{cc} \lambda a & b \\ \lambda c
& \lambda a + \mu b \end{array} \right)}
\\ {\hspace{3cm}} 
\ar@{~>}[r]_{\txt{\tiny \qquad new $a,c$}}
& 
{\left( \begin{array}{cc} a & b \\ c
& a + \mu b \end{array} \right)}
\ar@{~>}[r]_{\txt{\tiny ${}^{\text{elem.}}_{\text{oper.}}$\qquad \qquad \qquad }}
&
{\left( \begin{array}{cc} a + \mu b & b \\ c + \mu a + \mu^2 b
& a + \mu b \end{array} \right)}
\ar@{~>}[r]_>>{\txt{\tiny new $a,c$ \qquad \qquad \qquad}}
&
{\left( \begin{array}{cc} a  & b \\ c 
& a \end{array} \right)} .
}
$$
To perform these trasformations, you have to do the following steps.\\
Start with 
$\left( \begin{array}{cc} a & b \\ c
  & \lambda a + \mu b + \nu c \end{array} \right)$. \\ 
On the first step substract from the second column the first multiplied by 
$\nu$, get 
$\left( \begin{array}{cc} a & b-\nu a \\ c
  & \lambda a + \mu b  \end{array} \right)$.\\
Since $a,b,c$ are linearly independent, so are $a,b-\nu a,c$, and we may take
$b-\nu$ as a new vector $b$; in order to have the same vector in the low right
corner, we need to choose a new value for $\lambda$ as $\lambda+\mu \nu$; the
resulting matrix is then
$\left( \begin{array}{cc} a & b 
  \\ c & \lambda a + \mu b  \end{array} \right)$.\\
If $\lambda=0$, we have a column of forbidden form $(b,\mu b)$, hence 
$\lambda\ne 0$ and we may multiply the first column by $\lambda$ and get
$\left( \begin{array}{cc} \lambda a & \lambda b \\ c
  & \lambda a + \mu b \end{array} \right)$.\\
Scaling $a:=\lambda a,\, b:=\lambda b$, we arrive at a matrix
$\left( \begin{array}{cc} a & b \\ c & a + \mu b \end{array} \right)$;\\
In the next step, add $\mu$ times the second column to the first column and
obtain
$\left( \begin{array}{cc} a + \mu b & b \\ c + \mu a + \mu^2 b
  & a + \mu b \end{array} \right)$.\\[0.5ex]
Since $a + \mu b , b , c + \mu a + \mu^2 b$ are linearly independent because
$a,b,c$ are, we may call the three former vectors to be a new triple $a,b,c$ 
and finally get 
$\left( \begin{array}{cc} a  & b \\ c & a \end{array} \right)$
that completes the proof.
\qed

\begin{corollary} \label{ClassN}
A matrix $N$ defines a surjection 
$2\Omega^1(1) \to 2\ko$ if and only if  $N$ is equivalen modulo the 
$GL(B)\times GL(C) \times GL(V)$-action to one of the two matrices
$\left(\begin{array}{cc} e_1 & e_2 \\ e_3 & e_4
\end{array} \right)$ or 
$\left(\begin{array}{cc} e_1 & e_2 \\ e_3 & e_1 \end{array} \right)$.
\end{corollary}

{\bf Proof.}
From lemma \ref{OpenN} it is clear that if $N$ defines a surjection 
$2\Omega^1(1) \to 2\ko$ than $\dim\,\Span\, N\,\ge 3$. 
If $\dim\,\Span\, N\, = 4$, the conclusion is clear, otherwise 
apply remark \ref{remOpenN} and 
the lemma above. \qed

{\bf Syzygies of $N$'s}

\begin{thought} \label{SyzN} 
For both possible types of $N$, we can parametrize the set of all matrices
$M\in Mat_{2\times 2}(\bigwedge^3 V)$ such that $N\wedge M = 0$
by means of some number (in fact, 20) of scalar parameters $p_i,q_i\in k$
For both types of $N$ there exist matrices $M$ defining
a subbundle $2\Omega^4(4) \to 2\Omega^1(1)$
\end{thought}

Let $\Gamma \in Mat_{2\times (r+1)}(\bigwedge^3V) $ be a total matrix of 
syzygies of $N$, i.e. $N \wedge \Gamma $=0 and columns of $\Gamma$ 
generate all the syzygies of degree 3 of the matrix $N$. 
Then all possible matrices $M$ with the condition 
$N \wedge M = 0$ are of the form 
$$ 
M = \Gamma \cdot \left[ \left( \begin{array}{ccc} 
p_0, & \dots, & p_r \\
q_0, & \dots, & q_r 
\end{array} \right) ^T \right] . 
$$
 
If $N=\left(\begin{array}{cc} e_1 & e_2 \\ e_3 & e_4
\end{array} \right)$, then
$$\Gamma = \left(\begin{array}{cccccccccc} 
0 & e_{024} &  e_{023}+e_{014} & e_{013} & 0 & e_{124} & e_{123} & 0 & 
    e_{234} & e_{134} \\ 
e_{024} & e_{023}+e_{014} & e_{013} & 0 & e_{124} & e_{123} & 0 & e_{234} & 
    e_{134} & 0
\end{array} \right), $$
and an example of $M$ defining a subbundle is
$$\left(\begin{array}{cc} 
e_{023}+e_{014} & e_{134}+e_{024}+e_{013} \\
e_{124}+e_{013} & e_{023}+e_{014}
\end{array} \right); $$
and if $N=\left(\begin{array}{cc} e_1 & e_2 \\ e_3 & e_1
\end{array} \right)$, we get
$$\Gamma = \left(\begin{array}{cccccccccc} 
0 & e_{012} & e_{023} & e_{013} & 0 & e_{123} & 0 & e_{124} & 
    e_{234} & e_{134} \\
e_{012} & -e_{023} & e_{013} & 0 & e_{123} & 0 & e_{124} & 
    -e_{234} & e_{134} & 0 
\end{array} \right), $$
and, for example, the following matrix $M$ defines a subbundle
$$\left(\begin{array}{cc} 
e_{234}+e_{012} & e_{134}+e_{124}+e_{013} \\
e_{134}-e_{124}-e_{023} & -e_{234}+e_{123}
\end{array} \right). $$

\begin{remark}
We see that the case of instanton bundles given by a monad 
$$
0 \to k^2 \otimes \Omega^s(s) \stackrel{M}{\to}
    k^2 \otimes \Omega^1(1) \stackrel{N}{\to}
   k^2 \otimes \ko \to 0 $$
on $\P^4$ differs from that in $\P^3$  \cite{NarTrm} at the point
that the  $N$ with $\rank\,\Span\, N\, = 3$ are allowed in our
situation.
\end{remark}


{\bf The Parameter Space and its Irreducibility}

\begin{thought}
Denote:
$$Y_0 := \{ \langle N \rangle \in \P Mat_{2\times 2}(V) : 
         2\Omega^1(1) \stackrel{N}{\longrightarrow} 2\ko \ 
             \text{is surjective} \} . $$ 
$$ X := \{ (\langle M \rangle, \langle N \rangle) : N \wedge M = 0 \}  
\subset \P Hom(k^2,k^2 \otimes \bigwedge^3 V) \times 
\P Hom(k^2,k^2 \otimes V), $$
$$ X_0 := \{ (\langle M \rangle, \langle N \rangle) \in X : M \, 
            \text{is subbundle,\,}
            N \, \text{is surjective}\, \} . $$
\end{thought}

We refer to $X_0$ as to the {\bf parameter space}.

\begin{theorem} \label{Xirred}
$X_0$ is irreducible.
\end{theorem}

{\bf Proof}
We have an obvious projection map $X_0 \to Y_0$ and $X_0$ is open in 
$X\times_{\P Mat_{2\times 2}(V)} Y_0 =: X_{Y_0}$. 

The projection $X_{Y_0} \to Y_0$ is clearly a projective morphism with 
equidimensional fibers that are \linebreak 
19-dimensional linear projective subspaces 
of $\P Mat_{2\times 2}(\bigwedge^3V)$, hence $X_{Y_0}$ is irreducible by the 
proposition below, and so is $X_0$. \qed

\begin{prop} \label{IrredCrit}
Given a projective surjective morphism 
$p: S \twoheadrightarrow T$ where $T$ is 
irreducible and geometric fibers of $p$ are equidimensional and irreducible.
Then $S$ is irreducible.
\end{prop} \qed

\newpage

\section{The Dual Bundle}
\label{dual}

The aim of this section is to describe the duals to instanton bundles in 
our case of rank 4, quantum number 2 on $\P^4$ . The main result is:

\begin{theorem} \label{DualTh}
 Given a monad
$$ 
     2\Omega_{\P^4}^4(4) \stackrel{M}{\longrightarrow} 2\Omega_{\P^4}^1(1)
      \stackrel{N}{\longrightarrow} 2\ko_{\P^4} 
$$
with the cohomology bundle $\ke$ and a matrix $P$ defined as the kernel in
$$ 0 \to k^p \stackrel{P}{\longrightarrow}  
    k^2\otimes V \stackrel{\wedge M}{\longrightarrow} 
    k^2\otimes \bigwedge^4 V.
$$
where $\tilde M$ is the contraction with the matrix $M$.
Then the dual bundle $\ke^\vee$ is either 
the cohomology bundle of the monad
$$ 
2\Omega^4(4) \stackrel{M^\vee}{\longrightarrow} 2\Omega^1(1)
      \stackrel{P^\vee}{\longrightarrow} 2\ko 
$$
if $\rank\,\tilde M\,=8$, and is therefore an instanton bundle, too,
or, if $\rank\,\tilde M\,=7$, it is an extension
$$ 
0 \to \ko \to \ke^\vee \to \ke^\vee_1 \to 0
$$
where $\ke_1^\vee$ is a reflexive sheaf appearing as the cohomology of the 
monad
$$ 
2\Omega^4(4) \stackrel{M^\vee}{\longrightarrow} 2\Omega^1(1)
      \stackrel{P^\vee}{\longrightarrow} 3\ko.
$$
The values of $\rank\,\tilde M$ other than 8 or 7 do not appear. 
\end{theorem}

In the proof of the theorem, we have to distinguish different cases depending 
on $h^0\ke^\vee$.
To be able to control this invariant, we use the following

\begin{prop} \cite[1.15]{script}. The complex of vector spaces
$$
0 \to H^3(\ke \otimes \Omega^4) \to 
  H^3(\ke(-4))\otimes\bigwedge^4V^\vee \stackrel{\tilde M}{\longrightarrow}
  H^1(\ke(-1))\otimes V^\vee \stackrel{N}{\longrightarrow} 
  H^1\ke \to 0,
$$
where $\tilde M$ is contraction with $M$, is exact except at 
$H^1(\ke(-1))\otimes V^\vee$ where its cohomology is 
$H^1(\ke\otimes\Omega^4)$.
\end{prop}
\qed

\begin{corollary} \label{rkMtilde}
$h^0\ke^\vee = 8 - \rank\,\tilde M$
\end{corollary}

{\bf Proof}. Apply Serre duality and proposition above to get:
$$
h^0\ke^\vee \, = \, h^1(\ke\otimes\Omega^4) \, = \, 
   \dim\, h^1(\ke(-1))\otimes V^\vee \, - \, \rank\, N\, 
   - \, \rank\,\tilde M = 10 - 2 -  \rank\,\tilde M.  $$
\qed


\begin{lemma} The two equivalent inequalities $h^0\ke^\vee \le 1$ and
$\rank\,\tilde M \ge 7$ hold.
\end{lemma}

{\bf Proof.\ }
Firstly, we have to prove that cases $h^0\ke^\vee > 1$, 
i.e. $\rank\,\tilde M \le 6$ are impossible.

Recall that in \ref{OpenM} we have defined $P$ as the kernel 
of $\tilde M$,i.e. through the exact sequence 
$$ 0 \to k^p \stackrel{P}{\longrightarrow}  
    k^2\otimes V \stackrel{\tilde M}{\longrightarrow} 
    k^2\otimes \bigwedge^s V $$
If $\rank\, M\, = 6$, we get $p=4$ and $P$ is a $4\times 2$ matrix.
By the theorem quoted below
such a matrix can be assumed to be  
$P = \left( \begin{array}{cccc} e_0 & e_1 & e_2 & e_3 
                       \\ e_1 & e_2 & e_3 & e_4 \end{array} \right)$, 
because by lemma \ref{OpenM} $P$ does not generate a column of the form
$(\lambda v,\mu v)$.

A computation shows that 
$ 2\bigwedge^3V  \stackrel{P^\vee}{\longrightarrow} 4\bigwedge^4V$
is a surjection of vector spaces and 
therefore we have an exact sequence
$$
        k^2 \stackrel{M^\vee}{\longrightarrow} 2\bigwedge^3V 
       \stackrel{P^\vee}{\longrightarrow} 4\bigwedge^4V \longrightarrow 0
$$
which is impossible because 
$\dim_k(2\bigwedge^3V) = \dim_k ( 4\bigwedge^4V )  = 20$ and $M\ne 0$. 

If $\rank\,M < 6$, then no matrix $P$ of the size $(\ge 5)\times 2$ without 
generalized columns $\lambda\otimes x$ exists (cf. \cite[Lecture 9]{harris} ) 
hence, the cases with $\rank\,M < 7$ do not appear.
\qed

\begin{theorem} \label{HarrisThm} \cite[Lecture 9]{harris}
Consider a $k$-vector space $W$ and a matrix $A\in Mat_{2\times r}(W)$ that 
has no generalized column of the form $(\lambda x , \mu x)$ and such that 
$\Span\,A\,=W$. Then modulo the group operation 
$GL(2)\times GL(r)$ on $Mat_{2\times r}(W)$ the matrix $A$ is 
equivalent to a matrix of the form
$$
\left( \begin{array}{lllcl|llcl|ccc|lcl}
e_0 & e_1 & e_2 & \dots & e_{a_1-1} & e_{a_{1}+1} & e_{a_{1}+2} & \dots  
     & e_{a_2-1} & & \dots & & e_{a_l+1} & \dots & e_{n-1} \\
e_1 & e_2 & e_3 & \dots & e_{a_1} & e_{a_{1}+2} & e_{a_{1}+3} & \dots  
     & e_{a_2} & & \dots & & e_{a_l+2} & \dots & e_n 
\end{array} \right)
$$
where $(e_0,\dots,e_n)$ is a basis of $W$, $l=n-k$, and $a_1,\dots,a_l$ is a 
sequence of integers.
\end{theorem} \qed


Having done this, we proceed with the easier case of $h^0\ke^\vee =0$.

Construct the type-I monad for $\ke$, which is the left column of the 
display:
$$\xymatrix{
& 0 \ar[d] & 0 \ar[d] \\
0 \ar[r] & 2\ko(-1) \ar[d]_{A} \ar[r]^{M} & 2\Omega^1(1) \ar[d] \ar[r]^{N} 
         & 2\ko \ar@{=}[d] \ar[r] & 0 \\
0 \ar[r] & 8\ko \ar[d]_{B} \ar[r]  & 2V^\vee\otimes\ko  \ar[d] \ar[r] 
         & 2\ko \ar[r]  & 0 \\
& 2\ko(1) \ar[d] \ar@{=}[r] & 2\ko(1) \ar[d] \\
& 0 & 0 
}$$
Dualizing the left column, we get 
$$ 0 \longrightarrow 2\ko(-1) \stackrel{B^\vee}{\longrightarrow} 8\ko
      \stackrel{A^\vee}{\longrightarrow} 2\ko(1) \longrightarrow 0 $$
Since $H^0\ke^\vee = 0$, we use the similar diagram to pass back to the 
monad of type-II for $\ke^\vee$, getting
$$ 0 \longrightarrow 2\Omega^4(4) \stackrel{M_1}{\longrightarrow} 2\Omega^1(1)
      \stackrel{N_1}{\longrightarrow} 2\ko \longrightarrow 0 . $$

{\bf Claim. \ } $M_1=M^T$. \\
By the construction of the Beilinson monad, $M$ is given as the cup-product
$$ \bigwedge^3V^\vee \ \otimes H^3\ke(-4) \ = 
  \ H^0\Omega^2(3)\otimes H^3\ke(-4) \longrightarrow 
   H^3(\ke(-1)\otimes\Omega^2) \ \cong \
  H^2(\ke(-1)\otimes\Omega^1) \ \cong \ H^1(\ke(-1)) $$
Now by the naturality of the Serre duality
$$
\xymatrix{
H^3(\ke^\vee(-4))^\vee \ar[d]_{\text{Serre}}^{\cong} & 
   H^1\ke^\vee(-1))\otimes\bigwedge^3V^\vee \ar[d]^{\text{Serre}}_{\cong} 
                                                  \ar[l]_{M^\vee_1 \ \ } \\
   H^1\ke(-1) & H^3(\ke(-4))\otimes\bigwedge^3V^\vee \ar[l]_{M \qquad}
}$$
the claim follows.

From the sequence 
$$
  0 \longrightarrow k^2 \stackrel{P}{\longrightarrow} k^2 \otimes V  
    \ \stackrel{\tilde M = (-\wedge M)}{\longrightarrow} \ 
                                        k^2 \otimes \bigwedge^4V 
    \stackrel{N}{\longrightarrow} k^2\otimes \bigwedge^5V \longrightarrow 0  
$$
and the nondegeneracy of $P$ we see that
$k^2 \otimes \Omega^1(1) \stackrel{P^\vee}{\longrightarrow} k^N \otimes \ko$ \
is surjective and $N_1=P^\vee$, $H^0\ke^\vee = 0$.

This proves the theorem in the case $\rank\,\tilde M = 8$.


If $h^0\ke^\vee=1$ and $\rank\,\tilde M = 7$, then $n=3$ and $P$ is a 
$2\times 3$ matrix and we 
have a cokernel of $\tilde M$ bigger than $N$. The following commutative 
diagram with the exact row appears:
$$
\xymatrix{
{k^2\otimes V} \ar[r]^{\tilde M \ } 
   & {k^2\otimes \bigwedge^4V} \ar[r]^{N_1} \ar@{->>}[rd]_{N}
   & {k^3\otimes \bigwedge^5V} \ar[r]
   & 0 \\
&& {k^2\otimes \bigwedge^5V} \ar[u]
}$$
We can write $N_1$ in the form $N_1 = \left( \frac{N}{x\ y} \right)$ for some 
$x,y \in V$.  

The next important step of the proof is the 

\begin{lemma} In the notation introduced, 
$N_1 \,:\, 2\Omega^1(1) \to 3\ko$ is never surjective.
\end{lemma}

{\bf Proof} Suppose the contrary and get the commutative diagrams
$$ \xymatrix@C6ex{ 
&&& 0 \ar[d] \\
& 0 \ar[d]  && \ko \ar[d] &&&&& 0 \ar[d] & 0 \ar[d] \\
0 \ar[r] & \kk_1 \ar[d] \ar[r] & 2\Omega^1(1) \ar@{=}[d] \ar[r]^{N_1} 
         & 3\ko \ar[d] \ar[r] & 0 && 
    0 \ar[r] & 2\Omega^4(4) \ar@{=}[d] \ar[r]^{M} & \kk_1 \ar[d] \ar[r] 
         & \ke_1 \ar[d] \ar[r] & 0 \\
0 \ar[r] & \kk \ar[r] \ar[d] & 2\Omega^1(1) \ar[r]^{N} & 2\ko \ar[r] \ar[d] & 0 && 
    0 \ar[r] & 2\Omega^4(4) \ar[r]^{M} & \kk \ar[r] \ar[d] & \ke \ar[r] \ar[d] & 0 \\
& \ko \ar[d]  && 0 &&&&& \ko \ar[d] \ar@{=}[r]  & \ko \ar[d]   \\
& 0  &&&&&&& 0 & 0  
}$$ 
But then we can derive the relation of Chern polynomials 
$c(\ke_1) = c(\ke) = 1 + 2h^2 + 3h^4$ 
which is impossible because $c_4(\ke_1)=0$ as $\rank\,\ke_1\, =\, 3$, that 
proves the lemma.

Hence $N_1$ is degenerate, i.e. without loss of generality
$N_1 = \left( \frac{N}{x\ 0} \right)$ and the two diagrams look like

$$ \xymatrix@C5ex{ 
&&& 0 \ar[d] \\
& 0 \ar[d]  & {\ki} \ar@{^({-}>}[r] & \ko \ar[d] \ar[r] & \kc \ar[r] \ar@{=}[d] & 0 
              &&&& 0 \ar[d] & 0 \ar[d] \\
0 \ar[r] & \kk_1 \ar[d] \ar[r] & 2\Omega^1(1) \ar@{=}[d] \ar[r]^{N_1} 
         & 3\ko \ar[d] \ar[r] & \kc \ar[r] & 0 && 
    0 \ar[r] & 2\Omega^4(4) \ar@{=}[d] \ar[r]^{M} & \kk_1 \ar[d] \ar[r] 
         & \ke_1 \ar[d] \ar[r] & 0 \\
0 \ar[r] & \kk \ar[r] \ar[d] & 2\Omega^1(1) \ar[r]^{N} & 2\ko \ar[r] \ar[d] & 0 &&& 
    0 \ar[r] & 2\Omega^4(4) \ar[r]^{M} & \kk \ar[r] \ar[d] & \ke \ar[r] \ar[d] & 0 \\
& \ki \ar[d]  && 0 &&&&&& \ki \ar[d] \ar@{=}[r]  & \ki \ar[d]   \\
& 0  &&&&&&&& 0 & 0  
}$$
where $\ki$ is an ideal sheaf and $\kc = \ko / \ki$ is a sheaf of length 3 
as we shall see later.
Dualizing the right column of the right matrix, we get
$$
0 \longrightarrow \ki^\vee \longrightarrow \ke^\vee \longrightarrow \ke_1^\vee
  \longrightarrow \ke xt^1(\ki,\ko) \longrightarrow 0
$$
or, since $\ki^\vee \cong \ko \ , \ \ke xt^1(\ki,\ko) \cong \ko $ \ ,
$$ 
0 \longrightarrow \ko \longrightarrow \ke^\vee \longrightarrow \ke_1^\vee
                                               \longrightarrow 0
$$ 

As $H^0(\ke^\vee)=1$, we obtain $H^0(\ke_1^\vee)=0$ and hence we may dualize 
the monad for $\ke_1$ by intermediate passing to the monad of type I, 
dualizing that one and passing 
back to a monad of type II for $\ke_1^\vee$. So, we have a diagram
$$\xymatrix{ 
& 0 \ar[d] & 0 \ar[d] \\
0 \ar[r] & 2\ko(-1) \ar[r]^{M} \ar[d] & 2\Omega^1(1) \ar[r]^{N_1} \ar[d] 
                 & 3\ko \ar[r] \ar@{=}[d] & \kc \ar[r] & 0 \\
0 \ar[r] & 7\ko \ar[r] \ar[d] 
           & 2V^\vee\otimes\ko \ar[r]^{\ \lrcorner N_1} \ar[d] 
           & 3\ko \ar[r] & 0 \\
& 2\ko(1) \ar@{=}[r] \ar[d] & 2\ko(1) \ar[d] \\
& \kc \ar[d] & 0\\
& 0
}$$
Now we have a monad for $\ke_1$:
$$
0 \longrightarrow 2\ko(-1) \stackrel{A}{\longrightarrow} 7\ko
      \stackrel{B}{\longrightarrow} 2\ko(1) \longrightarrow \kc \longrightarrow 0    
$$
we dualize it ($B^\vee$ is no subbundle!) to get a monad for $\ke_1^\vee$
$$
0 \longrightarrow 2\ko(-1) \stackrel{B^\vee}{\longrightarrow} 7\ko
      \stackrel{A^\vee}{\longrightarrow} 2\ko(1) \longrightarrow 0
$$
and end up with a type-II monad
$$
 2\Omega^4(4) \stackrel{M^\vee}{\longrightarrow} 2\Omega^1(1)
      \stackrel{P^\vee}{\longrightarrow} 3\ko 
$$
where $M^\vee$ is no subbundle. 

As in the previous case, the left matrix is $M^\vee$ by the same argument 
using the naturality of the Serre duality and the right matrix is 
$P^\vee = \coker\,\tilde M$.

What remains is to prove the following 

\begin{prop} \label{ClassNprim} 
Let $(\langle M \rangle, \langle N \rangle)$ define a monad 
with a bundle as its cohomology with $\rank\,\tilde M = 7$ and 
$N_1 = \left( \frac{N}{x\ 0} \right) = \coker\,\tilde M$. Then: \\
(i) $\dim\, \Span\, N\,=3$ and $x \in \Span\, N$;\\
(ii) if $\kc = \kc oker ( 2\Omega^1(1) \stackrel{N_1}{\longrightarrow} 3\ko )$,
then $\dim\, \supp\, \kc\,=\,0$ and $\length\, \kc = 3$. 
\end{prop}

{\bf Proof of the proposition}
We do the proof in two steps. We begin with classifying all possible
matrices $N_1$ appearing as $\left( \frac{N}{x \ 0} \right)$ 
up to the natural action of the group $GL(2) \times GL(2) \times GL(5)$
under assumption, that $N$ generates no column and no row of the form 
$(\kappa_1 v, \kappa_2 v)$ for $v\in V,\, \kappa_1, \kappa_2\in k$ (cf. lemma
\ref{OpenN} and remark \ref{remOpenN}). 
Compare \ref{MotivClass} for the motivation.

Throughout the proof we denote 
$N_1 = \left( \begin{array}{cc} 
n_{11} & n_{12} \\ n_{21} & n_{22} \\ n_{31} & n_{32} 
\end{array} \right)$

{\bf Case 1.} $\dim \, \Span\, N = 4$

{\bf Case 1.1.} $x \in \Span\, N$

Obviously 
$$ N_{1} \sim \left( \begin{array}{cc} e_1 & e_2 \\ e_3 & e_4 \\ e_0 & 0 
\end{array} \right) \ =: \ T_1 $$ 

{\bf Case 1.2.}  $x \not\in \Span\, N$

{\it Case 1.2.1} $x, n_{1,2}, n_{2,1}, n_{2,2}$ form a basis
of $\Span\, N$. \\
Then we have the following sequence of transformations (explainations are given
after the diagram):
\scriptsize
$$
\xymatrix {
&&&
{\left( \begin{array}{cc} 
0 & e_2 \\ e_3 & e_4 \\ e_1 & 0
\end{array} \right)\ =: \ T_2}
\\
{\left( \begin{array}{cc} 
e_1+\lambda e_2 + \mu e_3 + \nu e_4 & e_2 \\ e_3 & e_4 \\ e_1 & 0
\end{array} \right)}
\ar@{~>}[r]
&
{\left( \begin{array}{cc} 
e_1 + \nu e_4 & e_2 \\ e_3 & e_4 \\ e_1 & 0
\end{array} \right)}
\ar@{~>}[r]
& 
{\left( \begin{array}{cc} 
\nu e_4 & e_2 \\ e_3 & e_4 \\ e_1 & 0
\end{array} \right)}
\ar@{~>}[ur]_{\nu = 0}
\ar@{~>}[dr]_{\nu \ne 0}
\\
&&&
{\left( \begin{array}{cc} 
e_4 & e_2 \\ e_3 & e_4 \\ e_1 & 0
\end{array} \right) \ =: \ T_3}
}$$
\normalsize
In this diagram, the first transformation from 
$\left( \begin{array}{cc} e_1+\lambda e_2 + \mu e_3 + \nu e_4 
  & e_2 \\ e_3 & e_4 \\ e_1 & 0 \end{array} \right)$ to 
$\left( \begin{array}{cc} e_1 + \nu e_4 & e_2 
     \\ e_3 & e_4 \\ e_1 & 0 \end{array} \right)$
 is done by substracting $\lambda$ times the 
second column from the first column, substracting $\mu$ times the second row 
from the first row and choosing new value for $\nu$;\\
then we substruct the third row from the first and get  
$\left( \begin{array}{cc} \nu e_4 & e_2 \\ e_3 & e_4 \\ e_1 & 0
  \end{array} \right)$.\\
If now $\nu=0$, we have already finished at the matrix $T_2$ from the diagram;
if $\nu\ne 0$, choose a new $e_1$ to be $\nu e_1$ and a new $e_3$ to 
be $\nu e_3$ and arrive at the matrix $T_3$. 

{\it Case 1.2.2.} Neither 
$x, n_{1,2}, n_{2,1}, n_{2,2}$ form a basis of $\Span\, N$, 
nor do $x, n_{1,1}, n_{1,2}, n_{2,2}$. \\
Then, without loss of geneality, 
$ x, n_{1,1}, n_{2,1}, n_{2,2} $ are a basis of $\Span\, N$.
$$ N_{1} \ = \ 
 \left( \begin{array}{cc} 
e_2 & e_1+\lambda e_2 + \mu e_3 + \nu e_4  \\ e_4 & e_3 \\ e_1 & 0
\end{array} \right) \ = \ 
 \left( \begin{array}{cc} 
e_2 & e_1+\lambda e_2  \\ e_4 & e_3 \\ e_1 & 0
\end{array} \right)
$$
because $\mu = \nu = 0$ due to our assumption. If $\lambda = 0$, we get
$$
   N_{1} \ = \
 \left( \begin{array}{cc} 
e_2 & e_1 \\ e_4 & e_3 \\ e_1 & 0
\end{array} \right) \ =: \ T_4
$$
otherwise:\\
add to the first row $1/\lambda$ times the third; \\
take $e_1+\lambda e_2$ to be a new $e_2$;\\
substract $1/\lambda$ times the second column from the first one; \\
and put $e_4-\frac{1}{\lambda} e_3$ to be a new $e_4$ finally arriving
at the matrix $\left( \begin{array}{cc} 
0 & e_2 \\ e_4 & e_3 \\ e_1 & 0 \end{array} \right) \sim T_2$,\\
and so we get nothing new.

The transformations that have been described just can be illustrated by the 
diagram: 
\scriptsize $$\xymatrix{
{ \left( \begin{array}{cc} 
e_2 & e_1+\lambda e_2  \\ e_4 & e_3 \\ e_1 & 0
\end{array} \right)}
\ar@{~>}[r]
&
{ \left( \begin{array}{cc} 
\frac{1}{\lambda} e_1 + e_2 & e_1+\lambda e_2  \\ e_4 & e_3 \\ e_1 & 0
\end{array} \right)}
\ar@{~>}[r]_{\text{\tiny \qquad \quad new $e_2$}}
&
{ \left( \begin{array}{cc} 
\frac{1}{\lambda} e_2 & e_2 \\ e_4 & e_3 \\ e_1 & 0
\end{array} \right)}
\ar@{~>}[r]_{\text{\tiny new $e_4$}\quad }
&
{ \left( \begin{array}{cc} 
0 & e_2 \\ e_4 & e_3 \\ e_1 & 0
\end{array} \right)
\sim T_2}
}$$ \normalsize
\\

{\bf Case 2.} $\dim \, \Span\, N = 3$.

{\bf Case 2.1.} $x \in \Span\, N$
Since we can bring $N_{1}$ to the form with 
$ n_{1,1} = n_{2,2}$
by elementary transformations  between first two rows, we get
$$ N_{1} \sim \left( \begin{array}{cc} e_1 & e_2 \\ e_3 & e_1 \\ e_0 & 0 
\end{array} \right) \ \sim \ T_3. $$  

{\bf Case 2.2.}  $x \not\in \Span\, N$.

First note that we cannot have 
$N_{1} \,\sim\, \left( \begin{array}{cc} w & 0 \\ w & 0 \\ u & v 
\end{array} \right)$
or 
$N_{1} \,\sim\, \left( \begin{array}{cc} w & 0 \\ 0 & w \\ u & v 
\end{array} \right)$
because such matrices do not generate a $2\times 2$-submatrix that 
can play the role of $N$. 

Denote by $S$ the set of points $v\in\P V$ such that $N_{1}$ has a 
generalized row of the form $(\alpha v, \beta v)$. By what we have remarked it 
is clear that $|S| \le 3$ for matrices $N_{1}$ of our type and the elements 
of $S$ are linearly independent.

If $S=\{e_1,e_2,e_3 \}$ then 
$N_{1} \,\sim\, \left( \begin{array}{cc} e_1 & 0 \\ 0 & e_2 \\ e_3 & e_3 
\end{array} \right) \ =: \ U_1$   

If $S=\{ e_1,e_2\}$, we get $N_{1} \,\sim\,
\left( \begin{array}{cc} e_1 & 0 \\ 0 & e_2 \\ 
     e_3+\lambda e_2 & \varepsilon e_3 + \mu ' e_1 \end{array} \right)
\, \sim \,
\left( \begin{array}{cc} e_1 & 0 \\ 0 & e_2 \\ 
     e_3 & \varepsilon  e_3 + \mu  e_1 \end{array} \right)$, 
where $\varepsilon=0,1$. We have $\mu \ne 0$ for otherwise $e_3\in S$, and
we also may assume $\varepsilon =0$ in order not to have 
$\varepsilon e_3 + \mu e_1 \in S$. We arrive at the matrix
$$
N_{1} \ \sim \left( \begin{array}{cc} e_1 & 0 \\ 0 & e_2 \\ e_3 & e_1 
\end{array} \right) \ =: \ U_2
$$

Finally, assume $S=\{e_1 \}$ and suppose that 
$ n_{3,1} = e_1,\  n_{3,2} = 0$. 
and that $n_{1,1} = n_{2,2}$ (this can be attained by elementary
transformations between the first two rows
Since $|S|=1$ we have $n_{2,1}$ is linearly independent of $e_1$. 
If now $n_{1,1} = n_{2,2} = \lambda e_1 + \mu e_2$ then
$N_{1} =\left( \begin{array}{cc} \lambda e_1 + \mu e_2 & e_3 
             \\ e_2 & \lambda e_1 + \mu e_2 \\ e_1 & 0 \end{array} \right) $   
and
$e_3 - \lambda (\lambda e_1 + \mu e_2) \in S$, contradiction.\\
Therefore we may suppose
$N_{1} =\left( \begin{array}{cc} e_3 & e_1 + \lambda e_2 + \mu e_3 
             \\ e_2 & e_3 \\ e_1 & 0 \end{array} \right)$. 
The condition that $(e_1,0)$ is the only degenerated generalized row of the matrix can be rewritten as an assumption that the equation 
$$
 (\, xe_3 + ye_2 + ze_1\, ) \,\wedge\, 
 (\, x(e_1+\lambda e_2 + \mu e_3 )\,+\, ye_3\, )\, = \, 0
$$
or, which is the same,
$$
\rank\, \left( \begin{array}{ccc} z & y & x \\ 
              x & \lambda x & \mu x + y  \end{array} \right) \ = \ 1
$$
has the only solution $(x,y,z) = (0,0,1) \in \P^2$. \\
If in this equation $x=0$, then also $y=0$ and we arrive at the solution
$(0,0,1)$. Otherwise put $x=1$ and get $y=\lambda z$ and 
$\lambda z^2 + \mu z - 1 = 0$. If $\lambda \ne 0$ or $\mu \ne 0$ then we have 
an additional solution. Hence $\lambda=\mu=0$ and
$$N_{1} \  \sim \ \left( \begin{array}{cc} e_3 & e_1 
             \\ e_2 & e_3 \\ e_1 & 0 \end{array} \right) \ =: \ U_3,
$$
that completes the classification.


We are starting now with the second part of the proof where
we shall see that for matrices $N_{1}$ of types $T_{1,2,3,4}$ there exist
no matrices $M$ defining a subbundle 
$2\Omega^4(4) \stackrel{M}{\longrightarrow} 2\Omega^1(1)$
and that for types $U_{1,2,3}$ they do exist.

Let $\Gamma \in Mat_{2\times (r+1)}(\bigwedge^3V) $ be a total matrix of 
syzygies of $N_{1}$, i.e. $N_{1} \Gamma $=0 and columns of $\Gamma$ 
generate all the syzygies of degree 3 of the matrix $N_{1}$. 
Then all possible matrices $M$ with the condition 
$N_{1} M = 0$ are of the form 
$$ 
M = G \cdot \left[ \left( \begin{array}{ccc} 
p_0, & \dots, & p_r \\
q_0, & \dots, & q_r 
\end{array} \right) ^T \right] 
$$

Here we shall several times use lemma \ref{ColWedge}.

If $N_{1}=T_1$, we get
$$
\Gamma = \left( \begin{array}{cccccc} 
0 & 0 & 0& e_{024} & e_{013} & e_{023} + e_{014} \\
e_{024} & e_{124} & e_{234} & e_{023} + e_{014} & 0 & e_{013}
\end{array}\right)
$$
\scriptsize $$
M = \left(
\begin{array}{cc}
 p_3 e_{024} + p_4 e_{013} + p_5 (e_{023} + e_{014}) & 
    q_3 e_{024} + q_4 e_{013} + q_5 (e_{023} + e_{014}) \\
 p_0 e_{024} + p_1 e_{124} + p_2 e_{234} 
             + p_3 (e_{023} + e_{014}) + p_5 e_{013} &
 q_0 e_{024} + q_1 e_{124} + q_2 e_{234} 
             + q_3 (e_{023} + e_{014}) + q_5 e_{013}  
\end{array}
\right)
$$ \normalsize
The exterior forms $\xi,\eta$ being as in lemma \ref{ColWedge}, we get:
$$ (\eta^{*2})^*\wedge\xi = - 2 (q_1 t + p_1 s) f e_{0123}
                            - 2 (q_2 t + p_2 s) f e_{0134} = 0 $$
$$ (\xi^{*2})^*\wedge\eta = - 2 (q_1 t + p_1 s) f e_{0124} 
                            - 2 (q_2 t + p_2 s) f e_{0234} = 0 $$
where
$$
f \ = \ ( q_3 q_4 - q_5^2) \cdot t^2 \ + \ 
  (p_4 q_3 + p_3 q_4 - 2 p_5 q_5) \cdot s t \ + \
  (p_3 p_4 - p_5^2 ) \cdot s^2
$$
and there always exists $(s,t)$ with  $f(s,t)=0$. 

If $N_{1}=T_2$, we get
$$
\Gamma = \left( \begin{array}{ccccccc} 
0 & 0 & 0 & e_{013} & e_{123} & e_{134} & e_{124} \\
e_{024} & e_{124} & e_{234} & 0 & 0 & 0 & e_{123}
\end{array}\right)
$$
\scriptsize $$
M = \left(
\begin{array}{cc}
p_3 e_{013} + p_4 e_{123} + p_5 e_{134} + p_6 e_{124} &
q_3 e_{013} + q_4 e_{123} + q_5 e_{134} + q_6 e_{124} \\
p_0 e_{024} + p_1 e_{124} + p_2 e_{234} + p_6 e_{123} &
q_0 e_{024} + q_1 e_{124} + q_2 e_{234} + q_6 e_{123}
\end{array}
\right)
$$ \normalsize
$$ (\eta^{*2})^*\wedge\xi =  2 (q_0 t+p_0s)(q_3 t + p_3 s)(q_6 t + p_6 s) 
                             \cdot e_{0123} 
                          -  2 (q_0 t+p_0s)(q_5 t + p_5 s)(q_6 t + p_6 s) 
                             \cdot e_{1234} = 0 $$
$$ (\xi^{*2})^*\wedge\eta =  2 (q_0 t+p_0s)(q_3 t + p_3 s)(q_6 t + p_6 s) 
                             \cdot  e_{0124} 
                           - 2 (q_2 t+p_2s)(q_3 t + p_3 s)(q_6 t + p_6 s) 
                             \cdot  e_{1234} = 0 $$
and the left hand sides of the equations have a common linear factor 
$(q_6 t + p_6 s)$.

In case $N_{1} = T_3$ 
$$
\Gamma = \left( \begin{array}{cccccc} 
  0 & 0 & 0 & e_{134} & -e_{123} & e_{124} \\
  e_{024} & e_{124} & e_{234} & 0 & e_{134} & e_{123} 
\end{array}\right)
$$
\scriptsize $$
M = \left( \begin{array}{cc} 
p_3 e_{134} - p_4 e_{123} + p_5 e_{124} &
q_3 e_{134} - q_4 e_{123} + q_5 e_{124} \\
p_0 e_{024} + p_1 e_{124} + p_2 e_{234} + p_4 e_{134} + p_5 e_{123} & 
q_0 e_{024} + q_1 e_{124} + q_2 e_{234} + q_4 e_{134} + q_5 e_{123} 
\end{array}\right)
$$ \normalsize
$$(\eta^{*2})^*\wedge\xi = f \cdot e_{1234} = 0 \quad ; 
  \quad (\xi^{*2})^*\wedge\eta = 0 $$
where $f$ is a polynomial of degree 3 in (s,t) and has a zero.

If $N_{1} = T_4$ then 
$$
\Gamma = \left( \begin{array}{ccccccc} 
0 & 0 & 0 & e_{014} & e_{124} & e_{134} & e_{123} \\
e_{014} & e_{124} & e_{134} & e_{013}+e_{024} & e_{123} & e_{234} & 0
\end{array}\right)
$$
\scriptsize $$
M = \left(
\begin{array}{cc}
p_3 e_{014} + p_4 e_{124} + p_5 e_{134} + p_6 e_{123} &
q_3 e_{014} + q_4 e_{124} + q_5 e_{134} + q_6 e_{123} \\
p_0 e_{014} + p_1 e_{124} + p_2 e_{134} + p_3 ( e_{013}+e_{024}) 
    + p_4 e_{123} + p_5 e_{234} &
q_0 e_{014} + q_1 e_{124} + q_2 e_{134} + q_3 ( e_{013}+e_{024}) 
    + q_4 e_{123} + q_5 e_{234} 
\end{array}
\right)
$$ \normalsize
$$ (\eta^{*2})^*\wedge\xi =  - 2 (q_6 t+p_6s)(q_3 t + p_3 s)^2
                             \cdot e_{0123} 
                          + 2 (q_6 t + p_6 s) F
                             \cdot e_{1234} = 0 $$
$$ (\xi^{*2})^*\wedge\eta = - 2 (q_6 t+p_6s)(q_3 t + p_3 s)^2
                             \cdot  e_{0124} 
                           + 2 (q_6 t+p_6s)(q_3 t + p_3 s)(q_5 t + p_5 s) 
                             \cdot  e_{1234} = 0 $$
where
$$
F\, =\, (q_2t-p_2s)(q_3t+p_3s)\, -\, (q_0t+p_0s)(q_5t+p_5s)  
$$
and the left hand sides of the equations have a common linear factor 
$(q_6 t + p_6 s)$.

We have the following examples of matrices $M$ defining subbundles 
$2\Omega^4(4) \stackrel{M}{\longrightarrow} 2\Omega^1(1)$ with 
$N_{1} \wedge M = 0$ for $N_{1}$ of the types $U_{1,2,3}$.

If $N_{1}=U_1$, we get 
$$\Gamma = \left( \begin{array}{cccccccc} 
0 & 0 & 0 & e_{013} & e_{123} & e_{134} & e_{012} & e_{124} \\  
e_{023} & e_{123} & e_{234} & 0 & 0 & 0 & -e_{012} & -e_{124}
\end{array}\right)
$$
and take, for example,
$$M=\left( \begin{array}{cc} 
e_{012}+e_{134} & e_{013}+e_{124}+e_{134} \\
-e_{012}+e_{234} & e_{023}-e_{124}        
\end{array} \right)  
$$

For $N_{1}=U_2$ we have the total matrix of syzygies 
$$\Gamma =\left( \begin{array}{cccccccc}  
e_{013} & e_{123} & e_{134} & e_{012} & e_{124} & 0 & 0 & 0 \\   
0 & 0 & 0 & -e_{023} & -e_{234} & e_{012} & e_{123} & e_{124}
\end{array} \right)  
$$
put
$$M=\left( \begin{array}{cc}  
e_{013}+e_{124} & e_{012}+e_{124}+e_{134} \\
e_{124}-e_{234} & -e_{023}-e_{234}     
\end{array} \right)  
$$

Finally, if $N_{1} = U_3$, then
$$\Gamma =\left( \begin{array}{cccccccc}  
0 & 0 & -e_{012} & 0 & e_{123} & -e_{124} & e_{013} & e_{134} \\
e_{013} & e_{134} & e_{023} & e_{123} & 0 & e_{234} & e_{012} & e_{124} 
\end{array} \right)  
$$
and an example of a M is
$$M=\left( \begin{array}{cc}  
-e_{012}+e_{013}-e_{124} & -e_{124}+e_{134} \\              
e_{012}+e_{023}+e_{234} & e_{013}+e_{123}+e_{124}+e_{234}
\end{array} \right) . 
$$

We see that $N_{1}$ allows a subbundle $M$ iff 
$r+1 = \rank\,\Gamma = 8$. 
But we always have an exact sequence of sheaves
$$
0 \longrightarrow (r+1)\Omega^4(5) \stackrel{\Gamma}{\longrightarrow} \Omega^1(2)
  \stackrel{N_{1}}{\longrightarrow} 3\ko(1) \longrightarrow \kc \longrightarrow 0
$$
where $\Gamma$ defines an injection of the two $\ko$-modules generated by their
global sections because it defines an injection of these global sections.
We obtain
$$\length\,\kc \, = \, \dim\, H^0(\kc) \, =\, 
   \dim\, H^0((r+1)\Omega^4(5)) - \dim\,H^0(\Omega^1(2)) + \dim\,H^0(3\ko(1)) \, = \, 
$$ $$
    =\, (r+1) - 2\cdot \dim\,\bigwedge^2 V + 3\cdot \dim\,V \, = \, 8-20+15 \, =\, 3,
$$
Q.E.D.

\newpage

\section{The Construction of the Moduli Space}
\label{formal}

Here we shall define the moduli problem we are working with
(\ref{ModuliProblem}) and show that this problem admits a coarse moduli space 
(theorem \ref{ExistModuli}).

In this chapter we are making use of the following general theory.

\begin{defin}
Let $X \subset \P^n$ be a quasi-projective variety, vhere $\P^n \cong \P V$,
where $V$ is an $(n+1)$-dimensional vector space. If $G$ is an affine 
algebraic group, acting on $X$, define a {\bf linearization} of this action
to be a group homomorphism $\varphi : G \to GL(V)$ such that 
$$
  g \circ \langle x \rangle \ = \ \langle \varphi(g) \circ x \rangle , 
$$
for any $g\in G\,,\, \langle x \rangle \in X$. 
\end{defin}

\begin{defin}
If $\lambda \in \kx_*(G)$ is a 1-parameter subgroup of $G$, it defines a 
1-parameter subgroup of $GL(V)$ by $\varphi \circ \lambda : k^* \to GL(V)$. 
As the image of 1-parameter subgroup is a torus and, in particular, 
diagonalizable, after a suitable basis change we have 
$\varphi \circ \lambda  (t) = Diag(t^{a_0},\dots,t^{a_n})$. For a 
vector $x\in V$ define
$$
  \mu(x,\lambda) := \min \{ a_k \ : \ x_k \ne 0 \} 
                  = \min \{\ m \ : \ \exists \lim_{t\to 0}\,[ t^{-m}
                                 \cdot \varphi\circ\lambda(t)\,(x)] \ \}
$$ 
where the second expression is independent on a diagonalization. 
\end{defin}

\begin{defin}
We call a point $x\in V$ {\bf stable} if $\mu(x,\lambda)<0$ for all 
nontrivial 1-parameter subgroups
$\lambda\in\kx_*(G)-\{ 1 \}$.
\end{defin}

\begin{defin}
We say that the (categorical) quotient $\pi : X \to Y$ (denoted: $Y=:X/G$) 
is the {\bf geometric quotient},
if the following conditions are satisfied:\\
$(i)$ \ $\pi$ is $G$-equivariant, surjective and affine; \\
$(ii)$ \ for any open subset $U\subseteq Y$ the ring homomorphism
     $\pi^\sharp : \ko_Y(U) \to \ko_X(\pi^{-1}U)^G$ is an isomorphism; \\
$(iii)$ \ for any closed $G$-equivariant subset $W\subseteq Y$ \ 
      the image $\pi(W)\subseteq  Y$ is closed;\\
$(iv)$ for every (closed) point $y\in Y$ the preimage $\pi^{-1}(y)$ is 
     precisely a $G$-orbit.
\end{defin}

\begin{theorem} \cite[theorem 3.14]{newstead} In the above 
situation, if $X$ is contained in the set
$\P^s$ of stable points with respect to some linearization $\varphi$ of the 
$G$-action, then there exists a geometric quotient $X/G$ which is again 
quasi-projective.   
\end{theorem}
 
And now we start applying this theory to our situation.

Consider $\P Hom(k^2,k^2 \otimes \bigwedge^3 V) \times 
\P Hom(k^2,k^2 \otimes V)$ with the action of the group
 $G = SL(2) \times  SL(2) \times  SL(2)$ given by 
$$ (g_1,g_2,g_3)(\langle M \rangle, \langle N \rangle) := 
  (\langle g_2^{-1} M g_1 \rangle , \langle g_3^{-1} M g_2 \rangle   ) .  $$

Recall the notation:
$$ X = \{ (\langle M \rangle, \langle N \rangle) : N \wedge M = 0 \}  
\subset \P Hom(k^2,k^2 \otimes \bigwedge^3 V) \times 
\P Hom(k^2,k^2 \otimes V), $$
$$ X_0 = \{ (\langle M \rangle, \langle N \rangle) \in X : M \ is \ subbundle,
\ N \ is \ surjective \} $$

We have the $G$-equivariant inclusions 
$$ X_{0} \subset X \subset \P Hom(k^2,k^2 \otimes \bigwedge^3 V) \times 
\P Hom(k^2,k^2 \otimes V) \stackrel{Segre}{\hookrightarrow} 
\P Hom(k^2,k^2 \otimes \bigwedge^3 V) \otimes 
Hom(k^2,k^2 \otimes V), $$
where we define the $G$-action on $\P Hom(k^2,k^2 \otimes \bigwedge^3 V) 
\otimes Hom(k^2,k^2 \otimes V)$ by 
$$ (g_1,g_2,g_3)(\langle M \otimes N \rangle) := 
  (\langle g_2^{-1} M g_1 \otimes g_3^{-1} M g_2 \rangle ). $$ 

\begin{prop} $X_0$ consists of stable points with respect to the last group 
action, 
i.e., $$X_0 \subset  \P \left( Hom(k^2,k^2 \otimes \bigwedge^3 V) \otimes 
Hom(k^2,k^2 \otimes V) \right) ^{s}$$
\end{prop}

{\bf Proof.}
In order to prove this, we are showing that $\mu (M \otimes N , \lambda ) < 0$
 for any nontrivial 1-parameter subgroup $\lambda \in \kx_{*}(G)$.
Any 1-parameter subgroup of $G$ is of the form 
$$
\lambda(t) = \left( U_1 \left( \begin{array}{cc} t^{r_1} & 0 \\ 0 & t^{-r_1} 
\end{array} \right) V_1 \ , \
U_2 \left( \begin{array}{cc} t^{r_2} & 0 \\ 0 & t^{-r_2} 
\end{array} \right) V_2 \ , \ 
U_3 \left( \begin{array}{cc} t^{r_3} & 0 \\ 0 & t^{-r_3} 
\end{array} \right) V_3 \right).
$$
Correspondently, the $k^*$-orbit of $(\langle M \rangle, \langle N \rangle)$ is 
then given by
\scriptsize
$$
(M, N) \stackrel{\lambda (t)}{\longmapsto}
\left( U_2^{-1} \left( \begin{array}{cc} t^{-r_2} & 0 \\ 0 & t^{r_2} 
                  \end{array} \right) 
\left( \begin{array}{cc} \alpha & \beta \\ \gamma & 
       \delta \end{array} \right) 
\left( \begin{array}{cc} t^{r_1} & 0 \\ 0 & t^{-r_1} 
       \end{array} \right) U_1 \ , \
U_3^{-1} \left( \begin{array}{cc} t^{-r_3} & 0 \\ 0 & t^{r_3} 
\end{array} \right) 
\left( \begin{array}{cc} x & y \\ z & w
       \end{array} \right) 
\left( \begin{array}{cc} t^{r_2} & 0 \\ 0 & t^{-r_2} 
\end{array} \right) U_2 \right)
$$
\normalsize
where 
$V_2^{-1} M V_1 = \left( \begin{array}{cc} \alpha & \beta \\ \gamma & 
                         \delta \end{array} \right) $ \ , \
$V_3^{-1} N V_2 = \left( \begin{array}{cc} x & y \\ z & w
                         \end{array} \right)  $.
Since $U_1$, $U_2$, $U_3$ are irrelevant for the computation of 
$\mu(( M,  N), \lambda)$, we shall omit them 
in the sequel. The coordinates of $\lambda (t) (M,N)$ in  
$ Hom(k^2,k^2 \otimes \bigwedge^3 V) \otimes Hom(k^2,k^2 \otimes V) $
will be the collection of $4 \cdot 4$ vectors of the form 
$$
(\text{entry of \ }M') \otimes (\text{entry of \ }N') 
    \in \bigwedge^3 V \otimes V \cong k^{50} ,
$$
where we abbreviate
$$ M' = \left( \begin{array}{cc} t^{-r_2} & 0 \\ 0 & t^{r_2} 
                  \end{array} \right) 
\left( \begin{array}{cc} \alpha & \beta \\ \gamma & 
       \delta \end{array} \right) 
\left( \begin{array}{cc} t^{r_1} & 0 \\ 0 & t^{-r_1} 
       \end{array} \right) \ ; \
N' = \left( \begin{array}{cc} t^{-r_3} & 0 \\ 0 & t^{r_3} 
\end{array} \right) 
\left( \begin{array}{cc} x & y \\ z & w
       \end{array} \right) 
\left( \begin{array}{cc} t^{r_2} & 0 \\ 0 & t^{-r_2} 
\end{array} \right) , $$
i.e. there will be $ 50 \cdot 4 \cdot 4 = 800 $ coordinates.
To keep notation easier, we shall think of these vectors in 
$\bigwedge^3 V \otimes V \cong k^{50}$ as of coordinates of the pair
$(M',N')$. 
Note that $\mu ((M,N),\lambda )=\mu (M,\lambda ) + \mu (N,\lambda ). $
Suppose $\mu ((M,N),\lambda ) \ge 0$. But 
$$N' = \left( \begin{array}{cc} x \cdot t^{r_2-r_3} & 
     y \cdot t^{-r_2-r_3} \\ z \cdot t^{r_2+r_3} 
      & w \cdot t^{-r_2+r_3} \end{array} \right).
$$
Since $N$ generates no row of the form $(\kappa_1 v, \kappa_2 v)$, 
all of the vectors $x,y,z,w$ are nonzero. If $\mu(N,\lambda)\ge 0$, then 
all the four expressions  
$(r_2-r_3),-(r_2-r_3),(r_2+r_3),-(r_2+r_3) \ \ge 0$, but this happens only if 
$r_2=r_3=0$. 
The same argument applied to $M$ shows that $r_1=r_2=0$, i.e. $\lambda$ is a 
trivial 1-parameter subgroup. So the claim follows.

By the theorem above we have a geometric quotient $X_0 / G =: \km$,
which is again a quasiprojective variety.

\begin{theorem} \label{ExistModuli}
The 29-dimensional irreducible quasi-projective 
variety $\km:=X_0 / G$ is a coarse moduli space 
(cf. definition \ref{ModuliSpaceDefin} and \ref{ModuliProblem}) 
of instanton bundles on $\P^4$ with quantum number 2 and of 
rank 4.
\end{theorem}

{\bf Proof.} Dimension and irreducibility follow immediately 
from the definition of $\km$, because $X_0$ is irreducible by 
theorem \ref{Xirred}.

The rest of this sections consists of showing 
that $\km$ fulfills the axioms of a coarse moduli space.

\begin{defin} \label{ModuliSpaceDefin}
Given a moduli problem defined by a contravariant functor 
$\Fgot: (Sch/k) \to (Sets)$, we call a scheme $\km$ a {\bf coarse moduli space} 
for our functor $\Fgot$, if it satisfies the following three axions:\\
$(i)$ There exists a natural transformations of functors $\dots$ 
$$
  \beta \ :\ \Fgot \longrightarrow \hom(-,\km)
$$
$(ii)$ \quad $\dots$ which is a bijection on a $k$-point:
$$
  \beta \ : \ \Fgot (Spec(k)) 
 \stackrel{1:1}{\longleftrightarrow} \text{(closed points of $\km$)} =
\hom(Spec(k),\km)
$$
$(iii)$ and universal in the sense that for any 
natural transformation of functors $\Fgot \to Mor(-,N)$, 
there exists a unique morphism $\psi : \km \to N$ such that the diagram 
commutes:
$$
\begin{array}{ccc}
  \Fgot & \stackrel{\beta}{\longrightarrow} & Mor(-,\km) \\
  & \searrow & \downarrow \lefteqn{Mor(-,\psi)} \\
  & & Mor(-,N). 
\end{array}
$$
\end{defin}

\begin{thought} \label{ModuliProblem} We are representing the functor
$$
  \Mgot (S) := \left\{ \kf : 
\begin{array}{c}
\kf \ \text{is a rank-4 vector bundle on \ } S \times \P^4 \\
\text{flat over \ } S \\
\text{and \ } \forall s \in S \ \text{closed \ } 
      \kf_{s} \ \text{is an instanton on \ } \P^4_s
\end{array}
\right\} \bigg{/} \sim
$$
where two such families are defined to be equivalent 
$$
   \kf \sim \kf' : \Leftrightarrow \ \kf'=\kf \otimes pr^*_S \kl 
$$
for a suitable line bundle $\kl$ on $S$.
\end{thought}

First we construct a natural transformation of functors 
$$
  \Mgot \longrightarrow \hom(-,\km)
$$
For  given $S$ we take a flat family $\bf{F}$ of locally free sheaves on $S$, i.e. a locally free sheaf on $S \times \P^4$ satisfying our conditions. $\bf{F}$ is then the cohomology of the relative Beilinson monad
$$
 {\bf K} \boxtimes \Omega^4(4) \longrightarrow  {\bf L} \boxtimes \Omega^1(1) 
        \longrightarrow  {\bf P} \boxtimes \ko   , 
$$ 
where ${\bf K}=R^3p_*({\bf F}\boxtimes\ko(-3)), 
{\bf L}=R^1p_*({\bf F}\boxtimes\ko(-1)), 
{\bf P}=R^1p_*({\bf F}\boxtimes\ko(-1))$ and $p$ is the product 
projection $S \times \P^4 \to S$. These three sheaves are locally free of rank 2
because for every $s\in S$
$$H^4p_*({\bf F}_s\boxtimes\ko(-3)) = 
H^2p_*({\bf F}_s\boxtimes\ko(-1)) =
H^2p_*({\bf F}_s\boxtimes\ko(-1)) = 0
$$

and we can apply the following theorem

\begin{theorem} \cite[III.12.11]{hartshorne} Let $P$ be a projective 
space, $p : S\times P \to S$ the product projection, and $\kf$ a coherent 
sheaf on $S\times P$ flat over $S$.  If the natural maps
$$(R^{j} p_* \kf)\otimes k(s) \stackrel{\sim}{\to} 
   H^j(P_s,\kf\otimes \ko_{P_s})$$
are surjective in all  points  $s\in S$,
then there are natural isomorphisms 
$$(R^{j-1} p_* \kf)\otimes k(s) \stackrel{\sim}{\to} 
   H^{j-1}(P_s,\kf\otimes \ko_{P_s})$$
for any point $s\in S$.
\end{theorem}

Over an open subset $U \subset S$ trivializing all of $\bf{K}, \bf{L}, \bf{P}$ we get 
$$
 2\ko_U \boxtimes \Omega^4(4) \stackrel{M}{\longrightarrow}  2\ko_U \boxtimes \Omega^1(1) 
        \stackrel{N}{\longrightarrow}  2\ko_U \boxtimes \ko   , 
$$ 
where the matrices $M,N$ depend on a point of $U$.
We get a composite map 
$$ 
  U \ni u \longmapsto (M(u),N(u)) \ \stackrel{\text{mod \ } G}{\longmapsto} 
  \ \km . 
$$
This composite map does not depend on the choice of $U$ and on a particular trivialization because the 
class $ (\langle M(u) \rangle , \langle N(u) \rangle) \ \text{mod \ } G$ is uniquely determined by the 
isomorphism class of $\bf{F}_u$. Therefore we get a global morphism $S \to \km$ (cf. below).

The next axiom:
$$
  \Mgot (pt) = \{ \text{isom classes of \ } \kf : \kf \text{\ is an instanton 
on\ } \P^4 \} 
 \stackrel{1:1}{\longleftrightarrow} \text{(closed points of $\km$)}
$$
To check this, we apply the \\ 

\begin{lemma} \cite[II.4.1.3]{OSS} Let $\ke, \ke'$ be the 
cohomology sheaves of two 
monads
$$
 M \ : \ 0 \longrightarrow \ka \longrightarrow \kb \longrightarrow \kc \longrightarrow 0 
$$ $$
 M' \ : \ 0 \longrightarrow \ka' \longrightarrow \kb' \longrightarrow \kc' \longrightarrow  0 
$$
over a smooth variety X. The maping 
$$
 \hom (M',M) \  \longrightarrow \hom(\ke,\ke')
$$
which associates to each homomorphism of monads the induced homomorphism of 
cohomology bundles is bijective if the following hypotheses are satisfied:
$$
 \hom(\kb,\ka') \ = \ \hom(\kc,\kb') \ = \ 0, 
$$ $$
 H^1(\kc^*\otimes\ka') \ = \ H^1(\kb^*\otimes\ka') \ = \ H^1(\kc^*\otimes\kb') \ = \ 
                             H^2(\kc^*\otimes\ka') \ = \ 0
$$ 
\end{lemma}

We put in the lemma
$$M \ = \ M' \ : 
\ 0 \longrightarrow 2\Omega^4(4) \longrightarrow \Omega^1(1) \longrightarrow 2\ko  \longrightarrow 0 
$$
the conditions become then \\
$(i)$ \ $\hom(\Omega^1(1),\Omega^4(4))=0, \hom(\ko,\Omega^1(1))=0$, that is 
obvious \\
$(ii)$ \ $H^1(\Omega^4(4))=0, H^1(\Omega^1(1))=0, H^2(\Omega^4(4))=0$, that is no less obvious \\
$(iii)$ \ $H^1\kh om(\Omega^1(1),\Omega^4(4))=0$, that can be proven by writing the isomorphisms
$$
H^1\kh om(\Omega^1(1),\Omega^4(4)) \cong H^1(\Omega^3(3) \otimes \ko (-1)) =  
H^1(\Omega^3(2) ) = 0 
$$
These conditions fulfilled, the isomorphism classes of monads of our type correspond 1:1 to the 
isomorphism classes of the cohomology sheaves.

{\it The 3rd axiom}. Given a natural transformation of functors $\Mgot \to Mor(-,N)$, we 
have to construct a morphism $\psi : \km \to N$ such that the following diagram commutes:
$$
\begin{array}{ccc}
  \Mgot & \stackrel{\beta}{\longrightarrow} & Mor(-,\km) \\
  & \searrow & \downarrow \lefteqn{Mor(-,\psi)} \\
  & & Mor(-,N) 
\end{array}
$$
By the properties of a quotient, constuction of $\psi$ is equivalent to construction of a $G$-equivariant 
morphism $\Psi: X_0 \to N$. 
We have the universal family $\bkf$  on $X_0$; take $\Psi$ to be its image under the map
 
$$
\begin{array}{cccc}
    \beta(X_0) \ : \ & \Mgot(X_0) & \longrightarrow & Mor(X_0,N) \\
                     & [\bkf] & \longmapsto & \Psi
\end{array}
$$
Since the diagram 
$$
\xymatrix{ \\  
    {[g^*\bkf]} \ar@{|->}@/^2pc/[rrr] \ar@{}[r]|{\in} & {\Mgot(X_0)} \ar[r] & Mor(X_0,N) \ar@{}[r]|{\ni}
              & {\Psi \circ g}       \\
    {[\bkf]} \ar@{}[r]|{\in} \ar@{|->}@/_2pc/[rrr] \ar@{|->}[u] & {\Mgot(X_0)} \ar[r] \ar[u]_{\Mgot(g)} & 
        Mor(X_0,N) \ar[u]_{Mor(g,N)} \ar@{}[r]|{\ni} & {\Psi} \ar@{|->}[u] 
\\
\\
}
$$
commutes, it suffices to show that $[g^*\kbf ]=[\kbf ]$, i.e. $g^*\kbf \cong \kbf \otimes pr^*_{X_0}\kl$ for a suitable line bundle $\kl$ on $X_0$.

It is clear that $(g^*\kbf)_x \cong \kbf_x$ for all $x \in X_0$. 
Further, $\hom_{\P^4}(\kbf_x,\kbf_x) \cong k$ for we know that automorphisms of $\kbf_x$ correspond to automorphisms of the monad, an the latters can be computed as follows. 

Suppose the diagram below is coomutative and the vertical arrows are isomorphisms:
$$
\xymatrix{ 
          {2\Omega^4(4)} \ar@{^{(}->}[r]^M & {2\Omega^1(1)} \ar[r]^{{e_1 \, e_2 \choose e_3 \, e_1}} 
  & {2\ko} 
\\
          {2\Omega^4(4)} \ar@{^{(}->}[r]^M \ar[u] & {2\Omega^1(1)} \ar[r]^{{e_1 \, e_2 \choose e_3 \, e_1}} 
\ar[u]_f & {2\ko} \ar[u]_h \\
} \\
$$
We have {\bf to show\ } that $g=h=\lambda \cdot Id$.\\
Write 
$$
f = \left( \begin{array}{cc} f_1 & f_2 \\ f_3 & f_4 \end{array} \right) \ ; \ 
h = \left( \begin{array}{cc} h_1 & h_2 \\ h_3 & h_4 \end{array} \right) . 
$$  
We get:
$$
\left( \begin{array}{cc} h_1 & h_2 \\ h_3 & h_4 \end{array} \right) \cdot
\left( \begin{array}{cc} e_1 & e_2 \\ e_3 & e_1 \end{array} \right) =
\left( \begin{array}{cc} e_1 & e_2 \\ e_3 & e_1 \end{array} \right) \cdot
\left( \begin{array}{cc} f_1 & f_2 \\ f_3 & f_4 \end{array} \right) ;
$$
$$
\left( \begin{array}{cc} h_1 e_1 + h_2 e_3 & h_1 e_2 + h_2 e_1 
       \\ h_3 e_1 + h_4 e_3 & h_3 e_2 + h_4 e_1 \end{array} \right) =
\left( \begin{array}{cc} f_1 e_1 + f_3 e_2 & f_2 e_1 + f_4 e_2 
       \\ f_1 e_3 + f_3 e_1 & f_2 e_3 + f_4 e_1 \end{array} \right).
$$
Comparison of the entries yields
$$
\begin{array}{rl}
 \text{ upper left :  }  & f_1=h_1, h_2=f_3=0   \\
 \text{ lower left :  }  & h_3=f_3=0, h_4=f_1=h_1=:\lambda  \\
 \text{ upper right : }  & f_4=h_1=\lambda, f_ 2=h_2=0
\end{array}
$$
i.e. 
$$ f = h = \left( \begin{array}{cc} \lambda & 0 \\ 0 & \lambda \end{array} \right) , $$
QED.\\

Define $\kl := \pi_*\kh om(\kbf, g^*\kbf)$ a line bundle on $X_0$, where $\pi=pr_{X_0}$.
Now we can check fiberwise that
$$
 \pi^*\kl \otimes \kbf = \left[ \pi^*\pi_*\kh om(\kbf, g^*\kbf) \right] \otimes \kbf \ 
    \stackrel{can \otimes id}{\longrightarrow} \kh om(\kbf, g^*\kbf) \otimes \kbf 
    \stackrel{ev}{\longrightarrow} g^*\kbf
$$
is an isomorphism.

So, we have cheched the three axioms and therefore realized the coarse moduli space as 
the geometric quotient quasi-projective variety. \qed

\newpage

\section{Smoothness of the Moduli Space}
\label{etale}

The fact we are proving in this section is:

\begin{theorem} \label{Msmooth} The moduli space $\km$ is smooth.
\end{theorem}

We start the proof by showing that $X_0$ is smooth.

\begin{theorem} \label{X0smooth}
$X_0$ is smooth.
\end{theorem}

{\bf Proof.}
We know that $X_0$ is 38-dimensional, i.e. has codimension 20 in 
$\P Mat_{2\times 2}(\bigwedge^3 V) \times \P Mat_{2\times 2}(V)$. 

Define the affine cone of $X_0$ to be 
$$ CX_0 = \{ (M,N) \, :\, M\ne 0 , \, N\ne 0, \, 
                       (\langle M \rangle, \langle N \rangle) \in X_0  \} 
 \ \subset \  Mat_{2\times 2}(\bigwedge^3 V) \times Mat_{2\times 2}(V)
$$
and analoguosly the affine cone $CX$ over $X$.

To prove that $X_0$ is smooth is the same as to prove that $CX$ is smooth 
at all points of $X_0$.

We shall write elements of 
$Mat_{2\times 2}(\bigwedge^2V)\times Mat_{2\times 2}(V)$ 
in the form 
$$\left( \ \left( \begin{array}{cc}
\sum_{i<j<k} A_{ijk} e_{ijk} & \sum_{i<j<k} B_{ijk} e_{ijk} \\
\sum_{i<j<k} C_{ijk} e_{ijk} & \sum_{i<j<k} D_{ijk} e_{ijk} 
\end{array} \right)\ , \ 
\left( \begin{array}{cc}
\sum_{i=0}^{4} a_i e_i & \sum_{i=0}^{4} b_i e_i \\
\sum_{i=0}^{4} c_i e_i & \sum_{i=0}^{4} d_i e_i  
\end{array} \right) \  \right)$$ 
and thus introduce coordinates 
$\{a_i, b_i, c_i, d_i, A_{ijk}, B_{ijk}, C_{ijk}, D_{ijk} \}$ 
in this 60-dimensional affine space. The equations of $CX$ can be written
in the matrix form as 
$$\begin{array}{l}
\left( \begin{array}{cc}
f_1\cdot e_{0123} +\dots +f_5\cdot e_{1234}   & 
     f_6\cdot e_{0123} +\dots +f_{10}\cdot e_{1234} \\
f_{11}\cdot e_{0123} +\dots +f_{15}\cdot e_{1234} & 
     f_{16}\cdot e_{0123} +\dots +f_{20}\cdot e_{1234}
\end{array} \right) \ = \\[15pt] \ = \ 
\left( \begin{array}{cc}
\sum_{i=0}^{4} a_i e_i & \sum_{i=0}^{4} b_i e_i \\
\sum_{i=0}^{4} c_i e_i & \sum_{i=0}^{4} d_i e_i  
\end{array} \right) \ \wedge \
\left( \begin{array}{cc}
\sum_{i<j<k} A_{ijk} e_{ijk} & \sum_{i<j<k} B_{ijk} e_{ijk} \\
\sum_{i<j<k} C_{ijk} e_{ijk} & \sum_{i<j<k} D_{ijk} e_{ijk} 
\end{array} \right) \ = \\[15pt] \ = \ 0 \ \in \ Mat_{2\times 2}(\bigwedge^4 V) 
\end{array}$$ 
This corresponds to 20 scalar equations $f_i(a_1,\dots,D_{234})=0$ for $i=1,\dots,20$.

The smoothness of $CX_0$ and hence of $X_0$ will be proven, if we show that the rank of the Jacobi
matrix 
$$
  J(M,N) = \left( \frac{\partial f_i}{\partial x} \right) (M,N)
$$
is 20 at any point $(M,N)\in CX_0$. Here $J(M,N)$ is a $20\times 60$-matrix, where 
$i \in \{1,\dots,20 \}$ and $x$ ranges through the 60 variables 
$\{ a_1,\dots,D_{234}\}$

The Leibnitz rule for matrices implies
$$\begin{array}{l}
\left( \begin{array}{cc}
\frac{\partial f_1}{\partial x}\cdot e_{0123} +\dots & \dots \\
\dots & \dots + \frac{\partial f_{20}}{\partial x}\cdot e_{1234}
\end{array} \right) \ = \\[15pt] 
\ = \ 
\frac{\partial\ }{\partial x}\left[ \, \left( \begin{array}{cc}
\sum_{i=0}^{4} a_i e_i & \sum_{i=0}^{4} b_i e_i \\
\sum_{i=0}^{4} c_i e_i & \sum_{i=0}^{4} d_i e_i  
\end{array} \right) \ \wedge \
\left( \begin{array}{cc}
\sum_{i<j<k} A_{ijk} e_{ijk} & \sum_{i<j<k} B_{ijk} e_{ijk} \\
\sum_{i<j<k} C_{ijk} e_{ijk} & \sum_{i<j<k} D_{ijk} e_{ijk} 
\end{array} \right) \, \right] \ = \\[15pt] 
\ = \ 
\frac{\partial\ }{\partial x}\left( \begin{array}{cc}
\sum_{i=0}^{4} a_i e_i & \sum_{i=0}^{4} b_i e_i \\
\sum_{i=0}^{4} c_i e_i & \sum_{i=0}^{4} d_i e_i  
\end{array} \right) \, \wedge \, M \, + \, N \, \wedge \,
\frac{\partial\ }{\partial x}
\left( \begin{array}{cc}
\sum_{i<j<k} A_{ijk} e_{ijk} & \sum_{i<j<k} B_{ijk} e_{ijk} \\
\sum_{i<j<k} C_{ijk} e_{ijk} & \sum_{i<j<k} D_{ijk} e_{ijk} 
\end{array} \right) .
\end{array}$$ 

To prove that $\rank\,J(M,N) = 20$ it suffices to check that $\rank\,J_1(M,N)=20$, 
where $J_1$ is a submatrix of $J$
$$
J_1(M,N) = \left( \frac{\partial f_i}{\partial x} \right)_{\tiny
\begin{array}{l}
1 \le i \le 20 \\
x \in \{A_{012},\dots,D_{234} \}
\end{array}
} (M,N)
$$

If $x \in \{A_{012},\dots,D_{234} \}$, the previous formula simplifies to
$$
\left( \begin{array}{cc}
\frac{\partial f_1}{\partial x}\cdot e_{0123} +\dots & \dots \\
\dots & \dots + \frac{\partial f_{20}}{\partial x}\cdot e_{1234}
\end{array} \right) \ = \ 
N \, \wedge \,
\frac{\partial\ }{\partial x}
\left( \begin{array}{cc}
\sum_{i<j<k} A_{ijk} e_{ijk} & \sum_{i<j<k} B_{ijk} e_{ijk} \\
\sum_{i<j<k} C_{ijk} e_{ijk} & \sum_{i<j<k} D_{ijk} e_{ijk} 
\end{array} \right) .
$$

By \ref{MotivClass} and \ref{ClassN} we may assume 
$N=\left( \begin{array}{cc} 
e_1 & e_2 \\ e_3 & e_4 \end{array} \right)$
or
$N=\left( \begin{array}{cc} 
e_1 & e_2 \\ e_3 & e_1 \end{array} \right)$.
In the former case $J_1$ is:
\tiny $$\left( \begin{array}{cccccccccccccccccccc}
0 &  0 &  0 &  0 &  0 &  0 &  0 &  0 &  0 &  0 &  -1 & 0 &  0 &  0 &  0 
&  0 &  0 &  0 &  0 &  0 \\
0 &  0 &  0 &  0 &  0 &  0 &  0 &  0 &  0 &  0 &  0 &  0 &  0 &  0 &  0 
&  0 &  0 &  0 &  0 &  0 \\ 
0 &  0 &  0 &  0 &  0 &  0 &  0 &  0 &  0 &  0 &  0 &  0 &  1 &  0 &  0 
&  0 &  0 &  0 &  0 &  0 \\  
-1 & 0 &  0 &  0 &  0 &  0 &  0 &  0 &  0 &  0 &  0 &  0 &  0 &  0 &  0 
&  0 &  0 &  0 &  0 &  0 \\ 
0 &  -1 & 0 &  0 &  0 &  0 &  0 &  0 &  0 &  0 &  0 &  0 &  0 &  1 &  0 
&  0 &  0 &  0 &  0 &  0 \\ 
0 &  0 &  -1 & 0 &  0 &  0 &  0 &  0 &  0 &  0 &  0 &  0 &  0 &  0 &  0 
&  0 &  0 &  0 &  0 &  0 \\ 
0 &  0 &  0 &  0 &  0 &  0 &  0 &  0 &  0 &  0 &  0 &  0 &  0 &  0 &  0 
&  0 &  0 &  0 &  0 &  0 \\ 
0 &  0 &  0 &  0 &  0 &  0 &  0 &  0 &  0 &  0 &  0 &  0 &  0 &  0 &  1 
&  0 &  0 &  0 &  0 &  0 \\ 
0 &  0 &  0 &  0 &  0 &  0 &  0 &  0 &  0 &  0 &  0 &  0 &  0 &  0 &  0 
&  0 &  0 &  0 &  0 &  0 \\ 
0 &  0 &  0 &  0 &  1 &  0 &  0 &  0 &  0 &  0 &  0 &  0 &  0 &  0 &  0 
&  0 &  0 &  0 &  0 &  0 \\ 
0 &  0 &  0 &  0 &  0 &  0 &  0 &  0 &  0 &  0 &  0 &  0 &  0 &  0 &  0 
&  -1 & 0 &  0 &  0 &  0 \\ 
0 &  0 &  0 &  0 &  0 &  0 &  0 &  0 &  0 &  0 &  0 &  0 &  0 &  0 &  0 
&  0 &  0 &  0 &  0 &  0 \\ 
0 &  0 &  0 &  0 &  0 &  0 &  0 &  0 &  0 &  0 &  0 &  0 &  0 &  0 &  0 
&  0 &  0 &  1 &  0 &  0 \\ 
0 &  0 &  0 &  0 &  0 &  -1 & 0 &  0 &  0 &  0 &  0 &  0 &  0 &  0 &  0 
&  0 &  0 &  0 &  0 &  0 \\ 
0 &  0 &  0 &  0 &  0 &  0 &  -1 & 0 &  0 &  0 &  0 &  0 &  0 &  0 &  0 
&  0 &  0 &  0 &  1 &  0 \\ 
0 &  0 &  0 &  0 &  0 &  0 &  0 &  -1 & 0 &  0 &  0 &  0 &  0 &  0 &  0 
&  0 &  0 &  0 &  0 &  0 \\ 
0 &  0 &  0 &  0 &  0 &  0 &  0 &  0 &  0 &  0 &  0 &  0 &  0 &  0 &  0 
&  0 &  0 &  0 &  0 &  0 \\ 
0 &  0 &  0 &  0 &  0 &  0 &  0 &  0 &  0 &  0 &  0 &  0 &  0 &  0 &  0 
&  0 &  0 &  0 &  0 &  1 \\ 
0 &  0 &  0 &  0 &  0 &  0 &  0 &  0 &  0 &  0 &  0 &  0 &  0 &  0 &  0 
&  0 &  0 &  0 &  0 &  0 \\ 
0 &  0 &  0 &  0 &  0 &  0 &  0 &  0 &  0 &  1 &  0 &  0 &  0 &  0 &  0 
&  0 &  0 &  0 &  0 &  0 \\ 
0 &  0 &  0 &  0 &  0 &  0 &  0 &  0 &  0 &  0 &  0 &  -1 & 0 &  0 &  0 
&  0 &  0 &  0 &  0 &  0 \\ 
1 &  0 &  0 &  0 &  0 &  0 &  0 &  0 &  0 &  0 &  0 &  0 &  -1 & 0 &  0 
&  0 &  0 &  0 &  0 &  0 \\ 
0 &  1 &  0 &  0 &  0 &  0 &  0 &  0 &  0 &  0 &  0 &  0 &  0 &  0 &  0 
&  0 &  0 &  0 &  0 &  0 \\ 
0 &  0 &  0 &  0 &  0 &  0 &  0 &  0 &  0 &  0 &  0 &  0 &  0 &  -1 & 0 
&  0 &  0 &  0 &  0 &  0 \\ 
0 &  0 &  0 &  0 &  0 &  0 &  0 &  0 &  0 &  0 &  0 &  0 &  0 &  0 &  0 
&  0 &  0 &  0 &  0 &  0 \\ 
0 &  0 &  0 &  -1 & 0 &  0 &  0 &  0 &  0 &  0 &  0 &  0 &  0 &  0 &  0 
&  0 &  0 &  0 &  0 &  0 \\ 
0 &  0 &  0 &  0 &  0 &  0 &  0 &  0 &  0 &  0 &  0 &  0 &  0 &  0 &  -1 
& 0 &  0 &  0 &  0 &  0 \\ 
0 &  0 &  0 &  0 &  0 &  0 &  0 &  0 &  0 &  0 &  0 &  0 &  0 &  0 &  0 
&  0 &  0 &  0 &  0 &  0 \\ 
0 &  0 &  0 &  0 &  -1 & 0 &  0 &  0 &  0 &  0 &  0 &  0 &  0 &  0 &  0 
&  0 &  0 &  0 &  0 &  0 \\ 
0 &  0 &  0 &  0 &  0 &  0 &  0 &  0 &  0 &  0 &  0 &  0 &  0 &  0 &  0 
&  0 &  0 &  0 &  0 &  0 \\  
0 &  0 &  0 &  0 &  0 &  0 &  0 &  0 &  0 &  0 &  0 &  0 &  0 &  0 &  0 
&  0 &  -1 & 0 &  0 &  0 \\ 
0 &  0 &  0 &  0 &  0 &  1 &  0 &  0 &  0 &  0 &  0 &  0 &  0 &  0 &  0 
&  0 &  0 &  -1 & 0 &  0 \\ 
0 &  0 &  0 &  0 &  0 &  0 &  1 &  0 &  0 &  0 &  0 &  0 &  0 &  0 &  0 
&  0 &  0 &  0 &  0 &  0 \\ 
0 &  0 &  0 &  0 &  0 &  0 &  0 &  0 &  0 &  0 &  0 &  0 &  0 &  0 &  0 
&  0 &  0 &  0 &  -1 & 0 \\ 
0 &  0 &  0 &  0 &  0 &  0 &  0 &  0 &  0 &  0 &  0 &  0 &  0 &  0 &  0 
&  0 &  0 &  0 &  0 &  0 \\ 
0 &  0 &  0 &  0 &  0 &  0 &  0 &  0 &  -1 & 0 &  0 &  0 &  0 &  0 &  0 
&  0 &  0 &  0 &  0 &  0 \\ 
0 &  0 &  0 &  0 &  0 &  0 &  0 &  0 &  0 &  0 &  0 &  0 &  0 &  0 &  0 
&  0 &  0 &  0 &  0 &  -1 \\ 
0 &  0 &  0 &  0 &  0 &  0 &  0 &  0 &  0 &  0 &  0 &  0 &  0 &  0 &  0 
&  0 &  0 &  0 &  0 &  0 \\  
0 &  0 &  0 &  0 &  0 &  0 &  0 &  0 &  0 &  -1 & 0 &  0 &  0 &  0 &  0 
&  0 &  0 &  0 &  0 &  0 \\ 
0 &  0 &  0 &  0 &  0 &  0 &  0 &  0 &  0 &  0 &  0 &  0 &  0 &  0 &  0 
&  0 &  0 &  0 &  0 &  0 
\end{array} \right) $$ \normalsize
and in the latter case 
$$
J_1 =\tiny \left( \begin{array}{cccccccccccccccccccc}
0&0&0&0&0&0&0&0&0&0&-1&0&0&0&0&0&0&0&0&0\\
0&0&0&0&0&0&0&0&0&0&0&0&0&0&0&0&0&0&0&0\\
0&0&0&0&0&0&0&0&0&0&0&0&1&0&0&0&0&0&0&0\\
-1&0&0&0&0&0&0&0&0&0&0&0&0&0&0&0&0&0&0&0\\
0&-1&0&0&0&0&0&0&0&0&0&0&0&1&0&0&0&0&0&0\\
0&0&-1&0&0&0&0&0&0&0&0&0&0&0&0&0&0&0&0&0\\
0&0&0&0&0&0&0&0&0&0&0&0&0&0&0&0&0&0&0&0\\
0&0&0&0&0&0&0&0&0&0&0&0&0&0&1&0&0&0&0&0\\
0&0&0&0&0&0&0&0&0&0&0&0&0&0&0&0&0&0&0&0\\
0&0&0&0&1&0&0&0&0&0&0&0&0&0&0&0&0&0&0&0\\
0&0&0&0&0&0&0&0&0&0&0&0&0&0&0&-1&0&0&0&0\\
0&0&0&0&0&0&0&0&0&0&0&0&0&0&0&0&0&0&0&0\\
0&0&0&0&0&0&0&0&0&0&0&0&0&0&0&0&0&1&0&0\\
0&0&0&0&0&-1&0&0&0&0&0&0&0&0&0&0&0&0&0&0\\
0&0&0&0&0&0&-1&0&0&0&0&0&0&0&0&0&0&0&1&0\\
0&0&0&0&0&0&0&-1&0&0&0&0&0&0&0&0&0&0&0&0\\
0&0&0&0&0&0&0&0&0&0&0&0&0&0&0&0&0&0&0&0\\
0&0&0&0&0&0&0&0&0&0&0&0&0&0&0&0&0&0&0&1\\
0&0&0&0&0&0&0&0&0&0&0&0&0&0&0&0&0&0&0&0\\
0&0&0&0&0&0&0&0&0&1&0&0&0&0&0&0&0&0&0&0\\
0&0&0&0&0&0&0&0&0&0&0&0&0&0&0&0&0&0&0&0\\
1&0&0&0&0&0&0&0&0&0&0&0&0&0&0&0&0&0&0&0\\
0&1&0&0&0&0&0&0&0&0&0&0&0&0&0&0&0&0&0&0\\
0&0&0&0&0&0&0&0&0&0&-1&0&0&0&0&0&0&0&0&0\\
0&0&0&0&0&0&0&0&0&0&0&-1&0&0&0&0&0&0&0&0\\
0&0&0&-1&0&0&0&0&0&0&0&0&-1&0&0&0&0&0&0&0\\
0&0&0&0&0&0&0&0&0&0&0&0&0&0&0&0&0&0&0&0\\
0&0&0&0&0&0&0&0&0&0&0&0&0&0&0&0&0&0&0&0\\
0&0&0&0&-1&0&0&0&0&0&0&0&0&0&0&0&0&0&0&0\\
0&0&0&0&0&0&0&0&0&0&0&0&0&0&1&0&0&0&0&0\\
0&0&0&0&0&0&0&0&0&0&0&0&0&0&0&0&0&0&0&0\\
0&0&0&0&0&1&0&0&0&0&0&0&0&0&0&0&0&0&0&0\\
0&0&0&0&0&0&1&0&0&0&0&0&0&0&0&0&0&0&0&0\\
0&0&0&0&0&0&0&0&0&0&0&0&0&0&0&-1&0&0&0&0\\
0&0&0&0&0&0&0&0&0&0&0&0&0&0&0&0&-1&0&0&0\\
0&0&0&0&0&0&0&0&-1&0&0&0&0&0&0&0&0&-1&0&0\\
0&0&0&0&0&0&0&0&0&0&0&0&0&0&0&0&0&0&0&0\\
0&0&0&0&0&0&0&0&0&0&0&0&0&0&0&0&0&0&0&0\\
0&0&0&0&0&0&0&0&0&-1&0&0&0&0&0&0&0&0&0&0\\
0&0&0&0&0&0&0&0&0&0&0&0&0&0&0&0&0&0&0&1
\end{array} \right) ,
$$ \normalsize
and we are done.

The rest of the work has to be postponed until the end of this section where 
all its ingredients will have been prepared. 

\begin{defin}
We call a morphism of schemes $f : S \to T$ {\bf locally isotrivial} with
fiber $F$, if there exists an open 
covering  $\ku = \{ U_i \}_{i \in I}$ of $T$ such that for each 
$i \in I$ there exists a diagram
$$
\xymatrix{
  {{\tilde U}_i \times F} \ar[r] \ar[d]^{pr_{{\tilde U}_i}} \ar@{}[rd]|{\Box} & 
    {U_i \times_{T} S} \ar[d]^{pr_{U_i}} \ar@{^{(}->}[r] & 
    {S} \ar[d]^{f}     \\
  {{\tilde U}_i} \ar[r]_{\text{\'etale}} & 
             {U_i} \ar@{^{(}->}[r]_{\text{open}} & {T}}
$$   
\end{defin}

The following results are presented in this form in \cite[2.5]{holger}

\begin{defin} An action of an algebraic group $G$ on a variety $X$ is 
called {\bf principal}, if 
the induced morphism is a closed immersion:
$$
\begin{array}{crcl}
   \omega \ : \  & G\times_k X & \longrightarrow & X_0 \times_k X_0 \\ 
    & (g,x) & \longmapsto & (y,g\circ x) 
\end{array}
$$
\end{defin}

\begin{defin} Let $X$ and $Y$ be separated algebraic varieties and let an 
affine algebraic 
group $G$ operate on $X$ with the geometric quotient morphism $\pi:X\to Y$. If $\pi$ is flat and the induced morphism
$$
\begin{array}{crclcc}
   \omega \ : \  & G\times_k X & 
     \stackrel{\cong}{\longrightarrow} & X \times_Y X & \subset & X \times_k X \\ 
    & (g,x) & \longmapsto & (y,g\circ x) 
\end{array}
$$
is an isomorphism, then  we call $(X,f,Y;F)$ a {\bf Mumford principal 
fibration}.
\end{defin}

\begin{theorem}
Let $X$ and $Y$ be separated algebraic varieties and let an affine algebraic group $G$ operate 
principally
on $X$ such that there exists a morphism $\pi:X\to Y$ making $Y$ into the geometrical quotient $X/G$. Then $(X,\pi,Y;G)$ is a Mumford principal fibration.
\end{theorem}

\begin{theorem}
Every Mumford principal fibration is locally isotrivial.
\end{theorem}

\begin{prop} \label{IsotrivProp}
The quotient map $X_0 \to X_0/G$ is a Mumford principial fibration
and, in particular, locally isotrivial.
\end{prop}

{\bf Proof} We have to show that: \\
$(i)$ \ $X_0, X_0/G$ are separated, what is true because they are both 
quasiprojective; \\
$(ii)$ \ $PGL(2)\times PGL(2)\times PGL(2)$ is an affine algebraic group. 
This is true because 
$PGL(2) \cong SL(2)/\{\pm\Id\}$ and there is the following theorem:

\begin{theorem} \cite[5.2.5]{Springer}
Let $H$ be a closed normal subgroup of the linear algebraic group $G$. Then $G/H$ with the 
induced group structure is a linear algebraic group.
\end{theorem}
 
$(iii)$ \ the action of $G$ on $X_0$ is principal, i.e. 
$$
\begin{array}{crcl}
   \omega \ : \  & PG\times_k X_0 & \longrightarrow & X_0 \times_k X_0 \\ 
    & (g,x) & \longmapsto & (y,g\circ x) \\
    & (f,g,h)(\langle M \rangle, \langle N \rangle) & \longmapsto &
       ((\langle M \rangle, \langle N \rangle),
           (\langle gMf^{-1} \rangle, \langle hNg^{-1} \rangle))  
\end{array}
$$
is an immersion, where $PG = PGL(2) \times PGL(2) \times PGL(2)$. \\
{\bf Step 1} We shall show first, that 
$$
\begin{array}{crcl}
   \omega_Y \ : \  & PG'\times_k Y_0 & \longrightarrow & Y_0 \times_k Y_0 \\ 
    & (g,y) & \longmapsto & (y,g\circ y) \\
    & (f,g,h)(\langle M \rangle, \langle N \rangle) & \longmapsto &
       ((\langle M \rangle, \langle N \rangle),
           (\langle gMf^{-1} \rangle, \langle hNg^{-1} \rangle))  
\end{array}
$$
where $PG' = PGL(2) \times PGL(2)$ is quotient of $G$ modulo the third factor,
and 
$$Y_0 = \{ \langle N \rangle \in \P Mat_{2\times 2}(V) : 
         2\Omega^1(1) \stackrel{N}{\longrightarrow} 2\ko \to 0 \} . $$ 

{\it Set-theoretical injectivity.} 
Suppose 
$\left\langle \left( 
\begin{array}{cc} g_1 & g_2 \\ g_3 & g_4 \end{array}
\right) \right\rangle , 
\left\langle \left( 
\begin{array}{cc} h_1 & h_2 \\ h_3 & h_4 \end{array}
\right) \right\rangle \in PGL(2) 
$ and
$$\left( \begin{array}{cc} g_1 & g_2 \\ g_3 & g_4 \end{array} \right)
\left( \begin{array}{cc} e_1 & e_2 \\ e_3 & e_1 \end{array} \right)
\left( \begin{array}{cc} h_1 & h_2 \\ h_3 & h_4 \end{array} \right) =
\lambda \cdot
\left( \begin{array}{cc} e_1 & e_2 \\ e_3 & e_1 \end{array} \right) \ , \
\lambda \in k^*. 
$$
The comparison of the coefficients standing by $e_2$ and $e_3$ gives the matrix
equalities:
$$
e_2 \ : \left( 
\begin{array}{cc} g_1 h_3 & g_1 h_4 \\ g_3 h_3 & g_3 h_4 
  \end{array} \right) = 
\left( \begin{array}{cc} 0 & \lambda \\ 0 & 0 \end{array} \right) \quad ; 
\quad e_3 \ : \left( 
\begin{array}{cc} g_2 h_1 & g_2 h_2 \\ g_4 h_1 & g_4 h_2 
  \end{array} \right) = 
\left( \begin{array}{cc} 0 & 0 \\ \lambda & 0 \end{array} \right) \ , 
$$
hence $g_1 h_4 = \lambda, g_3 = h_3 = 0$ and 
$g_4 h_1 = \lambda, g_2 = h_2 = 0$. And now the coefficients of $e_1$ give 
$$\left( \begin{array}{cc} g_1 & 0 \\ 0 & g_4 \end{array} \right)
\left( \begin{array}{cc} \lambda / g_4 & 0 
       \\ 0 & \lambda / g_1 \end{array} \right) = 
\left( \begin{array}{cc} \lambda & 0 \\ 0 & \lambda \end{array} \right)
$$
and we deduce $g_1/g_4=1$, that means 
$$
\left\langle \left( 
\begin{array}{cc} g_1 & g_2 \\ g_3 & g_4 \end{array}
\right) \right\rangle = 
\left\langle \left( 
\begin{array}{cc} h_1 & h_2 \\ h_3 & h_4 \end{array}
\right) \right\rangle = \langle \Id \rangle .
$$

The case of 
$N=\left( \begin{array}{cc} e_1 & e_2 \\ e_3 & e_4 \end{array} \right)$
can be done similarly.

{\it Smoothness.} It suffices to check that for each $N$ the map of the tangent spaces has the full rank equal to $\dim Y_0 + \dim G = \dim Y_0 + 6$:
$$
\begin{array}{crcl}
\tau \ : \ & T_{(\langle N \rangle , \Id ,\Id) }(Y_0 \times G') & 
\longrightarrow & T_{(\langle N \rangle), \langle N \rangle}(Y_0 \times Y_0) \\
& (N_1, \gg_1, \gh_1) & \longmapsto & (N_1, N_1+\gg_1N-N\gh_1),
\end{array}
$$
where $N_1 \in Mat_{2\times 2}(V);\ \gg_1,\gh_1 \in \sl$. 
Again, we do the proof for 
$N=\left( \begin{array}{cc} e_1 & e_2 \\ e_3 & e_1 \end{array} \right)$ 
and the other  case is considered even easier.

Clearly it suffices to check that
$$
\begin{array}{crclc}
\rank (  & (\sl \oplus \sl) & \longrightarrow & Mat_{2 \times 2}(V) & ) = 6 \\
  & (\gg_1, \gh_1) & \longmapsto & \gg_1 N - N \gh_1
\end{array}
$$
Studying the images of the basis vectors of $\sl \oplus \sl$ under this map, 
we get:
\footnotesize
$$\begin{array}{cccccccccccccc}
\left( \begin{array}{cc} 1 & 0 \\ 0 & -1 \end{array} \right) &
\left( \begin{array}{cc} e_1 & e_2 \\ e_3 & e_1 \end{array} \right) & = &
\left( \begin{array}{cc} e_1 & e_2 \\ -e_3 & -e_1 \end{array} \right)
 &&& ; &&& 
\left( \begin{array}{cc} e_1 & e_2 \\ e_3 & e_1 \end{array} \right) &  
\left( \begin{array}{cc} 1 & 0 \\ 0 & -1 \end{array} \right) & = &
\left( \begin{array}{cc} e_1 & -e_2 \\ e_3 & -e_1 \end{array} \right) 
\\  \\
\left( \begin{array}{cc} 0 & 1 \\ 0 & 0 \end{array} \right) &
\left( \begin{array}{cc} e_1 & e_2 \\ e_3 & e_1 \end{array} \right) & = &
\left( \begin{array}{cc} e_3 & e_1 \\ 0 & 0 \end{array} \right) 
&&& ; &&& 
\left( \begin{array}{cc} e_1 & e_2 \\ e_3 & e_1 \end{array} \right) & 
\left( \begin{array}{cc} 0 & 1 \\ 0 & 0 \end{array} \right) & = & 
\left( \begin{array}{cc} 0 & e_1 \\ 0 & e_3 \end{array} \right) 
\\  \\
\left( \begin{array}{cc} e_1 & e_2 \\ e_3 & e_1 \end{array} \right) &
\left( \begin{array}{cc} 0 & 0 \\ 1 & 0 \end{array} \right) & = & 
\left( \begin{array}{cc} 0 & 0 \\ e_1 & e_2 \end{array} \right) 
&&& ; &&& 
\left( \begin{array}{cc} e_1 & e_2 \\ e_3 & e_1 \end{array} \right) &
\left( \begin{array}{cc} 0 & 0 \\ 1 & 0  \end{array} \right) &  = & 
\left( \begin{array}{cc} e_2 & 0 \\ e_1 & 0 \end{array} \right)   
\end{array}
$$
\normalsize 
Now choose as a bassi in $Mat_{2\times 2}(V)$ he matrices
$$
\left( \begin{array}{cc} e_i & 0 \\ 0 & 0 \end{array} \right) \ , \
\left( \begin{array}{cc} 0 & e_i \\ 0 & 0 \end{array} \right) \ , \
\left( \begin{array}{cc} 0 & 0 \\ e_i & 0 \end{array} \right) \ , \
\left( \begin{array}{cc} 0 & 0 \\ 0 & e_i \end{array} \right) \ , \
i=\overline{0,4}
$$
and write the submatrix of the matrix of $\tau$ corresponding to the values 
$i=1,2,3$.
$$
\begin{array}{c || r r r | r r r | r r r | r r r ||}
\text{{ \tiny 
$ \begin{pmatrix} 1 & 0 \\ 0  & -1 \end{pmatrix} \cdot - $ 
}}
& 1 & 0 & 0 & 0 & 1 & 0 & 0 & 0 & -1 & -1 & 0 & 0 
\\
\text{{ \tiny 
$ \begin{pmatrix} 0 & 1 \\ 0  & 0 \end{pmatrix} \cdot - $ 
}}
& 0 & 0 & 1 & 1 & 0 & 0 & 0 & 0 & -0 & -0 & 0 & 0 
\\
\text{{ \tiny 
$ \begin{pmatrix} 0 & 0 \\ 1  & 0 \end{pmatrix} \cdot - $ 
}}
& 0 & 0 & 0 & 0 & 0 & 0 & 1 & 0 & -0 & 0 & 1 & 0 
\\
\text{{ \tiny 
$ - \cdot \begin{pmatrix} 1 & 0 \\ 0  & -1 \end{pmatrix} $ 
}}
& 1 & 0 & 0 & 0 & -1 & 0 & 0 & 0 & 1 & -1 & 0 & 0 
\\
\text{{ \tiny 
$ - \cdot \begin{pmatrix} 0 & 1 \\ 0  & 0 \end{pmatrix} $ 
}}
& 0 & 0 & 0 & 1 & 0 & 0 & 0 & 0 & 0 & 0 & 0 & 1 
\\
\text{{ \tiny 
$ - \cdot \begin{pmatrix} 0 & 0 \\ 1  & 0 \end{pmatrix} $ 
}}
& 0 & 1 & 0 & 0 & 0 & 0 & 1 & 0 & 0 & 0 & 0 & 0 
\\ 
\multicolumn{1}{c}{} & e_1 & e_2 & \multicolumn{1}{c}{e_3} 
 & e_1 & e_2 & \multicolumn{1}{c}{e_3} & e_1 & e_2 & \multicolumn{1}{c}{e_3} 
 & e_1 & e_2 & \multicolumn{1}{c}{e_3} \\
\multicolumn{1}{c}{} & \multicolumn{3}{c} 
{\text{{ \tiny $ \begin{pmatrix} \bullet & \cdot \\ \cdot  
& \cdot \end{pmatrix}$ }}}
 & \multicolumn{3}{c} 
{\text{{ \tiny $ \begin{pmatrix} \cdot & \bullet \\ \cdot  & \cdot 
\end{pmatrix}$ }}}
 & \multicolumn{3}{c} 
{\text{{ \tiny $ \begin{pmatrix} \cdot & \cdot \\ \bullet  & \cdot 
\end{pmatrix}$ }}}
 & \multicolumn{3}{c} 
{\text{{ \tiny $ \begin{pmatrix} \cdot & \cdot \\ \cdot  & \bullet 
\end{pmatrix}$ }}} 
\end{array} =: A
$$
whose rank is obviously 6.

{\bf Step 2.} Now we return to our $\omega$. \\
By lemma \ref{span3}, if
$(\langle M \rangle, \langle N \rangle) \in X_0$ 
and $\rank\, \Span\, M = 3$, then up to elementary transformations of colums
$M$ is of the form $\left( \begin{array}{cc} a & b \\ c
& a \end{array} \right)$, where $a,b,c$ are linearly independent in 
$\bigwedge^3 V$.\\

{\it Set-theoretical injectivity} follows by step 1 applied to the equalities 
$\langle fMg^{-1} \rangle = \langle M \rangle$ and \linebreak
$\langle gNh^{-1} \rangle = \langle N \rangle$.

{\it Smoothness.} We have to check that the map of tangent spaces 
$$
\begin{array}{rcl}
 T_{((\langle M \rangle,\langle N \rangle) , (\Id ,\Id, \Id)) }(X_0 \times G') & 
\longrightarrow & T_{(\langle M \rangle,\langle N \rangle),
                     (\langle M \rangle,\langle N \rangle))}(X_0 \times X_0) \\
 ((M_1,N_1), (\gf_1,\gg_1,\gh_1)) & \longmapsto & ((M_1,N_1),
                    (M_1+\gf_1M-M\gg_1,N_1+\gg_1N-N\gh_1)) ,
\end{array}
$$
has the full rank, where $\gf_1,\gg_1,\gh_1 \in \sl , 
M_1 \in Mat_{2 \times 2}(\bigwedge^3 V),  
N_1 \in Mat_{2 \times 2}(V)$. \\  
Clearly it suffices to check that the linear map
$$
 (\gf_1,\gg_1,\gh_1))  \longmapsto  (\gf_1M-M\gg_1, \gg_1N-N\gh_1)
$$
has rank 9. But this map has a (generalized sub)matrix of the form
$$
\left( \
\begin{array}{|cc|cc|}
 \cline{1-2} 
   && \multicolumn{2}{c}{\quad 0_{3 \times 12}} \\
  \cline{3-4} 
  \multicolumn{2}{|c|}{\text{\LARGE ${}^{A}$}} &  &  \\
 \cline{1-2} 
  \multicolumn{2}{c}{\quad 0_{3 \times 12}}  
     & \multicolumn{2}{|c|}{\text{\LARGE ${}^{A}$}}   \\
 \cline{3-4}
 \end{array} \
\right) 
$$

which is of rank 9 as A is of rank 6.\\
Hence the proposition. \qed

{\bf Proof of the theorem \ref{Msmooth} .}
By the previous proposition the quotient projection $X_0 \to X_0 / G$ 
is locally isotrivial, that means by definition, that there exists an open 
covering  $\ku = \{ U_i \}_{i \in I}$ of $X_0 / G$ such that for each 
$i \in I$ there exists a diagram
$$
\xymatrix{
  {{\tilde U}_i \times G} \ar[r] \ar[d]^{pr_{{\tilde U}_i}} \ar@{}[rd]|{\Box} & 
    {U_i \times_{X_0/G} X_0} \ar[d]^{pr_{U_i}} \ar@{^{(}->}[r] & 
    {X_0} \ar[d]^{\text{mod\ }G}     \\
  {{\tilde U}_i} \ar[r]_{\text{\'etale}} & {U_i} \ar@{^{(}->}[r]_{\text{open}} & {X_0 / G}   
}
$$
Now we have the following chain of equivalences :
\begin{eqnarray*}
 X_0/G \text{\ is smooth } & \Leftrightarrow & \\
 \text{all\ } U_i \text{\ are smooth} & \Leftrightarrow & 
          \text{\it \footnotesize ( by \cite[VI, 4.9, p.117]{AK} )} \\ 
 {\tilde U}_i \text{\ are smooth} & \Leftrightarrow & \\
 {\tilde U}_i \times G \text{\ are smooth} & \Leftrightarrow & 
           \text{\it \footnotesize ( by loc cit, because 
   ${\tilde U}_i \times G \to {U_i \times_{X_0/G} X_0}$ is \'etale, too )}   \\
  U_i \times G \text{\ are smooth} & \Leftrightarrow & \\
  X_0 \text{\ is smooth,} 
\end{eqnarray*}
and the last statement was the contents of theorem \ref{X0smooth}.  \qed

\newpage

\section{The Jumping Lines}
\label{jump}

As the first Chern class of our bundles is 0, we may expect that restricting
an instanton bundle $\ke$ on a generic line $\ell \subset \P^4$ gives
$\ke|_\ell \cong 4\ko_\ell$, and we shall see that this is indeed true. Then we
define a {\bf jumping line} as a line with $\ke|_\ell \ne 4\ko_\ell$. The aim
of this section is to understand to some extend the geometry of the subsets
of the Grassmanninan of lines on $\P^4$ corresponding to different nontrivial 
splittings of $\ke|_\ell$.

The theorems \ref{JumpDefect},\ref{ConicGrass},\ref{SurfGrass} and corollary
\ref{ConicDGrass} can be summarized in the 

\begin{theorem} \label{JumpTh}
Let an instanton bundle $\ke$ of rank 4 with quantum number 2 on $\P^4=\P V$ 
($V$ is a 5-dimensional vector space)
be given by a monad
$$ 
     2\Omega_{\P^4}^4(4) \stackrel{M}{\longrightarrow} 2\Omega_{\P^4}^1(1)
      \stackrel{N}{\longrightarrow} 2\ko_{\P^4} 
$$
and let the linear map
$$     k^2\otimes V \stackrel{\tilde M}{\longrightarrow} 
    k^2\otimes \bigwedge^4 V.
$$
be the contraction with the matrix $M$.
Let $\ell$ be a line in $\P^4$ identified via the Pl\"ucker embedding with
an element of $\P\bigwedge^2 V \, \supset \, \G (2,5)\, \ni\, \ell$. 
If $\langle \ell \rangle = \langle x \wedge y \rangle$ in $\P\bigwedge^2 V$ and
$M=(m_{i,j})_{i,j=1,2}$, denote by $M\wedge \ell$ a matrix 
$(m_{i,j}\wedge x \wedge y)_{i,j=1,2} \in Mat_{2\times 2}(\bigwedge^5 V) \cong
 Mat_{2\times 2}(k)$ (well defined up to  proportionality).
Then:\\
(i) if $\rank\, M\wedge\ell = 2$, then $\ell$ is not a jumping line; \\
(ii) if $\rank\, M\wedge\ell = 1$, then $\ell$ is a jumping line with splitting
$$ \ke_\ell \cong 
          \ko_\ell(-1) \oplus \ko_\ell^{\oplus 2} \oplus \ko_\ell(1); $$
(iii) if $ M\wedge\ell = 0$, then $\ell$ is a jumping line with
$$
\ke_\ell \cong \ko_\ell(-1)^{\oplus 2} \oplus \kf
          \text{\quad or\quad } 
\ke_\ell \cong\ko_\ell(-2) \oplus \ko_\ell \oplus \kf
$$
where $\kf \cong \ko_\ell(1)^{\oplus 2}$ \quad or 
          \quad $\kf\cong\ko_\ell \oplus \ko_\ell(2)$; \\
(iv) in the situation of (iii), the jumping lines $\ell$ with the property 
$\ke_\ell=\ko_\ell(-2)\oplus\ko_\ell\oplus\kf$
form a smooth conic on the Grassmannian $\G (2,5)$; \\
(v) in the situation of (iii), the jumping lines $\ell$ with the property 
$\kf = \ko_\ell\oplus\ko_\ell(2)$ form either a smooth conic, if 
$\rank\,\tilde M\,=8$, or, if $\rank\,\tilde M\,=7$, a surface 
 on the Grassmannian $\G (2,5)$.
\end{theorem}

Restricting a monad
$$ 2\Omega^4(4) \stackrel{M}{\longrightarrow} 
   2\Omega^1(1) \stackrel{N}{\longrightarrow} 2\ko $$
on $\P^4$ with the cohomology sheaf $\ke$
to a line $\ell=\P W  \subset \P^4 = \P V$, we get a monad of the restriction
$\ke|_\ell$:
$$
  2 \left[ \Omega^1_\ell(1) \otimes \bigwedge^3 (V/W)^\vee \right] 
    \stackrel{M}{\longrightarrow} 
   2 \left[ \Omega^1_\ell(1) \ \oplus \ (V/W)^\vee \otimes \ko(1)|_\ell \right] 
 \stackrel{N}{\longrightarrow} 2\ko 
$$
Clearly $\ell$ is a jumping line iff $H^1(\ke|_\ell(-1)) \ne 0$.\\
The monad on $\ell$ twisted down once decomposes into two exact sequences:
$$ 0 \longrightarrow 2\Omega^4(3)|_\ell \longrightarrow \kk(-1)|_\ell 
       \longrightarrow \ke(-1)|_\ell \longrightarrow 0 $$
$$ 0 \longrightarrow \kk(-1)|_\ell \longrightarrow 2\Omega^1|_\ell 
       \longrightarrow 2\ko|_\ell(-1) \longrightarrow 0 $$
whose long exact cohomology sequences give rise to the diagram:
$$
\xymatrix{
 & & 0 \ar[d] \\
 & {2H^1(\Omega^1_\ell \otimes \bigwedge^3 (V/W)^\vee)} 
   \ar[r] \ar[rd] \ar@{=}[ld]  &
   {H^1(\kk|_\ell(-1))} \ar[r] \ar[d]^{\cong} & {H^1(\ke|_\ell(-1))} \ar[r] & 0 \\
 {2\bigwedge^3(V/W)^\vee} \ar[rd]_{\lrcorner M} && 
   {2H^1(\Omega^1_\ell \ \oplus \ (V/W)^\vee \otimes \ko|_\ell(-1))} 
          \ar[d] \ar@{=}[ld] \\
 & 2k & 0
}
$$ 
therefore 
$$ \dim_k H^1\ke_\ell(1) \ =\ 
2 -\rank(2\bigwedge^3(V/W)^\vee \stackrel{M}{\to} 2k)
    \ = \ 2 - \rank(M\wedge\ell) . $$
To justify the last equality sign, we use the diagram 
$$\xymatrix{
& {2\bigwedge^3(V/W)^\vee \ar[r]^{\quad \lrcorner M} \ar[d]_{\cong}}
     & 2k \ar[d]^{\cong} \\
2k \ar@{}[r]|<<{\quad \cong} & {2\bigwedge^3(V/W)^\vee \otimes \bigwedge^2W^\vee} 
       \ar[r]_{\qquad \qquad \lrcorner (M\wedge\ell)}
 & 2k
}$$
which can be designed to be comutative by choosing the isomorphism in such a way
that the diagram 
$$\xymatrix{
& {\bigwedge^3(V/W)^\vee \ar[r]^{\quad \lrcorner w} \ar[d]_{\cong}}
     & k \ar[d]^{\cong} \\
k \ar@{}[r]|<<{\quad \cong} & {\bigwedge^3(V/W)^\vee \otimes \bigwedge^2W^\vee} 
       \ar[r]_{\qquad \qquad \lrcorner (w\wedge\ell)}  & k
}$$
commutes for contractions with any $w \in \bigwedge^3 V$

So, we have proved

\begin{theorem} \label{JumpDefect} 
Depending on $M \wedge \ell$, the following statements hold:
$$
\begin{array}{rcl}
\rank\, M\wedge\ell = 2 & \Rightarrow & \ke_\ell \cong \ko_\ell^{\oplus 4} \\ 
\rank\, M\wedge\ell = 1 & \Rightarrow & \ke_\ell \cong 
          \ko_\ell(-1) \oplus \ko_\ell^{\oplus 2} \oplus \ko_\ell(1) \\ 
\rank\, M\wedge\ell = 0 & \Rightarrow & \ke_\ell 
     \cong \ko_\ell(-1)^{\oplus 2} \oplus \kf
          \text{\quad or\quad } \ko_\ell(-2) \oplus \ko_\ell \oplus \kf
\end{array}
$$
where $\kf \cong \ko_\ell(1)^{\oplus 2}$ \quad or 
          \quad $\ko_\ell \oplus \ko_\ell(2)$.
\end{theorem}

Clearly $\ko_\ell(-2) \oplus \ko_\ell$ appears in the splitting iff 
$H^1\ke_\ell \ne 0$
The long exact cohomology sequences on $\ell$ yield
$$
\xymatrix{
& 0 \ar[d] \\
& 2H^0\left( \Omega^1_\ell(1) \ \oplus \ (V/W)^\vee \otimes \ko|_\ell \right) 
        \ar@{=}[r] \ar[d] 
         & 2(V/W)^\vee \ar[d]^{\lrcorner N} \\
& 2H^0(\ko_\ell) \ar@{=}[r] \ar[d] & 2k \\
0 \ar[r] & H^1\kk_\ell \ar[r]^{\cong} \ar[d] & H^1\ke_\ell \ar[r] & 0 \\
& 0
}$$
hence $ H^1\ke_\ell \cong \coker(\ 2(V/W)^\vee \stackrel{N}{\to} 2k \ ) $. 
Therefore $H^1\ke_\ell \ne 0$ iff $N: 2(V/W)^\vee \to 2k$ is not surjective iff
$N^\vee: 2k \to 2(V/W)$ is not injective. 

In the case $N=\left( \begin{array}{cc} e_1 & e_2 \\ e_3 & e_4 \end{array} \right)$
it means that the intersection of linear subspaces of $2V$
$$
 2W \ \cap \ 
\left( k \cdot {e_1 \choose e_2} + k \cdot {e_3 \choose e_4} \right) \ \ne \ 0. 
$$
i.e. there exists $(s:t)\in\P^1$ such that
$ \left\{ {{te_1 + se_3 \in W} \atop {te_2 + se_4 \in W}} \right.$ \\
So, to describe all jumping lines with $\ko_\ell(-2) \oplus \ko_\ell$ take 
the ruling lines of the scroll between the lines $\overline {e_1 e_3}$ and
 $\overline {e_2 e_4}$. This system of lines corresponds to a smooth conic in 
$\G(2,5) \subset \P^5$. 

If $N=\left( \begin{array}{cc} e_1 & e_2 \\ e_3 & e_1 \end{array} \right)$, 
we get the conditions 
$ \left\{ {{te_1 + se_3 \in W} \atop {te_2 + se_1 \in W}} \right.$ \\
and then the jumping lines fill the plane $\overline {e_1 e_2 e_3}$ and form 
again a smooth conic 
$$ 
(te_1+se_3)\wedge(te_2+se_1) \quad = \quad t^2e_{12} - tse_{23} - s^2e_{13}
$$
in the Grassmannian $\G(2,5)$. 

We formulate the result as

\begin{theorem} \label{ConicGrass}
The jumping lines $\ell$ with the property 
$\ke_\ell=\ko_\ell(-2)\oplus\ko_\ell\oplus\kf$
form a smooth conic on the Grassmannian $\G (2,5)$.
\end{theorem}

Since in the case $H^0\ke^\vee = 0$ the dual bundle $\ke^\vee$ is an 
instanton bundle again, we get

\begin{corollary} \label{ConicDGrass}
If $H^0\ke^\vee=0$, then the locus of jumping lines $\ell$ 
with $\ke_\ell=\kg\oplus\ko_\ell\oplus\ko_\ell(2)$ is also a smooth conic on
$\G (2,5)$.
\end{corollary}

If $h^0\ke^\vee=1$, then we have an extension
$$ 
0 \longrightarrow \ko \longrightarrow \ke^\vee \longrightarrow \ke_1^\vee
                                               \longrightarrow 0,
$$
where $\ke^\vee_1$ is a cohomology sheaf (not a bundle) of the monad
$$
 2\Omega^4(4) \stackrel{M^\vee}{\longrightarrow} 2\Omega^1(1)
      \stackrel{P^\vee}{\longrightarrow} 3\ko .
$$

If we restrict the above extension to a line $\ell$, we get 
$$\ko_\ell \longrightarrow \ke^\vee_\ell \longrightarrow 
\ke_{1,\ell}^\vee  \longrightarrow 0$$
where $\ke_{1,\ell}^\vee$ is reflexive hence locally free. 
One can choose an isomorphism of $\ke_\ell^\vee \cong \kl \oplus \ko(p)$ such 
that this sequence becomes
$$\xymatrix{
0 \ar[r] & \kl \ar[r] \ar@{}[d]|{\bigoplus} & \kf \ar[r] \ar@{}[d]|{\bigoplus} & 0\\
\ko_\ell \ar[r] & \ko_\ell(a) \ar[r] & \kt \ar[r] & 0
}$$
where $a\ge0$\,,\,$\ke_{1,\ell}^\vee \cong \kf \oplus \kt$, and $\kt$ is a torsion 
sheaf. 
Since $\kt$ has to be $0$, we get $a=0$ and again 
$\ke^\vee_\ell \cong \ko_\ell \oplus \ke_{1,\ell}^\vee$

We are interested in lines $\ell$ with the property that $\ko_\ell(-2)$ is contained in 
the splitting of $\ke_{1,\ell}^\vee$, i.e. $H^1(\ke_{1,\ell}^\vee) \ne 0$. Such $\ell$'s 
will be precisely those giving $\ko_\ell \oplus \ko_\ell(2)$ in the splitting 
of $\ke_\ell$. 

But $ H^1(\ke^\vee_{1,\ell}) \cong \coker(\ 2(V/W)^\vee \stackrel{P^\vee}{\to} 3k \ ) $. 
Therefore $H^1\ke_\ell \ne 0$ iff $N: 2(V/W)^\vee \to 2k$ is not surjective iff
$P: 3k \to 2(V/W)$ is not injective. \\
If $P=\left( \begin{array}{ccc} a_1 & a_2 & a_3 \\ b_1 & b_2 & b_3 \end{array} \right)$, 
it will mean that the intersection of linear subspaces of $2V$:
$$
2W \ \cap \ 
\left( k \cdot {a_1 \choose b_1} + k \cdot {a_2 \choose b_2} 
+ k \cdot {a_3 \choose b_3} \
                               \right) \ \ne \ 0. 
$$
i.e. there exists $(s:t:u)\in\P^2$ such that
$ \left\{ {{sa_1 + ta_2 + ua_3 \in W} \atop {sb_1 + tb_2 + ub_3 \in W}} \right.$ \\
This system of lines corresponds to a surface in 
$\G(2,5) \subset \P^5$ given by a parametrization  
$(sa_1 + ta_2 + ua_3 ) \wedge (sb_1 + tb_2 + ub_3)$.
Note that this expression is nowhere zero because of the nondegeneracy of $P$.
Also note that by theorem \ref{HarrisThm} after a suitable basis choice in $V$  \ 
$P$ can be written in one of the (possibly) two forms 
 $P\sim \left( \begin{array}{ccc} e_0 & e_1 & e_2 \\ e_1 & e_2 & e_3 \end{array} \right)$\
or\  $P\sim\left( \begin{array}{ccc} e_0 & e_1 & e_3 \\ e_1 & e_2 & e_4 \end{array} \right)$.

The argument presented proves the

\begin{theorem} \label{SurfGrass}
If $h^0\ke^\vee=1$ then the locus of jumping lines with
$\ke_\ell=\kg\oplus\ko_\ell\oplus\ko_\ell(2)$ is a surface in $\G (2,5)$.
\end{theorem}

\newpage

\section{Subsets of the Moduli Space}
\label{subsets}

In this section we are studying the following three special subsets of $\km$:\\
$\km^3\ :=\ $ sheaves in whose monads $\dim\,\Span\,N\,=\,3$  \\   
$\km^7\ :=\ $ sheaves $\ke$ in whose monads $rank\,\tilde M\, =\, 7$\, or,
equivalently (cf. corollary \ref{rkMtilde}), $H^0 \ke^\vee = 1$ \\
$\km^{sd}\ :=\ $ self-dual bundles $\ke$, i.e. bundles with
$\ke \cong \ke^\vee$.

\begin{theorem} \label{dimM3}
$\km^3$ is irreducible of codimension 2.
\end{theorem}

{\bf Proof} Denote :
$$\begin{array}{rcl}
Y^3 := \{ \langle N \rangle : \dim\Span\,N\,=3 \} \subset Y & ; &  
                 Y^3_0 := Y^3 \cap Y_0 , \\
X^3 := \{ (\langle M \rangle , \langle N \rangle ) \in X \ : 
              \ \langle N \rangle \in Y^3 \} \subset X & ; & 
          X^3_0 := X^3 \cap X_0 .
\end{array}$$  

Since $Y_3^0$ is dense in $Y^3$ and $X_0$ is a fibration over $Y_0$ with 
irreducible fibers of constant dimension 19, we get 
$$ \codim_X X^3 = \codim_Y Y^3 . $$
If $\langle N \rangle \in Y^3$, it means that its entries are linearly 
dependent. But the entries of $N$ are a priori arbitrary vectors from the 
5-dimensional vector space $V$, and it is classically known that the condition 
of the linear dependence of 4 vectors in $k^5$ has codimension 2.

Irreducibility can be verified by construcing the surjection from an 
irreducible set onto $Y_0^3 / G'$ as $\km$ is a fibration over it with 
irreducible fibers. The composite map
$$\begin{array}{ccccc}
V^0 \subset 3V & \longrightarrow & Y_3^0 & \longrightarrow & Y_3^0/G' \\
(v_1,v_2,v_3) & \longmapsto & 
   \left( \begin{array}{cc} v_1 & v_2 \\ v_3 & v_1 \end{array} \right)
   & \longmapsto &
  \left[ \left( \begin{array}{cc} v_1 & v_2 \\ v_3 & v_1 \end{array} 
          \right) \right]
\end{array}$$
where $V^0$ is the set of triples of linear independent vectors in $V$ and, 
clearly, an open subset of $3V$.
does the job. \qed


\begin{theorem} \label{M7thm}
$\km^7$ is irreducible of codimension 3 and contained in 
$\km^3$.
\end{theorem}

{\bf Proof.} Denote by $X^7_0 \subset X_0$ a set of all pairs 
$(\langle M \rangle,\langle N \rangle)$ with $H^0\ke^\vee=1$.
By the classification done above $X^7_0$ is an open nonempty subset in
$$X^7 = \{ (\langle M \rangle,\langle N \rangle) : 
           \left( \frac{N}{\lambda x \, \mu x} \right) \wedge M = 0 , 
           \exists x\in \P V, \, (\lambda,\mu)\in \P^1 \}
           \,\subset\, \P Mat_{2 \times 2}(\bigwedge^3 V) \times Y^3 $$
which, in turn, is an image of the set
$$\begin{array}{rcl}
\tilde X^7 & = & \{ (\langle M \rangle,\langle N \rangle,x,(\lambda,\mu)) :
        \left( \frac{N}{\lambda x \, \mu x} \right) \wedge M = 0 ,
           \dim\,\Span\, N\, =3,\,  
           \exists x\in \P \Span\,N \, , \, (\lambda,\mu)\in \P^1 \} \\
  & \subset & \P Mat_{2 \times 2}(\bigwedge^3 V) \times Y^3 \times 
                                                  \P V \times \P^1 . 
\end{array}$$
Now consider the morphism $\pi$ from the following diagram:
$$\xymatrix{
{\tilde X^7} \ar@{^{(}->}[r]^{\text{\scriptsize closed}\qquad\qquad\qquad} 
             \ar@{->>}[dr]_{\pi} &
    {\P Mat_{2 \times 2}(\bigwedge^3 V) \times Y^3 \times \P V \times \P^1 } 
       \ar[d]   \\
& { Y^3 \times \P V \times \P^1 }
}$$
Clearly $\pi$ is a projective morphism onto an irreducible subvariety
$$\{ \dim\,\Span\,N=3; \, x\in \Span\,N \} \subset
Y^3 \times \P V \times \P^1$$
with equidimensional irreducible fibers (which are, in fact, 
$(8\cdot 2-1)$-dimensional
linear subspaces of $\P Mat_{2 \times 2}(\bigwedge^3 V)$ \,), 
hence $\tilde X^7$ and so are $X^7, X^7_0$ and $\km^7$.

Further, $\km^7 \subset \km^3$, the latter has codimenion 2 and both are 
irreducible. Therefore, to prove the statement about the dimension it suffices 
to show that $\codim_{X_0}\, X^7_0 \, \le 3$.

The general form of the matrix $M$ as given in the introduction allows us to 
write the general form of $\tilde M$ as well. 
If $N=\left( \begin{array}{cc} e_1 & e_2 \\ e_3 & e_4
\end{array} \right)$, a computation shows that
$$ \tilde M = \left( \begin{array}{cccccccccc}
p_6 & p_5 & p_9 & p_8, 0 &  q6 & q_5 & q_9 & q_8 & 0 \\ 
-p_2 &-p_1 &0 &  0 &  p_8 & -q_2 &-q_1 &0 &  0 &  q_8 \\ 
p_3 & p_2 & 0 &  0 &  -p_9 &q_3 & q_2 & 0 &  0 &  -q_9 \\
0 &  0 &  p_2 & p_1 & p_5 & 0 &  0 &  q_2 & q_1 & q_5 \\
0 &  0 &  -p_3 &-p_2 &-p_6 &0 &  0 &  -q_3 &-q_2 &-q_6 \\
p_5 & p_4 & p_8 & p_7 & 0 &  q_5 & q_4 & q_8 & q_7 & 0 \\ 
-p_1 &-p_0 &0 &  0 &  p_7 & -q_1 &-q_0 &0 &  0 &  q_7 \\
p_2 & p_1 & 0 &  0 &  -p_8 &q_2 & q_1 & 0 &  0 &  -q_8 \\
0 &  0 &  p_1 & p_0 & p_4 & 0 &  0 &  q_1 & q_0 & q_4 \\
0 &  0 &  -p_2 &-p_1 &-p_5 &0 &  0 &  -q_2 &-q_1 &-q_5 
\end{array} \right)
$$
and in the case $N=\left( \begin{array}{cc} e_1 & e_2 \\ e_3 & e_1
\end{array} \right)$  
$$ \tilde M = \left( \begin{array}{cccccccccc}
p_5 & p_7 & p_9 & p_8 & 0 &  q_5 & q_7 & q_9 & q_8 & 0 \\ 
-p_2 &0 &  0 &  0 &  p_8 & -q_2 &0 &  0 &  0 &  q_8 \\
p_3 & 0 &  0 &  0 &  -p_9 &q_3 & 0 &  0 &  0 &  -q_9 \\
-p_1 &0 &  0 &  0 &  p_7 & -q_1 &0 &  0 &  0 &  q_7 \\
0 &  -p_1 &-p_3 &-p_2 &-p_5 &0 &  -q_1 &-q_3 &-q_2 &-q_5 \\
p_4 & p_6 & p_8 & -p_7 &0 &  q_4 & q_6 & q_8 & -q_7 &0 \\
p_1 & 0 &  0 &  0 &  -p_7 &q_1 & 0 &  0 &  0 &  -q_7 \\
p_2 & 0 &  0 &  0 &  -p_8 &q_2 & 0 &  0 &  0 &  -q_8 \\
-p_0 &0 &  0 &  0 &  p_6 & -q_0 &0 &  0 &  0 &  q_6 \\
0 &  -p_0 &-p_2 &p_1 & -p_4 &0 &  -q_0 &-q_2 &q_1 & -q_4 
\end{array} \right)
$$ 
We see that in the $10\times 10$-matrices $\tilde M$ there are two pairs of 
equal rows. If we delete the redundant rows, we end up with $8\times 10$ 
matrices. The condition that $\rank\,\tilde M\, =7$ can be then locally 
written by three equations as for all $M$'s we consider 
$\rank\,\tilde M\, \ge 7$. This implies that $\codim_{X_0}\,X^7 \le 3$, Q.E.D.
\qed


\begin{theorem} The subset $\km^{sd}$ of the self-dual bundles is an
irreducible closed subvariety of codimension 4
\end{theorem}

{\bf Proof.}
Denote by $X_0^{sd}$ the set of the parameter space sorresponding to the 
self-dual bundles. We shall prove that $X_0^{sd}$ is a smooth irreducible 
closed $G$-invariant subvariety of codimension 4, and by the axioms of a 
geometric quotient the same is then true for $\km^{sd}$. First we give a more 
handier criterion of self-duality.

\begin{lemma} A monad
$$ 2\Omega^4(4) \stackrel{M}{\longrightarrow}  2\Omega^1(1) 
                 \stackrel{N}{\longrightarrow}  2\ko  $$
defines a self-dual bundle $\ke$ if and only if $\dim\,\Span\, M\, = 3$.
In this case, there exist a matrix $M_0=M_0^\vee$ and a group element 
$g\in GL(2)$ such that $M=M_0\cdot g$.
\end{lemma}

{\bf Proof of the lemma.}  
We know from the chapter on the dual bundles that $\ke\cong\ke^\vee$ iff
$M \sim M^\vee$. Obviously, if $M$ is a symmetric matrix, then $\dim\Span\,M=3$.
Conversely, if $\dim\Span\,M=3$ and there is no generalized column of the form
$\lambda\otimes x$, then by lemma \ref{span3} there exist $M_0=M_0^\vee$ and 
$g\in GL(2)$ such that $M=M_0\cdot g$, in particular, 
$\dim\Span\,M\,=\dim\Span\,M_0\,=3$, hence the lemma. \qed

We can proove that $\codim_{X_0}\, X_0^{sd} \le 4$ simply by looking at the 
equations. For any $\langle N \rangle \in Y_0$ 
consider the fiber $X_{0,N}$ of $X$ over $\langle N \rangle$. Then $X_{0,N}$ is isomorphic to an open subset of $P^{19}$ with homogeneous coordinates 
$p_0,\dots,q_9$. The condition $\dim\Span\,M=3$ can be written as
$$
\rank\,\left( \begin{array}{ccccccc}
p_1 & p_2 & p_3 & p_5 & p_6 & p_8 & p_9 \\
p_0 & p_1 & p_2 & p_4 & p_5 & p_7 & p_8 \\
q_1 & q_2 & q_3 & q_5 & q_6 & q_8 & q_9 \\
q_0 & q_1 & q_2 & q_4 & q_5 & q_7 & q_8 \\
\end{array} \right)\ = \ 3
$$
for the case $\dim\Span\,N=4$ and
$$
\rank\,\left( \begin{array}{ccccccc}
p_1 & p_2 & p_3 & p_5 & p_7 & p_8 & p_9 \\
p_0 & p_1 & p_2 & p_4 & p_6 & p_7 & p_8 \\
q_1 & q_2 & q_3 & q_5 & q_7 & q_8 & q_9 \\
q_0 & q_1 & q_2 & q_4 & q_6 & q_7 & q_8 
\end{array} \right)\ = \ 3
$$
if $\dim\Span\,N=3$. As these ranks are always $\ge 3$, the codimension of the
subvariety $X^{sd}_{0,N} \subset X_{0,N}$ defined by these equations is $\le 4$.

We shall construct a surjection of an irreducible $38-4=34$-dimensional
variety onto $X_0^{sd}$ and thus prove the inverse inequality
$\codim_{X_0}\, X_0^{sd} \ge 4$ and the irreducibility of $X_0^{sd}$.

Consider the variety $Z \in X_0$ given by the condition 
$\langle M \rangle = \langle M^\vee \rangle$. In each $X_{0,N}$ it is given by
7 linear equations hence each $Z \cap X_{0,N}$ a $19-7=12$-dimensional linear 
projective space and by proposition \ref{IrredCrit} $Z$ is an irreducible 
$38-7=31$-dimensional variety.  We define a map
$$ Z \times SL(2) \to X_0^{sd} \ : \ 
  (\, (\langle M \rangle,\langle N \rangle), g\, ) \mapsto 
      (\langle M\cdot g \rangle,\langle N \rangle) $$
which is surjective by the lemma above and where $Z \times SL(2)$ is 
$31+3=34$-dimensional and irreducible.
This completes the proof. \qed

The same argument shows that

\begin{corollary} The intersection $\km^{sd} \cap \km^3$ is an irreducible 
subvariety of codimension 4 in $\km^3$.
\end{corollary}

\begin{remark} From the section on the dual bundle we see, 
that $\km^7 \cap \km^{sd} = \emptyset$.
\end{remark}

\newpage

\section{A Description of the Moduli Space as a Fibration}
\label{fibration} 

From the proofs done in the sections \ref{formal} and \ref{etale} we 
can derive also the following

\begin{theorem} Consider the group $G'=SL(2) \times SL(2)$ operating on $Y_0$. 
Then:\\
(i) There exists a geometric quotient $Y_0/G'=:Q$ which is a smooth 
quasi-projective variety of dimension 13; \\
(ii) The quotient map $Y_0 \to Q$ is a Mumford principal fibration;\\
(iii) $Y^3_0$ is a $G'$-invariant closed subset in $Y_0$ and therefore defines 
a closed subset $Q^3 \subset Q$ of codimension 2;\\
(iv) The projection $X_0 \to Y_0$  induces the surjective map $\km \to Q$ 
whose geometric fibres are isomorphic to open subsets in the Grassmannian   
$\G (2,10)$.
\end{theorem} 

{\bf Proof\ } of (i) and (ii) was in fact done in the step 1 of the proof
of proposition \ref{IsotrivProp} together with the proof of theorem 
\ref{Msmooth}. Part (iii) follows from the proof of \ref{dimM3}. Part (iv)
is obvious from definitions of $X_0$, $Y_0$ and of group actions. \qed

\begin{prop}
There exists a generically 2:1 map $\delta : Q \to Z$ where $Z=\P S^2V$ 
is the space of quadratic hypersurfaces in $\P V^\vee$ and $\delta$ has 
the properties:
(i) the ramification locus of $delta$ is precisely $Q^3$
(ii) the image $\delta (Q)$ is contained in the set of sigular quadric 
hypersurfaces, i.e. in the set of quadratic cones . 
\end{prop}

{\bf Sketch of the proof} (cf. \cite{NarTrm} for more details). 
Define the morphism
$$\begin{array}{crcl} \det : & Y_0 & \longrightarrow & Z=\P S^2V \\
  & \langle N \rangle & \longmapsto & \langle \det (N) \rangle
\end{array}$$
This morphism is $G'$-equivariant (in the situation of reduced 
quasi-projective varieties it suffices to check this set-theoretically on 
closed points) hence factors through a map $\delta : Q \to Z$.

Calculations of the primages show that 
$$\delta^{-1}(e_1e_4-e_2e_3) = \left\{ 
\left[ \left( \begin{array}{cc} e_1 & e_2 \\ e_3 & e_4 \end{array} 
\right) \right] \ , \
\left[ \left( \begin{array}{cc} e_1 & e_3 \\ e_2 & e_4 \end{array} 
\right) \right] 
\right\} $$
$$ \delta^{-1}(e_1^2-e_2e_3) = \left\{ \quad
\left[ \left( \begin{array}{cc} e_1 & e_2 \\ e_3 & e_1 \end{array} 
\right) \right] \quad \right\}$$ 
where $[N]$ denotes the class of $\langle N \rangle$ modulo the action of $G'$.
This implies the property (i), while the property (ii) is obvious.
\qed

There is also another observation that helps to understand the structure of
$Q$ and hence that of $\km$. The construction we are presenting now will be 
implicitly used in the section \ref{plane}. 

\begin{construction} of a rational map $\sigma : Q \to \P V^\vee$.
\end{construction}

Similarly to the last proposition, induce $\sigma$ by defining a 
$G'$-equivariant morphism 
$$\begin{array}{crcl} \Span : & Y_0^4 & \longrightarrow & \P V^\vee \\
  & \langle N \rangle & \longmapsto & \Span\, N
\end{array}$$
where $Y_0^4 = Y_0 - Y_0^3$. Then $\sigma$ is singular along $Q^3$. 

However, we can construct a birational modification $\pi : \tilde Q \to Q$
such that $\sigma$ lifts to $Q$ as a regular map.
Define $\tilde Q$ as the graph of an incidence correspondence
$$\tilde Q := \{ (\, [ N ] , H \, ) : \Span\,N \subseteq H \} 
              \subset Q \times \P V^\vee $$
and $\pi$ as the natural projection.

Clearly, $\pi$ is birational and an isomorphism outside $\pi^{-1}Q^3$. 
The morphism $\tilde \sigma$ appears from the commutative diagram:
$$\xymatrix{
{\tilde Q} \ar@{}[r]|{\subset \quad} \ar[d]_{\pi} 
           \ar@/^1pc/[rr]^{\tilde \sigma} & 
          Q \times \P V^\vee  \ar[r]_{\ pr_2} \ar[dl]^{pr_1} &
          \P V^\vee \ar@{=}[d] \\
Q \ar@{-->}[rr]_{\sigma} & & \P V^\vee
}$$
It is also easy to see that fibers of $\pi$ over $Q^3$ are isomorphic to 
$\P^1$.

\newpage

\section{Restriction to a 3-Plane }
\label{plane}

\begin{prop} Let $\ke$ be the instanton bundle defined by a monad
$$ 0 \to 2\Omega^4(4) \stackrel{M}{\longrightarrow} 2\Omega^1(1) 
    \stackrel{N}{\longrightarrow} 2\ko \to 0  $$
and $H=\P W \subset \P V$ a projective 3-plane. Then $h^0(\ke |_H)=2$ iff
$W \supset \Span\, N$.
\end{prop}

{\bf Proof.} Restricting the monad to $H$ and denotinh by $\kn$ its 
kernel bundle, we get:
$$0 \to 2\Omega^3_H(3) \otimes (V/W)^\vee 
   \stackrel{M(H)}{\longrightarrow} 
   2 [\, \Omega^1_H(1) \, \oplus \, (V/W)^\vee\otimes\ko_H \, ]
    \stackrel{N(H)}{\longrightarrow} 2\ko_H \to 0 
$$
$$
 H^0(\ke|_H) = H^0(\kn|_H) = \ker\, 
       (\, 2(V/W)^\vee \stackrel{\lrcorner N}{\longrightarrow} 2k \, ) 
$$    
But $\dim\,H^0(\ke|_H) = 2$ means precisely that 
$2(V/W)^\vee  \stackrel{\lrcorner N}{\longrightarrow} 2k$ is the zero map, or, 
equivalently, that the dula map $2k \stackrel{N^\vee}{\longrightarrow} 2(V/W)$ 
is zero, i.e. $\Span\, N \subset W$. Hence the proposition., \qed

\begin{theorem}  \label{PlaneTh} 
If $\dim\,\Span\, N=4$, denote $W = \Span\, N$ and 
$H=\P W \subset \P V$ will be a projective 3-plane.
Then $\ke|_H \cong \kf \oplus 2\ko_H$. 
\end{theorem}

{\bf Proof.} Put $N=\left(\begin{array}{cc} e_1 & e_2 \\ e_3 & e_4
\end{array} \right)$ and write the matrix $M$ as $M=M'\wedge e_0 + M''$ where
\linebreak
$M'\in Mat_{2\times 2} (\bigwedge^2 W)$ and 
$M''\in Mat_{2\times 2} (\bigwedge^3 W)$. 

The restricted monad gives us
$$\xymatrix{
0 \ar[r] & 2\Omega_H^3(3) \ar@{=}[ld] \ar[rr]^{{M' \choose M''}} \ar[rd] && 
     2\Omega_H^1(1) \oplus 2\ko_H \ar[rr]^{(N\ 0)} && 2\ko_H \ar[r] & 0 \\
2\ko_H(-1) && \kn'\oplus 2\ko_H \ar[ur] \ar[dr] \\
& 0 \ar[ru] & & \ke_H \ar[r] & 0
}$$
in particular, $\kn'=\ker\,(N : 2\Omega^1_H(1)\to 2\ko)$,
and hence we have a diagram
$$\xymatrix{
&& 0 \ar[d] \\
&& 2\ko_H \ar[d] \ar[d] \ar@{=}[r] & 2\ko_H \ar[d] \\
0 \ar[r] & 2\ko_H(-1) \ar@{=}[d] \ar[r] & \kn'\oplus 2\ko_H \ar[r] \ar[d]
   & \ke_H \ar[r] \ar[d] & 0 \\
& 2\ko_H(-1) \ar[r] & \kn' \ar[r] \ar[d]
   & \kf \ar[r] \ar[d] & 0 . \\
&& 0 & 0
}$$

{\bf Claim.} $M'$ defines an injection of sheaves 
$2\Omega_H^3(3) \to 2\Omega^1(1)$.

Indeed, if $M'$ is not injective as a sheaf homomorphism, we may assume 
$M'=\left(\begin{array}{cc} 0 & b' \\ 0 & d' \end{array} \right)$. Then 
the first column of $M$ has two elements $a,c\in \bigwedge^3H$. But if so,
$a^*,c^*$ are divisible by $e_0$ and hence  $a^{*2}=c^{*2}=0$ and by lemma
\ref{WedgeLemma} $M$ doesn't define a subbundle, contradiction.

We are looking for a splitting $\rho$ as shown in the diagram
$$\xymatrix{
&& 0 \ar[d] & 0 \ar[d] \\
&& 2\ko_H \ar[d]_{{0 \choose Id}=j} \ar@{=}[r] & 2\ko_H \ar[d]_{\varepsilon} \\
0 \ar[r] & 2\ko_H(-1) \ar@{=}[d] \ar[r]_{\text{\tiny ${M'\choose M''}$}} 
   & \kn'\oplus 2\ko_H \ar[r]_{\quad\pi} \ar[d]_{(\Id \ 0)} 
                                          \ar@/_/@{.>}[u]_{\sigma}`
   & \ke_H \ar[r] \ar[d] \ar@/_/@{.>}[u]_{\rho} & 0 \\
0 \ar[r] & 2\ko_H(-1) \ar[r] & \kn' \ar[r] \ar[d]
   & \kf \ar[r] \ar[d] & 0 \\
&& 0 & 0
} $$
If $\rho$ exists, then we also have a splitting $\sigma=\rho\pi$ as 
$ \sigma j = \rho\pi j = \rho \varepsilon = \Id_{2\ko_H}$.

Denote $\sigma=(S' \ S'')$. Since $\sigma$ is a splitting, we get 
$$ \Id = (S' \ S'')\cdot {0 \choose \Id} = S'' $$
and  $\sigma=(S' \ \Id )$.

Given any such splitting $\sigma$, it induces a splitting $\rho$ if and only if
$(S' \ \Id )\cdot {M' \choose M''}=0$. Therefore the question is reduced 
to the following:\\
When does a sheaf homomorphism $S' : \kn' \to 2\ko_H$ exist such that 
$S'M'+M''=0$?

Since we have an exact sequence 
$$ 
 0 \to \kn'_H \longrightarrow 2\Omega^1(1) \stackrel{N}{\longrightarrow} 
        2\ko_H \to 0
$$
and $\ext^1(2\ko_H,2\ko_H)=0$, every $S'$ lifts to some 
$S_1: 2\Omega_H^1(1) \to 2\ko_H$.

If we take two different liftings $S_{1,2}: 2\Omega_H^1(1) \to 2\ko_H$ of $S'$
then $S_1-S_2=R\cdot N$ for some \linebreak $R:2\ko_H\to 2\ko_H$, then 
$S_2\wedge M' = S_1 \wedge M' +  R \cdot N \wedge M'' = S_1 \wedge M'$.

We see now that the existence of $\rho$ is equivalent to existence of a matrix 
$S_1\in Mat_{2\times 2}(W)$ such that $S \wedge M' + M''$.

Let 
$$ S_1=\left( \begin{array}{cc} s_{11} & s_{12} \\ s_{21} & s_{22} 
\end{array} \right) , \
M'=\left( \begin{array}{cc} m'_{11} & m'_{12} \\ m'_{21} & m'_{22} 
\end{array} \right), \
M''=\left( \begin{array}{cc} m''_{11} & m''_{12} \\ m''_{21} & m''_{22} 
\end{array} \right); $$
then the conditions $S_1\wedge M'+M''=0$ can be rewritten as
$$\left. \begin{array}{l}
s_{11}\wedge m'_{11} + s_{12}\wedge m'_{21} = m''_{11}\\
s_{11}\wedge m'_{12} + s_{12}\wedge m'_{22} = m''_{12}\\
s_{21}\wedge m'_{11} + s_{22}\wedge m'_{21} = m''_{21}\\
s_{21}\wedge m'_{12} + s_{22}\wedge m'_{22} = m''_{22} \end{array} \right.
$$
Put $s_{11}=\sum_{i=1}^{4} x_i e_i$ and 
$s_{12}=\sum_{i=1}^{4} y_i e_i$, 
$$m''_{11}=\mu_1 e_{123} + \mu_2 e_{124} + \mu_3 e_{134} + \mu_4 e_{234} $$
$$m''_{12}=\nu_1 e_{123} + \nu_2 e_{124} + \nu_3 e_{134} + \nu_4 e_{234} $$
   
A computation shows that the first two conditions give in fact a non-homogeneous
system of linear equations 
$$\left(\begin{array}{cccc|cccc}  
-p_2 & p_3 & 0 & 0 & -p_1 & p_2 & 0 & 0 \\
-p_1 & p_2 & 0 & 0 & -p_0 & p_1 & 0 & 0 \\
0 & 0 & -p_2 & p_3 & 0 & 0 & -p_1 & p_2 \\
0 & 0 & -p_1 & p_2 & 0 & 0 & -p_0 & p_1 \\
\hline
-q_2 & q_3 & 0 & 0 & -q_1 & q_2 & 0 & 0 \\
-q_1 & q_2 & 0 & 0 & -q_0 & q_1 & 0 & 0 \\
0 & 0 & -q_2 & q_3 & 0 & 0 & -q_1 & q_2 \\
0 & 0 & -q_1 & q_2 & 0 & 0 & -q_0 & q_1 \\
\end{array} \right) 
\cdot \left( \begin{array}{c}
x_1 \\ x_2 \\ x_3 \\ x_4 \\ y_1 \\ y_2 \\ y_3 \\ y_4 
\end{array} \right) = \left( \begin{array}{c} 
\mu_1 \\ \mu_2 \\ \mu_3 \\ \mu_4 \\ \nu_1 \\ \nu_2 \\ \nu_3 \\ \nu_4 
\end{array} \right)
$$
with the matrix nondegenerated on an open dense subset of $X_0 - X^3_0$. 
Hence the set of monads with the property that $\ke|_{\Span\,N}$ splits
form an open dense subset of $X_0 - X^3_0$.

Now take $(X_0 - X^3_0)\times \P^3$ and consider the diagram of flat families
over $X_0 - X^3_0$, where $p$ is the product projection on $\P^3$.
$$\xymatrix{
0 \ar[r] & 2p^*\Omega_H^3(3) \ar@{=}[ld] \ar[rr] \ar[rd] && 
     2p^*\Omega^1(1) \oplus 2\ko \ar[rr] && 2\ko \ar[r] & 0 \\
2p^*\ko(-1) && \bkn \oplus 2\ko \ar[ur] \ar[dr] \\
& 0 \ar[ru] & & \bke \ar[r] & 0
}$$
which gives the diagram from the very beginning of the proof over each 
$(\langle M \rangle, \langle N \rangle)$ with $W=\Span\,N$. By what we have 
proven over each closed point, the composition $2\ko \to \bkn\oplus 2\ko 
\to \bke$ is an injection
and we may consider a short exact sequence
$$
0 \to 2\ko \to \bke \to \bkf \to 0
$$
where $\bkf$ is a cokernel sheaf. Since the families are flat over 
$X_0 - X^3_0$, this sequence splits over its closed subset. 

But $X_0 - X^3_0$ is irreducible, hence, putting things together, we see, 
that a splitting exists for each pair 
$(\langle M \rangle, \langle N \rangle) \in (X_0 - X^3_0)$,
Q.E.D. \qed

\begin{remark} As you can see from the proof, $\bkf$ is a family of instanton 
bundles on $\P^3$  with $c_2=2$. Such bundles have been studied, e.g., 
in \cite{NarTrm}.
\end{remark}

\newpage

\section{Sections of the Kernel Bundle }
\label{scroll}

In this section we describe the sections of the twisted kernel bundle $\kn(1)$ 
of the monad
$$
0\to 2\Omega^4(4) \stackrel{M}{\longrightarrow} 2\Omega^1(1)
  \stackrel{N}{\longrightarrow} \ko \to 0
$$ 
which vanish on sume 2-plane in $\P^4$.

We know that $h^0\kn(1)=10$ and 
$H^0(\kn(1)) \subset 2\bigwedge^2 V^\vee \cong 2\bigwedge^3 V$

It is easy to see that if $\xi \in H^0\Omega^1(2) = 2\bigwedge^3 V$ then the 
zero scheem of this section is precisely $\P V_\xi$ with the reduced structure,
where $V_\xi=\{ v\in V : v\wedge \xi = 0 \}$ is the space of linear factors of 
the 3-form $\xi$
Therefore a section $(\xi,\eta) \in H^0(\kn(1)) \cong 2\bigwedge^3 V$
vanishes if and only if 
$$\xi = \lambda \eta,\, \exists\lambda\in k-\{ 0\}; 
       \quad \xi^{*2}=\eta^{*2}=0 $$
because, as we have seen in the proof of lemma \ref{WedgeLemma}, 
$\xi^{*2}=0$ is 
eqiuvalent to $\dim\,V_\xi=3$.

If $\Gamma$ is the total matrix of syzygies of $N$, then 
$${\xi \choose \eta} =\Gamma \cdot (p_0,\dots,p_9)^T. $$

If $N=\left(\begin{array}{cc} e_1 & e_2 \\ e_3 & e_4
\end{array} \right)$
then 
the condition $\xi = \lambda \eta$ can be written as
$$\rank\,\left( \begin{array}{ccccccc}
p_1 & p_2 & p_3 & p_5 & p_6 & p_8 & p_9 \\
p_0 & p_1 & p_2 & p_4 & p_5 & p_7 & p_8 
\end{array} \right) = 1
$$
and 
$$\xi^{*2}=2 \cdot [\, (p_2 p_8 - p_1 p_9)e_{0123}
 + (p_3 p_8 - p_2 p_9)e_{0124} + (p_1 p_6 - p_2 p_5)e_{0134} 
 + (p_2 p_6 - p_3 p_5)e_{0234} + (p_2^2 - p_1 p_3)e_{1234} \, ] .  ]
$$

If $N=\left(\begin{array}{cc} e_1 & e_2 \\ e_3 & e_1
\end{array} \right)$, the condition $\xi = \lambda \eta$ becomes
$$\rank\,\left( \begin{array}{ccccccc}
p_1 & p_2 & p_3 & p_5 & p_7 & p_8 & p_9 \\
p_0 & p_1 & p_2 & p_4 & p_6 & p_7 & p_8 
\end{array} \right) = 1
$$
and $$\xi^{*2}=2 \cdot [\,
(p_3 p_8 - p_2 p_9)e_{0124}+ (p_1 p_8-p_2 p_7) e_{0134} + 
(p_1 p_9 - p_3 p_7)e_{0234} \,] . $$

In both cases we see that, if the condition on the rank is fulfilled, 
$\xi^{*2}=\eta^{*2}=0$ automatically.
The conditions on the rank define a variety which is classically known as a 
3-dimensional rational normal scrolls in $\P^9 = \P  H^0(\kn(1))$. 
See \cite[8.26]{harris} for other equivalent definitions and geometry 
of rational normal scrolls.

We have proved 

\begin{theorem} The subset of $\P  H^0(\kn(1))$ corresponding to sections 
vanishing on a 2-plane, is a 3-dimensional rational normal scroll.
\end{theorem}
\qed

\newpage
\addcontentsline{toc}{section}{References}

\vspace{8cm}
Hiermit ekl\"are ich, dass ich die vorliegende Arbeit selbstst\"andig 
erstellt und keine anderen als die angegebenen Hilfsmittel verwendet habe.

\end{document}